\documentclass[11pt]{article}
\usepackage{amsmath,amsthm,amssymb,color,epic,eepic,cite}
\newlength{\minitwocolumn}
\newcommand{\Z}{{\Bbb Z}} %??
\newcommand{\C}{{\Bbb C}} %??
\newcommand{\F}{{\cal F}}

\newcommand{\cA}{{\cal A}}

\renewcommand{\H}{{\cal H}}
\newcommand{\la}{\lambda}

\newcommand{\nn}{{\nonumber}}
\newcommand{\bea}{\begin{eqnarray}}
\newcommand{\ena}{\end{eqnarray}}
\newcommand{\beit}{\begin{itemize}}
\newcommand{\enit}{\end{itemize}}

\newcommand{\be}{\begin{eqnarray*}}
\newcommand{\en}{\end{eqnarray*}}
\newcommand{\lb}[1]{\label{#1}}
\newcommand{\id}{\hbox{id}}

\newcommand{\BW}[5]
{\left(\begin{array}{cc}#1 & #2 \cr #3 & #4 \cr\end{array}
\Biggl| #5\right)}

\def\infq4p#1{{(#1;q^4,p)_\infty}}

\newcommand{\EXP}[1]{{\exp\biggl\{#1\biggr\}}}

\newcommand{\al}{\alpha}

\newcommand{\vep}{\varepsilon}

\font\teneufm=eufm10
\font\seveneufm=eufm7
\font\fiveeufm=eufm5
\newfam\eufmfam
\textfont\eufmfam=\teneufm
\scriptfont\eufmfam=\seveneufm
\scriptscriptfont\eufmfam=\fiveeufm

\let\goth\frak
\newcommand{\slth}{\widehat{\goth{sl}}_2}

\newcommand{\uq}{U_q\bigl(\slth\bigr)}

\makeatletter
\@addtoreset{equation}{section}
\makeatother

\newtheorem{thm}{Theorem}[section]
\newtheorem{prop}[thm]{Proposition}

\newtheorem{conj}[thm]{Conjecture}

\newtheorem{df}{Definition}[section]

\newcommand{\nc}{\newcommand}
\nc{\mref}[1]{(\ref{#1})}
\nc{\vt}{v_{2\gL_0}}
\nc{\vo}{v_{\gL_0}}
\nc{\vot}{v_{\gL_1+\gL_0}}
\nc{\vw}{v_{\gL_1}}
\nc{\ppmm}{\genfrac{}{}{-10pt}{10pt}{++}{--}}
\nc{\wom}[5]{\Omega\left(\left.\begin{array}{ll}{#1}&{#2}\\{#3}&{#4}\end{array}\right|{#5}\right)}
\nc{\com}[5]{\chi\left(\left.\begin{array}{ll}{#1}&{#2}\\{#3}&{#4}\end{array}\right|{#5}\right)}
\nc{\we}[5]{W\left(\left.\begin{array}{ll}{#1}&{#2}\\{#3}&{#4}\end{array}\right|{#5}\right)}
\nc{\lmat}[6]{L_{#6}\left(\left.\begin{array}{ll}{#1}&{#2}\\{#3}&{#4}\end{array}\right|{#5}\right)}
\nc{\lmats}[5]{L\left(\left.\begin{array}{ll}{#1}&{#2}\\{#3}&{#4}\end{array}\right|{#5}\right)}
\nc{\lmatk}[6]{L^{(k)}_{#6}\left(\left.\begin{array}{ll}{#1}&{#2}\\{#3}&{#4}\end{array}\right|{#5}\right)}
\nc{\lmatsk}[5]{L^{(k)}\left(\left.\begin{array}{ll}{#1}&{#2}\\{#3}&{#4}\end{array}\right|{#5}\right)}
\nc{\hmat}[6]{h_{#6}\left(\left.\begin{array}{ll}{#1}&{#2}\\{#3}&{#4}\end{array}\right|{#5}\right)}
\nc{\hmats}[5]{H\left(\left.\begin{array}{ll}{#1}&{#2}\\{#3}&{#4}\end{array}\right|{#5}\right)}
\nc{\web}[5]{\overline{W}\left(\left.\begin{array}{ll}{#1}&{#2}\\{#3}&{#4}\end{array}\right|{#5}\right)}
\nc{\wep}[5]{W'\left(\left.\begin{array}{ll}{#1}&{#2}\\{#3}&{#4}\end{array}\right|{#5}\right)}
\nc{\wes}[5]{W^*\left(\left.\begin{array}{ll}{#1}&{#2}\\{#3}&{#4}\end{array}\right|{#5}\right)}
\nc{\wess}[5]{W^{**}\left(\left.\begin{array}{ll}{#1}&{#2}\\{#3}&{#4}\end{array}\right|{#5}\right)}
\nc{\cet}[7]{C^{#6}_{#7}\left(\left.\begin{array}{ll}{#1}&{#2}\\{#3}&{#4}\end{array}\right|{#5}\right)}
\nc{\bcet}[7]{\bar{C}^{#6}_{#7}\left(\left.\begin{array}{ll}{#1}&{#2}\\{#3}&{#4}\end{array}\right|{#5}\right)}
\nc{\wet}[7]{W^{#6}_{#7}\left(\left.\begin{array}{ll}{#1}&{#2}\\{#3}&{#4}\end{array}\right|{#5}\right)}
\nc{\bwet}[7]{\overline{W}^{#6}_{#7}\left(\left.\begin{array}{ll}{#1}&{#2}\\{#3}&{#4}\end{array}\right|{#5}\right)}
\nc{\wec}[7]{\widetilde{W}^{#6}_{#7}\left(\left.\begin{array}{ll}{#1}&{#2}\\{#3}&{#4}\end{array}\right|{#5}\right)}
\nc{\wgen}[6]{W^{#6}\left(\left.\begin{array}{ll}{#1}&{#2}\\{#3}&{#4}\end{array}\right|{#5}\right)}
\nc{\wgenp}[6]{W^{*{#6}}\left(\left.\begin{array}{ll}{#1}&{#2}\\{#3}&{#4}\end{array}\right|{#5}\right)}
\nc{\wo}[5]{\Omega\left(\left.\begin{array}{ll}{#1}&{#2}\\{#3}&{#4}\end{array}\right|{#5}\right)}
\nc{\wsgen}[8]{{#8}^{#6}_{#7}\left(\left.\begin{array}{ll}{#1}&{#2}\\{#3}&{#4}\end{array}\right|{#5}\right)}
\nc{\qbinom}[2]{{\genfrac{[}{]}{0pt}{}{{#1}}{{#2}}}_{q}}
\nc{\hhg}[4]{\phi\left({{{#1}\,\,\,{#2}}\atop{{#3}}};
                     {#4}\right)}
\nc{\fullhhg}[5]{{_1}\phi_2\left({{{#1}\,\,\,{#2}}\atop{{#3}}};
                     {#4},{#5}\right)}
\nc{\qp}[2]{({#1}\, ; \, {#2})_{\infty}}
\nc{\qpf}[1]{({#1}\, ; \, q^4)_{\infty}}
\nc{\pp}[1]{({#1}\, ; \, p)_{\infty}}
\nc{\qpp}[1]{({#1}\, ; \, p, q^4)_{\infty}}
\nc{\sect}{\section}
\nc{\ssect}{\subsection}
\nc{\sssect}{\subsubsection}
\nc{\ud}[1]{\underline{{#1}}}
\nc{\isomo}{\buildrel {\sim} \over \longrightarrow}
\nc{\Aff}{\operatorname{Aff}}
\nc{\ot}{\otimes}
\nc{\er}{\end{array}}
\nc{\bev}[1]{\begin{equation}\begin{array}{#1}}
\nc{\eeq}{\end{equation}}
\nc{\ee}{\end{eqnarray}}
\nc{\ben}{\begin{eqnarray*}}
\nc{\een}{\end{eqnarray*}}
\nc{\bec}{\begin{equation}\begin{array}{lll}}
\nc{\eec}{\end{array}\end{equation}}
\nc{\ed}{\end{document}}
\nc{\half}{\ensuremath{\frac{1}{2}}}
\nc{\Hom}{\operatorname{Hom}}
\nc{\vac}{|\textrm{vac}\rangle}
\nc{\dvac}{\langle\textrm{vac}|}
\nc{\ra}{\rightarrow}  
\nc{\lra}{\longrightarrow}
\nc{\uqp}{U^{\prime}_q (\widehat{sl}_2)}
\nc{\uqbp}{U_q (b_+)}
\nc{\uqbm}{U_q (b_-)}
\nc{\ub}{U^{\prime}_q (b_+)}
\nc{\vsl}{V(\sigma(\lambda))}
\nc{\vl}{V(\lambda)}  
\nc{\bu}{\bullet}
\nc{\an}{{\ell}}
\nc{\ws}{\;\;}
\nc{\qu}{{1\ov 4}}
\nc{\hif}{\hb{ if }}
\nc{\hev}{\hb{ is even }}
\nc{\hod}{\hb{ is odd }}
\nc{\Tr}{{\rm Tr}}
\nc{\ad}{{\rm Ad}}
\nc{\hb}{\hbox}
\nc{\curlra}{\buildrel{\sim}\over\longrightarrow}
\nc{\epp}{\varepsilon^{\prime}} 
\nc{\ol}{\overline}
\nc{\pl}{\prod\limits} 
\nc{\sli}{\sum\limits} 
\nc{\nin}{\noindent}
\nc{\ga}{\alpha}
\nc{\gb}{\beta}
\nc{\gd}{\delta}
\nc{\gep}{\varepsilon}
\nc{\gz}{\zeta}
\nc{\gt}{\theta}
\nc{\gk}{\kappa}
\nc{\gl}{\lambda}
\nc{\gp}{\phi}
\nc{\gs}{\sigma}
\nc{\go}{\omega}
\nc{\gn}{\nu}
\nc{\gr}{\rho}

\nc{\s}{\sigma}
\nc{\ep}{\varepsilon}
\nc{\zi}{\zeta^{-1}}

\nc{\gG}{\Gamma}
\nc{\gD}{\Delta}
\nc{\gT}{\Theta}
\nc{\gL}{\Lambda}
\nc{\gO}{\Omega}
\nc{\gP}{\Phi}

\nc{\cF}{\mathcal{F}}
\nc{\cP}{\mathcal{P}}
\nc{\cS}{\mathcal{S}}
\nc{\cN}{\mathcal{N}}
\nc{\cH}{\mathcal{H}}
\nc{\cO}{\mathcal{O}}
\nc{\cT}{\mathcal{T}}
\nc{\cQ}{\mathcal{Q}}
\nc{\cW}{\mathcal{W}}

\nc{\N}{\mathbb{N}}

\nc{\fg}{\mathfrak{g}}

\nc{\bi}{\bar{i}}
\nc{\bj}{\bar{j}}
\nc{\bp}{\bar{p}}
\nc{\bgr}{\bar{\rho}}
\nc{\bA}{\bar{\alpha}}
\nc{\bB}{\bar{\beta}}
\nc{\bC}{\bar{\gamma}}
\nc{\by}{\bar{y}}
\nc{\tf}{\tilde{f}}
\nc{\te}{\tilde{e}}
\nc{\ts}{\tilde{s}}
\nc{\tgP}{\widetilde{\Phi}}
\nc{\tgPs}{\tilde{\Psi}}
\nc{\tgn}{\tilde{\nu}}
\nc{\tgl}{\tilde{\lambda}}
\nc{\tge}{\tilde{\eta}}
\nc{\txi}{\tilde{\xi}}
\nc{\tep}{\tilde{\epsilon}}

\nc{\cB}{\check{b}}

\nc{\goto}{\mapsto}
\nc{\embed}{\hookrightarrow}
\nc{\rien}{\emptyset}
\nc{\Nt}{\frac{N}{2}}
\nc{\vn}{\hspace*{-33truemm}}
\nc{\vm}{\hspace*{-0truemm}}
\nc{\ti}{t^{-1}}
\nc{\vb}{v^{(1)}}
\nc{\vbn}{v^{(n)}}
\nc{\us}{\underline{s}}
\nc{\vz}{V^{(1)}_z}
\nc{\vzn}{V^{(n)}_z}
\nc{\vzo}{V^{(1)}_1}
\nc{\piz}{\pi_z^{(1)}}
\nc{\pizn}{\pi_z^{(n)}}
\nc{\pis}{\pi_{(z,\us)}}
\nc{\bW}{\overline{W}}
\nc{\bQ}{\overline{Q}}
\nc{\tQ}{\widetilde{Q}}
\nc{\bT}{\overline{T}}
\nc{\note}[1]{\vspace*{-5mm}\marginpar[left]{\scriptsize\bf{#1}}}
\nc{\ft}{\footnotesize}

\newcommand{\vth}[3]{\vartheta_{#1}\left( \frac{#2}{r} \bigl| {#3} \right)}
\newcommand{\vtf}[3]{\vartheta_{#1}\left( {#2} \Bigl| {#3} \right)}

\newcommand{\tp}{{\tilde{p}}}

\newcommand{\mmatrix}[1]{\begin{matrix} #1 \end{matrix}}

\newcommand{\baj}{{\bar{j}}}
\newcommand{\bao}{{\bar{1}}}

\newcommand{\bak}{{\bar{k}}}

\newcommand{\noi}{\noindent}
\newcommand{\bell}{\bar{\ell}}
\nc{\tgL}{\widetilde{\Lambda}}
\newcommand{\trace}{{\rm Tr}}
\newcommand{\cZ}{{\mathcal Z}}

\begin{document}
\bibliographystyle{unsrt}
%%%%%%%%%%%%%%%%%%%%%%%%%%%%%%%
%}}}
%{{{ title
\begin{flushright}
%{\bf \Large FINAL DRAFT  }  \today\\[10mm]
\end{flushright}
\begin{center}
{\Large \bf The Vertex-Face Correspondence
   and Correlation Functions of the Fusion Eight-Vertex Model \\I: The General Formalism\\[10mm] }
{\large \bf Takeo Kojima$^1$,  Hitoshi Konno$^2$ and Robert Weston$^3$}\\[6mm]
{\it $^1$ Department of Mathematics, College of Science and Technology, \\ 
     Nihon University, Chiyoda-ku, Tokyo 101-0062, Japan.\\ kojima@math.cst.nihon-u.ac.jp} \\[5mm] 
{\it $^2$Department of Mathematics, Faculty of Integrated Arts $\&$ Sciences, \\Hiroshima University, Higashi-Hiroshima 739-8521, Japan. \\
       konno@mis.hiroshima-u.ac.jp} \\[5mm]
{\it $^3$Department of Mathematics, Heriot-Watt University,\\
Edinburgh EH14 4AS, UK. R.A.Weston@ma.hw.ac.uk}\\[5mm]
April 2005\\[10mm]
\end{center}
%}}}
%{{{ abstract ...
\begin{abstract}
\noindent 
Based on the vertex-face correspondence, we give an algebraic analysis 
 formulation of correlation functions of the $k\times k$ fusion eight-vertex model 
 in terms of the corresponding fusion SOS model. Here $k\in \Z_{>0}$.
A general formula for correlation functions is derived as a trace over 
the space of states of lattice operators such as the 
 corner transfer matrices,  the half transfer matrices (vertex operators)   
and the tail operator. 
 We give a realization of these lattice operators as well as 
  the space of states 
  as objects in the level $k$ representation
   theory of the elliptic algebra $U_{q,p}(\slth)$. 
\end{abstract}

\newpage
\section{Introduction}
The eight-vertex model was solved by Baxter in the series
of seminal papers \cite{Bax72a,Bax73aI,Bax73aII,Bax73aIII}. One of the key insights in these papers was the
realization that by a suitable change of basis it was possible to 
map the model to an `Ising-like' model \cite{Bax73aIII}. This model, which we now
refer to as an SOS model, possessed the property of charge
conservation through a vertex, and its transfer matrix could be
diagonalised using a conventional Bethe ansatz. The height
restricted versions of the SOS model, now called RSOS or ABF models,
later achieved independent fame, largely due to their connection
with conformal field theory models \cite{ABF,Huse84}.

The method of fusion, leading to higher-spin analogues of the
eight-vertex model, was developed in \cite{KRYS81} and \cite{KS82}. 
Baxter's vertex-face correspondence between the eight-vertex and 
SOS models was then extended to these fusion models in \cite{DJMO86}.
As a result, fusion, or higher-spin, SOS models were defined and 
studied in great detail \cite{DJKMO87,DJKMO88}.

In the early 1990s a new approach to solvable lattice models was
developed by Jimbo, Miwa and their collaborators \cite{Daval,JM}.
The approach
was applied originally and most fully to the 6-vertex model. 
The central idea was to exploit to the fullest possible extent
an underlying $\uq$ symmetry of the infinite lattice model.
The transfer matrix and its associated vector space, the eigenstates,
and the local operators of the model were all constructed in terms of
this algebra and its associated vertex operators.
This enabled this group to express all correlation
functions of the six-vertex model in terms of traces of algebraic 
objects over highest-weight representations of $\uq$. One final 
ingredient, a free-field realization of the algebra, was then used 
to compute these traces, thus yielding multiple integral expressions
for correlation functions \cite{JM}.

The success of this method which we shall call the `algebraic
analysis' approach
led to a great deal of subsequent work in which a whole variety of
models were considered from a similar point of view. In particular,
a parallel discussion of the eight-vertex model in terms of an
elliptic algebra, $\cA_{q,p}(\slth)$\cite{FIJKMYa,FIJKMYb,JKOS1}, 
was presented in
\cite{JMN,JKKMW}. However, a free-field realization
of $\cA_{q,p}(\slth)$ was and still is lacking. The difficulty is
again essentially related to the lack of charge conservation for the
eight-vertex model.
Thus, while expressions for correlation functions as 
traces over highest-weight modules exist, it has not
proved possible in general to evaluate
these traces.

RSOS models were also considered using the algebraic analysis approach
\cite{JMO93,LP96}. In this case, a free-field realization 
of the vertex operators appearing in the trace formula was constructed
and the trace
was computed \cite{LP96}. What is more, this was originally done, perhaps
surprisingly, in the absence of a full understanding of the underlying
symmetry algebra.

The correlation functions of the eight-vertex model were finally
computed in a beautiful piece of work by Lashkevich and Pugai
 \cite{LaP98,La02}. Their approach
was to use and extend Baxter's  vertex-face correspondence in order
to map the trace expression for correlation functions of the
eight-vertex model into one for SOS models. They then computed the 
latter using the free-field realization of \cite{LP96}. An 
interesting aspect of this work was that correlation functions 
of the eight-vertex model were found to correspond
to correlation functions
of the SOS model incorporating a certain dislocation 
or `tail' operator. Furthermore, this tail operator had a surprisingly
simple realization in terms of the SOS free-field realization.

The elliptic algebra $U_{q,p}(\slth)$ associated with fusion SOS
models was first defined in terms of elliptic Drinfeld currents in
\cite{Ko98}. This algebra, or more precisely the closely related
algebra ${\cal{B}}_{q,\gl}(\slth)$\cite{JKOS1},
was later interpreted as a quasi-Hopf twisting of $\uq$ in
\cite{JKOS}. A free-field realization of $U_{q,p}(\slth)$ was constructed
in \cite{Ko98,JKOS} and the level one case was 
shown to be equivalent to the `phenomenologically' derived
realization
of \cite{LP96}.

In this paper, we revisit the approach of Lashkevich and Pugai \cite{LaP98} in
light of these recent developments in the understanding of the
underlying
symmetry algebra $U_{q,p}(\slth)$. Our twin motivations to carry
out this work were: to generalise the approach of \cite{LaP98} to the higher-spin
fusion vertex models, and obtain as many explicit results as possible;
and to construct the objects appearing in the correlation function traces,
such as vertex and tail operators, directly in terms of the algebra
$U_{q,p}(\slth)$. 

In Section 2, we construct the higher-spin
vertex and SOS weights, and the higher-spin intertwiners that relate
them. In Section 3, we review the relevant aspects of the algebraic analysis 
approach for vertex and SOS models. In Section 4, we
generalise the graphical arguments of \cite{LaP98} in
order to
connect correlation functions of the higher-spin vertex models to 
those of SOS models with an associated tail operator. In Section 5, 
we review the construction of $U_{q,p}(\slth)$ and give a direct 
algebraic construction of the vertex operators and the tail operators
occurring in the SOS trace expression for vertex model correlation
functions. One of our key results is contained in Conjecture 5.9, 
which gives a remarkably simple algebraic picture of the tail operator
as a simple integer powers of one of the half-currents occurring
in the definition of $U_{q,p}(\slth)$.
In a subsequent paper \cite{KKW04b}, we shall consider the 
spin-$1$ generalisation of the eight-vertex model in detail. This
corresponds to level $k=2$ case of the general formalism given in
the present paper.  

\section{The Fusion Vertex and SOS Models}
In this section, we define and relate the fusion vertex and SOS models
that are of interest to us in this paper. They are statistical-mechanical
models whose Boltzmann weights are expressed in terms of elliptic
theta functions.

\subsection{Notation}
First, let us fix our notation for theta functions.
 Let ${p}=e^{-\frac{\pi K'}{K}}$, $q=-e^{-\frac{\pi \lambda}{2K}}$ and $\zeta=e^{-\frac{\pi \la u}{2K}}$.
 We introduce $x$, $\tau$ and $r$ by $x=-q$, 
 $p=e^{-\frac{2\pi i}{ \tau}}=x^{2r}$, i.e., $\tau=\frac{2iK}{K'}$ and $r=\frac{K'}{\la}$. 
 Throughout this paper ${\rm Im}\tau >0$.

We use the theta functions defined in terms of 
$\tp=e^{2\pi i \tau}$ by
\be
&&\vartheta_1(u|\tau)=2\tilde{p}^{1/8}(\tilde{p};\tilde{p})_\infty\sin\pi u
\prod_{n=1}^\infty(1-2\tilde{p}^n\cos2\pi u+\tilde{p}^{2n}),\\
&&\vartheta_0(u|\tau)=-ie^{\pi i(u+\tau/4)}\vartheta_1\left(u+\frac{\tau}{2}\Big|\tau\right),\\
&&\vartheta_2(u|\tau)=\vartheta_1\left(u+\frac{1}{2}\Big|\tau\right),\\
&&\vartheta_3(u|\tau)=e^{\pi i(u+\tau/4)}\vartheta_1\left(u+\frac{\tau+1}{2}\Big|\tau\right).
\en
We define the symbols $[u]^{(s)}$, $[u]$ and $[u]^*$,  by
\be
&&[u]^{(s)}=x^{\frac{u^2}{s}-u}\Theta_{x^{2s}}(x^{2u})=C\vth{1}{u}{\tau},\quad
C= x^{-\frac{ r}{4}}e^{-\frac{\pi  i }{4}}\ \tau^{1/2},\\
&&[u]=[u]^{(r)},\quad [u]^*=[u]^{(r-k)}.\en
with\\[-10mm]
\be
&&\Theta_p(z)=(z;p)_\infty(p/z;p)_\infty(p;p)_\infty,\quad
(z;p_1,p_2,\cdots,
p_m)_\infty=\hspace*{-3mm}\prod_{n_1,n_2,\cdots,n_m=0}^\infty
\hspace*{-5mm} (1-zp_1^{n_1}p_2^{n_2}\cdots p_m^{n_m}).
\en
Furthermore, we also make use of functions defined in terms of $[u]$ by
\be
&&[A]_M=[A][A-1]\cdots [A-M+1],\\
&&[A,B]=[A][A+1]\cdots [B] \quad (A<B),\qquad [A,A-1]=1,\\
&&\left[\mmatrix{A\cr B\cr}\right]=\frac{[A][A-1]\cdots [A-B+1]}{[B][B-1]\cdots
[1]},\\
&&(a,b)_M=(b,a)_M=\left[\mmatrix{M\cr
                                 \frac{a-b+M}{2}\cr}\right]^{-1}
                                 \frac{\left[\frac{a+b-M}{2}, \frac{a+b+M}{2}
                                  \right]}{\sqrt{[a][b]}}.
\en

In order to distinguish the above notation from that of the $q$-integer, we use the following notation for the later:
\ben 
[\![n]\!]_x=\frac{x^{n}-x^{-n}}{x-x^{-1}}.\een

 Finally, the relation of our $\tau,\ \lambda,\ r,\ u$ to $\tau, \lambda, L,\ u$ in Date et.al.\cite{DJKMO88} is as follows. 
 $\tau_{DJKMO}=-\frac{1}{\tau},\ \la_{DJKMO}=i \la,\ L_{DJKMO}=-\tau r$ and $u_{DJKMO}=-\la u$.

\subsection{The Fusion Vertex Models}
\subsubsection{The eight-vertex model}
The eight-vertex model is a two-dimensional square lattice model. The dynamical variable $\vep_j$ takes 
the values $+$ or $-$. For each vertex, we associate a variable $\vep_j$ with each edge $j$. We allow only the 
eight-configurations for each vertex as depicted in Figure 2.1 $(b)$.
\begin{center}
\setlength{\unitlength}{0.008mm}%
\begingroup\makeatletter\ifx\SetFigFont\undefined%
\gdef\SetFigFont#1#2#3#4#5{%
  \reset@font\fontsize{#1}{#2pt}%
  \fontfamily{#3}\fontseries{#4}\fontshape{#5}%
  \selectfont}%
\fi\endgroup%
\begin{picture}(4212,3648)(601,-3340)
\put(3451,-3286){\makebox(0,0)[lb]{\smash{\SetFigFont{12}{14.4}{\rmdefault}{\mddefault}{\updefault}{\color[rgb]{0,0,0}$(a)$}%
}}}
\thinlines
{\color[rgb]{0,0,0}\put(4801,-961){\vector(-1, 0){2400}}
}%
\put(4726,-811){\makebox(0,0)[lb]{\smash{\SetFigFont{12}{14.4}{\rmdefault}{\mddefault}{\updefault}{\color[rgb]{0,0,0}$\vep_2$}%
}}}
\put(3751,164){\makebox(0,0)[lb]{\smash{\SetFigFont{12}{14.4}{\rmdefault}{\mddefault}{\updefault}{\color[rgb]{0,0,0}$\vep_1$}%
}}}
\put(3751,-2000){\makebox(0,0)[lb]{\smash{\SetFigFont{12}{14.4}{\rmdefault}{\mddefault}{\updefault}{\color[rgb]{0,0,0}$\vep_1'$}%
}}}
\put(2501,-811){\makebox(0,0)[lb]{\smash{\SetFigFont{12}{14.4}{\rmdefault}{\mddefault}{\updefault}{\color[rgb]{0,0,0}$\vep_2'$}%
}}}
\put(2001,-961){\makebox(0,0)[lb]{\smash{\SetFigFont{12}{14.4}{\rmdefault}{\mddefault}{\updefault}{\color[rgb]{0,0,0}$v$}%
}}}
\put(3601,-2461){\makebox(0,0)[lb]{\smash{\SetFigFont{12}{14.4}{\rmdefault}{\mddefault}{\updefault}{\color[rgb]{0,0,0}$u$}%
}}}
\put(-1500,-1036){\makebox(0,0)[lb]{\smash{\SetFigFont{12}{14.4}{\rmdefault}{\mddefault}{\updefault}{\color[rgb]{0,0,0}$R(u-v)^{\vep_1\vep_2}_{\vep_1'\vep_2'}$}%
}}}
\put(1426,-1036){\makebox(0,0)[lb]{\smash{\SetFigFont{12}{14.4}{\rmdefault}{\mddefault}{\updefault}{\color[rgb]{0,0,0}=}%
}}}
{\color[rgb]{0,0,0}\put(3601,239){\vector( 0,-1){2400}}
}%
\end{picture}{}\qquad\qquad\qquad{}

\vspace{5mm}
 \setlength{\unitlength}{0.012mm}%
\begingroup\makeatletter\ifx\SetFigFont\undefined%
\gdef\SetFigFont#1#2#3#4#5{%
  \reset@font\fontsize{#1}{#2pt}%
  \fontfamily{#3}\fontseries{#4}\fontshape{#5}%
  \selectfont}%
\fi\endgroup%
\begin{picture}(11724,2466)(589,-2215)
\put(11176,-286){\makebox(0,0)[lb]{\smash{\SetFigFont{12}{14.4}{\rmdefault}{\mddefault}{\updefault}{\color[rgb]{0,0,0}$-$}%
}}}
\thinlines
{\color[rgb]{0,0,0}\put(1801,-361){\vector(-1, 0){1200}}
}%
{\color[rgb]{0,0,0}\put(3301,-361){\vector(-1, 0){1200}}
}%
{\color[rgb]{0,0,0}\put(2701,239){\vector( 0,-1){1200}}
}%
{\color[rgb]{0,0,0}\put(4801,-361){\vector(-1, 0){1200}}
}%
{\color[rgb]{0,0,0}\put(4201,239){\vector( 0,-1){1200}}
}%
{\color[rgb]{0,0,0}\put(6301,-361){\vector(-1, 0){1200}}
}%
{\color[rgb]{0,0,0}\put(5701,239){\vector( 0,-1){1200}}
}%
{\color[rgb]{0,0,0}\put(7801,-361){\vector(-1, 0){1200}}
}%
{\color[rgb]{0,0,0}\put(7201,239){\vector( 0,-1){1200}}
}%
{\color[rgb]{0,0,0}\put(9301,-361){\vector(-1, 0){1200}}
}%
{\color[rgb]{0,0,0}\put(8701,239){\vector( 0,-1){1200}}
}%
{\color[rgb]{0,0,0}\put(10801,-361){\vector(-1, 0){1200}}
}%
{\color[rgb]{0,0,0}\put(10201,239){\vector( 0,-1){1200}}
}%
{\color[rgb]{0,0,0}\put(12301,-361){\vector(-1, 0){1200}}
}%
{\color[rgb]{0,0,0}\put(11701,239){\vector( 0,-1){1200}}
}%
\put(1351, 89){\makebox(0,0)[lb]{\smash{\SetFigFont{12}{14.4}{\rmdefault}{\mddefault}{\updefault}{\color[rgb]{0,0,0}$+$}%
}}}
\put(1626,-286){\makebox(0,0)[lb]{\smash{\SetFigFont{12}{14.4}{\rmdefault}{\mddefault}{\updefault}{\color[rgb]{0,0,0}+}%
}}}
\put(676,-286){\makebox(0,0)[lb]{\smash{\SetFigFont{12}{14.4}{\rmdefault}{\mddefault}{\updefault}{\color[rgb]{0,0,0}+}%
}}}
\put(1351,-961){\makebox(0,0)[lb]{\smash{\SetFigFont{12}{14.4}{\rmdefault}{\mddefault}{\updefault}{\color[rgb]{0,0,0}+}%
}}}
\put(2851, 14){\makebox(0,0)[lb]{\smash{\SetFigFont{12}{14.4}{\rmdefault}{\mddefault}{\updefault}{\color[rgb]{0,0,0}+}%
}}}
\put(3126,-286){\makebox(0,0)[lb]{\smash{\SetFigFont{12}{14.4}{\rmdefault}{\mddefault}{\updefault}{\color[rgb]{0,0,0}+}%
}}}
\put(2851,-961){\makebox(0,0)[lb]{\smash{\SetFigFont{12}{14.4}{\rmdefault}{\mddefault}{\updefault}{\color[rgb]{0,0,0}$-$}%
}}}
\put(2176,-286){\makebox(0,0)[lb]{\smash{\SetFigFont{12}{14.4}{\rmdefault}{\mddefault}{\updefault}{\color[rgb]{0,0,0}$-$}%
}}}
\put(5926,-2161){\makebox(0,0)[lb]{\smash{\SetFigFont{12}{14.4}{\rmdefault}{\mddefault}{\updefault}{\color[rgb]{0,0,0}$(b)$}%
}}}
\put(4276, 14){\makebox(0,0)[lb]{\smash{\SetFigFont{12}{14.4}{\rmdefault}{\mddefault}{\updefault}{\color[rgb]{0,0,0}+}%
}}}
\put(4626,-286){\makebox(0,0)[lb]{\smash{\SetFigFont{12}{14.4}{\rmdefault}{\mddefault}{\updefault}{\color[rgb]{0,0,0}$-$}%
}}}
\put(4276,-961){\makebox(0,0)[lb]{\smash{\SetFigFont{12}{14.4}{\rmdefault}{\mddefault}{\updefault}{\color[rgb]{0,0,0}+}%
}}}
\put(3676,-286){\makebox(0,0)[lb]{\smash{\SetFigFont{12}{14.4}{\rmdefault}{\mddefault}{\updefault}{\color[rgb]{0,0,0}$-$}%
}}}
\put(5706, 89){\makebox(0,0)[lb]{\smash{\SetFigFont{12}{14.4}{\rmdefault}{\mddefault}{\updefault}{\color[rgb]{0,0,0}+}%
}}}
\put(6051,-286){\makebox(0,0)[lb]{\smash{\SetFigFont{12}{14.4}{\rmdefault}{\mddefault}{\updefault}{\color[rgb]{0,0,0}$-$}%
}}}
\put(5776,-961){\makebox(0,0)[lb]{\smash{\SetFigFont{12}{14.4}{\rmdefault}{\mddefault}{\updefault}{\color[rgb]{0,0,0}$-$}%
}}}
\put(5106,-286){\makebox(0,0)[lb]{\smash{\SetFigFont{12}{14.4}{\rmdefault}{\mddefault}{\updefault}{\color[rgb]{0,0,0}+}%
}}}
\put(7276, 89){\makebox(0,0)[lb]{\smash{\SetFigFont{12}{14.4}{\rmdefault}{\mddefault}{\updefault}{\color[rgb]{0,0,0}$-$}%
}}}
\put(7551,-286){\makebox(0,0)[lb]{\smash{\SetFigFont{12}{14.4}{\rmdefault}{\mddefault}{\updefault}{\color[rgb]{0,0,0}+}%
}}}
\put(7276,-961){\makebox(0,0)[lb]{\smash{\SetFigFont{12}{14.4}{\rmdefault}{\mddefault}{\updefault}{\color[rgb]{0,0,0}$-$}%
}}}
\put(6600,-286){\makebox(0,0)[lb]{\smash{\SetFigFont{12}{14.4}{\rmdefault}{\mddefault}{\updefault}{\color[rgb]{0,0,0}+}%
}}}
\put(8776, 89){\makebox(0,0)[lb]{\smash{\SetFigFont{12}{14.4}{\rmdefault}{\mddefault}{\updefault}{\color[rgb]{0,0,0}$-$}%
}}}
\put(9051,-286){\makebox(0,0)[lb]{\smash{\SetFigFont{12}{14.4}{\rmdefault}{\mddefault}{\updefault}{\color[rgb]{0,0,0}+}%
}}}
\put(8776,-961){\makebox(0,0)[lb]{\smash{\SetFigFont{12}{14.4}{\rmdefault}{\mddefault}{\updefault}{\color[rgb]{0,0,0}+}%
}}}
\put(8176,-286){\makebox(0,0)[lb]{\smash{\SetFigFont{12}{14.4}{\rmdefault}{\mddefault}{\updefault}{\color[rgb]{0,0,0}$-$}%
}}}
\put(10276, 89){\makebox(0,0)[lb]{\smash{\SetFigFont{12}{14.4}{\rmdefault}{\mddefault}{\updefault}{\color[rgb]{0,0,0}$-$}%
}}}
\put(10526,-286){\makebox(0,0)[lb]{\smash{\SetFigFont{12}{14.4}{\rmdefault}{\mddefault}{\updefault}{\color[rgb]{0,0,0}$-$}%
}}}
\put(10276,-961){\makebox(0,0)[lb]{\smash{\SetFigFont{12}{14.4}{\rmdefault}{\mddefault}{\updefault}{\color[rgb]{0,0,0}+}%
}}}
\put(9676,-286){\makebox(0,0)[lb]{\smash{\SetFigFont{12}{14.4}{\rmdefault}{\mddefault}{\updefault}{\color[rgb]{0,0,0}+}%
}}}
\put(11776, 89){\makebox(0,0)[lb]{\smash{\SetFigFont{12}{14.4}{\rmdefault}{\mddefault}{\updefault}{\color[rgb]{0,0,0}$-$}%
}}}
\put(12101,-286){\makebox(0,0)[lb]{\smash{\SetFigFont{12}{14.4}{\rmdefault}{\mddefault}{\updefault}{\color[rgb]{0,0,0}$-$}%
}}}
\put(11776,-961){\makebox(0,0)[lb]{\smash{\SetFigFont{12}{14.4}{\rmdefault}{\mddefault}{\updefault}{\color[rgb]{0,0,0}$-$}%
}}}
{\color[rgb]{0,0,0}\put(1201,239){\vector( 0,-1){1200}}
}%
\end{picture}\lb{evm}

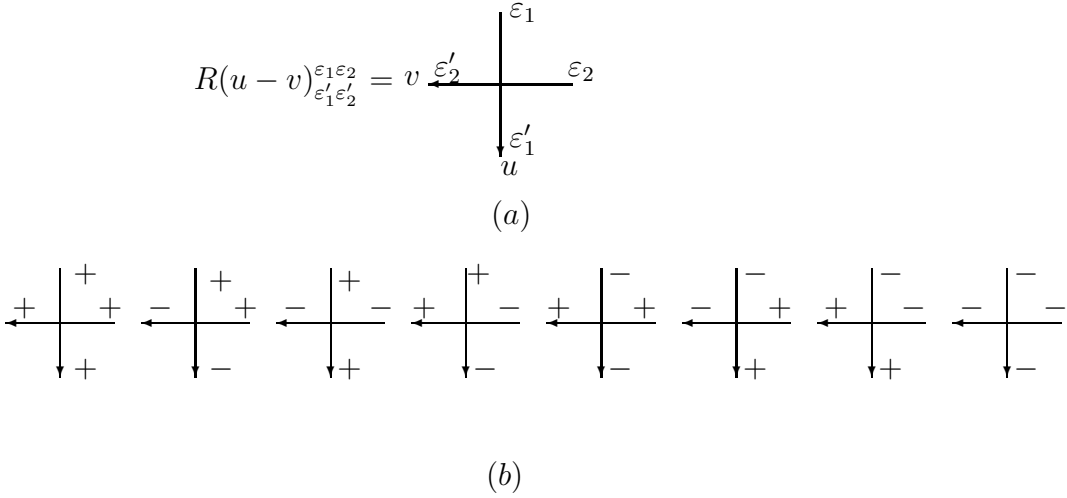
{\footnotesize Figure 2.1:  The eight-vertex model: $(a)$ the $R$-matrix;
$(b)$ the eight possible configurations.}
\end{center}

We assign the following     
Boltzmann weight ( the $R$ matrix ) $R(u)^{\vep_1\vep_2}_{\vep'_1\vep'_2}$ to each configuration: 
\bea
&&{R}(u)=
{{R}_0(u)}\left(\mmatrix{a(u)&&&d(u)\cr&b(u)&c(u)&\cr&c(u)&b(u)&\cr d(u)&&&a(u)\cr}\right)
,\lb{DJKMO-R4}
\ena
where
\bea
R_0(u)&=&z^{-\frac{r-1}{2r}}\frac{(pq^2z;q^4,p)_\infty(q^2z;q^4,p)_\infty
(p/z;q^4,p)_\infty(q^4/z;q^4,p)_\infty}{(pq^2/z;q^4,p)_\infty
(q^2/z;q^4,p)_\infty(pz;q^4,p)_\infty(q^4z;q^4,p)_\infty},\lb{evR}\\[2mm]
a(u)&=&\frac{\vtf{2}{\frac{1}{2r}}{\frac{\tau}{2}}\vtf{2}{\frac{  u}{2r}}{\frac{\tau}{2}}}{\vtf{2}{0}{\frac{\tau}{2}}\vtf{2}{\frac{1+u}{2r}}{\frac{\tau}{2}}},\qquad
b(u)=\frac{\vtf{2}{\frac{1 }{2r}}{\frac{\tau}{2}}\vtf{1}{\frac{  u}{2r}}{\frac{\tau}{2}}}
{\vtf{2}{0}{\frac{\tau}{2}}\vtf{1}{\frac{1+u}{2r}}{\frac{\tau}{2}}},\\
c(u)&=&\frac{\vtf{1}{\frac{ 1}{2r}}{\frac{\tau}{2}}\vtf{2}{\frac{  u}{2r}}{\frac{\tau}{2}}}
{\vtf{2}{0}{\frac{\tau}{2}}\vtf{1}{\frac{1+u}{2r}}{\frac{\tau}{2}}},\qquad
d(u)=-\frac{\vtf{1}{\frac{ 1}{2r}}{\frac{\tau}{2}}\vtf{1}{\frac{  u}{2r}}{\frac{\tau}{2}}}
{\vtf{2}{0}{\frac{\tau}{2}}\vtf{2}{\frac{1+u}{2r}}{\frac{\tau}{2}}},
\ena
with $z=\zeta^2=x^{2u}$. The extra parameter $u$ is called the spectral parameter. 
Let $V$ denote a two-dimensional vector space  spanned by $v_+,\ v_-$. We regard $R(u)$ as an operator 
acting on $V \otimes V$. 
\be
&&R(u)v_{\vep_1}\otimes v_{\vep_2}=\sum_{\vep'_1,\, \vep'_2=+,-} R(u)^{\vep_1\vep_2}_{\vep'_1\vep'_2}
v_{\vep'_1}\otimes v_{\vep'_2}.
\en
The $R$-matrix satisfies the Yang-Baxter equation, unitarity and crossing symmetry relations
 given as follows.
\bea
&&\hspace*{-10mm}R_{12}(u_1-u_2)R_{13}(u_1-u_3)R_{23}(u_2-u_3)=R_{23}(u_2-u_3)R_{13}(u_1-u_3)R_{12}(u_1-u_2),\lb{ybeq}\\
&&\hspace*{-10mm}R(u)PR(-u)P=I,\lb{uv}\\
&&\hspace*{-10mm}(PR(u)P)^{t_1}=(\sigma^y\otimes 1)R(-u-1)(\sigma^y\otimes 1)^{-1},\lb{cross}
\ena
where 
$P(v_{\vep_1}\otimes v_{\vep_2})= v_{\vep_2}\otimes v_{\vep_1}$, and ${}^{t_1}$ denotes transposition with 
respect to the first vector space in the tensor product. 
Note that $R_0(u)$ satisfies the following inversion relations.
\be
&&R_0(u)R_0(-u)=1,\\
&&R_0(u)R_0(u+1)=-\frac{[u+1]}{[u]}.
\en

\subsubsection{Fusion of the eight-vertex model}

Define the operator $\Pi_{1 \cdots k}$ by
\be
&&\Pi_{1\cdots k}=\frac{1}{k!}(P_{1k}+\cdots+P_{k-1k}+I)\cdots(P_{13}+P_{23}+I) (P_{12}+I).
\en
This yields the projection operator onto the space $V^{(k)}$ of symmetric tensors in $V^{\otimes k}$. A basis for $V^{(k)}$ is given by $\{ v^{(k)}_\ep \}_{\ep=-k,-k+2,..,k}$, where $\ep=\ep_1+\ep_2+\cdots \ep_k$ (with $\ep_i\in\{+,-\}$), and 
\bea
v^{(k)}_\ep&=&\Pi_{1\cdots k} v_{\vep_1}\otimes v_{\vep_2} \otimes \cdots \otimes v_{\vep_k}\nn\\
&=&\frac{1}{k!}\sum_{\sigma \in S_k}v_{\vep_{\sigma(1)}}\otimes v_{\vep_{\sigma(2)}}\otimes \cdots \otimes v_{\vep_{\sigma(k)}},\lb{fusionv}
\ena
where $S_k$ is the symmetric group. 

We then define an operator  
\be
&&R_{1\cdots k, \baj}(u)=\Pi_{1\cdots k}R_{1\baj}(u+k-1)\cdots R_{k-1\baj}(u+1)R_{k\baj}(u)\ \in {\rm End}(V^{(k)} \otimes V_{\baj}).
\en
This satisfies 
\be
&&R_{1\cdots k, \baj}(u)\Pi_{1\cdots k}=R_{1\cdots k, \baj}(u).
\en
The $k\times k$ fusion of the $R$-matrix is then given by
\bea
&&R^{(k,k)}(u)=\Pi_{\bao\cdots \bak}R_{1\cdots k, \bak}(u)R_{1\cdots k, \overline{k-1}}(u-1)\cdots R_{1\cdots k, \bao}(u-k+1). \lb{kkfusionR}
\ena
This is an operator in ${\rm End}(V^{(k)} \otimes V^{(k)})$. It
satisfies 
\bea
&&R^{(k,k)}(u)=R^{(k,k)}(u)\Pi_{1\cdots k}=R^{(k,k)}(u)\Pi_{\bao\cdots \bak}.\lb{projkk}
\ena
Using the YBE \eqref{ybeq} repeatedly and \eqref{projkk}, we verify that $R^{(k,k)}(u)$ satisfies the YBE 
on $V^{(k)}\otimes V^{(k)} \otimes V^{(k)}$. 
$R^{(k,k)}(u)$ also satisfies the unitarity condition which is the
simple higher $k$ version of \mref{uv}. 

The situation with
crossing symmetry is somewhat more complicated. For general $k$ we 
have in general the following relation\cite{JM}. 
\bea
 (P^{(k)}R^{(k,k)}(u) P^{(k)})^{t_1} =
(Q\ot \id) R^{(k,k)}(-u-1) (Q^{-1}\ot \id),\lb{qconj}
\ena
where $P^{(k)}$ is the permutation operator acting on $V^{(k)}\ot
V^{(k)}$, and $Q$ is a $k+1$ dimensional matrix whose entries are
independent
of $u$. Clearly, for $k=1$, it follows from \mref{cross} that $Q=\gs^y$. For $k=2$, we find by
explicit
calculation that\cite{Konno04}
\be Q&=&\frac{1}{2}\left(
\begin{matrix}1+y^2&0&1-y^2\cr 0&x^2&0\cr
      1-y^2&0&1+y^2\end{matrix}\right)\quad \hb{where}\\[3mm]
x^2&=&{-\frac{1}{2}\frac{\vtf{0}{0}{\tau}\vtf{3}{\frac{1}{\tau}}{\tau}}
{\vtf{0}{\frac{1}{\tau}}{\tau}\vtf{3}{0}{\tau}}},\quad 
y^2={-\frac{\vtf{2}{0}{\tau}\vtf{3}{\frac{1}{\tau}}{\tau}}{\vtf{2}{\frac{1}{\tau}}{\tau}\vtf{3}{0}{\tau}}}.
\en

The matrix elements of
$R^{(k,k)}(u)$ then define the $k\times k$ fusion  
eight-vertex model whose dynamical variables takes values in $\{-k, -k+2,..,k\}$.

\subsubsection{The ground states}

We consider the principal regime specified by 
\be
 &&0<\sf{k}<1,\quad 0<\la<K',\quad -1<u<0,
 \en
where $\sf{k}$ denotes the elliptic modulus. Hence $r>1$.
In this regime, we find that for the $k=1$ eight-vertex model we have
$c>a,|b|,d$, with $a>0, b<0, d>0$.
Maximal Boltzmann weight configurations, which we
refer to as ground states, will
therefore only involve the weight $c$. More generally, for arbitrary
$k>0$, we find that a maximal-weight configuration of edge variables 
is labelled by $\ell\in\{0,1,\cdots,k\}$, and is 
a periodic repetition of the pattern in Figure 2.2, where $\ep_\ell=k-2\ell$.

\setlength{\unitlength}{0.0004in}
\begingroup\makeatletter\ifx\SetFigFont\undefined%
\gdef\SetFigFont#1#2#3#4#5{%
  \reset@font\fontsize{#1}{#2pt}%
  \fontfamily{#3}\fontseries{#4}\fontshape{#5}%
  \selectfont}%
\fi\endgroup%
{\renewcommand{\dashlinestretch}{30}
\begin{picture}(3644,3659)(-4000,500)
%\thinlines
\thinlines
\path(1222,3622)(1222,22)
\path(1162.000,262.000)(1222.000,22.000)(1282.000,262.000)
\path(2422,3622)(2422,22)
\path(2362.000,262.000)(2422.000,22.000)(2482.000,262.000)
\path(262.000,2482.000)(22.000,2422.000)(262.000,2362.000)
\path(22,2422)(3622,2422)
\path(262.000,1282.000)(22.000,1222.000)(262.000,1162.000)
\path(22,1222)(3622,1222)
\put(2572,622){\makebox(0,0)[lb]{{\SetFigFont{10}{16.8}{\rmdefault}{\mddefault}{\updefault}$\ep_{\ell}$}}}
\put(2572,3022){\makebox(0,0)[lb]{{\SetFigFont{10}{16.8}{\rmdefault}{\mddefault}{\updefault}$\ep_{\ell}$}}}
\put(1297,1822){\makebox(0,0)[lb]{{\SetFigFont{10}{16.8}{\rmdefault}{\mddefault}{\updefault}$\ep_{\ell}$}}}
\put(2572,1822){\makebox(0,0)[lb]{{\SetFigFont{10}{16.8}{\rmdefault}{\mddefault}{\updefault}$-\ep_{\ell}$}}}
\put(1297,622){\makebox(0,0)[lb]{{\SetFigFont{10}{16.8}{\rmdefault}{\mddefault}{\updefault}$-\ep_{\ell}$}}}
\put(1372,3022){\makebox(0,0)[lb]{{\SetFigFont{10}{16.8}{\rmdefault}{\mddefault}{\updefault}$-\ep_{\ell}$}}}
\put(547,2497){\makebox(0,0)[lb]{{\SetFigFont{10}{16.8}{\rmdefault}{\mddefault}{\updefault}$-\ep_{\ell}$}}}
\put(1747,2497){\makebox(0,0)[lb]{{\SetFigFont{10}{16.8}{\rmdefault}{\mddefault}{\updefault}$\ep_{\ell}$}}}
\put(2947,2497){\makebox(0,0)[lb]{{\SetFigFont{10}{16.8}{\rmdefault}{\mddefault}{\updefault}$-\ep_{\ell}$}}}
\put(2947,1297){\makebox(0,0)[lb]{{\SetFigFont{10}{16.8}{\rmdefault}{\mddefault}{\updefault}$\ep_{\ell}$}}}
\put(1747,1297){\makebox(0,0)[lb]{{\SetFigFont{10}{16.8}{\rmdefault}{\mddefault}{\updefault}$-\ep_{\ell}$}}}
\put(547,1297){\makebox(0,0)[lb]{{\SetFigFont{10}{16.8}{\rmdefault}{\mddefault}{\updefault}$\ep_{\ell}$}}}
\lb{vgs}\end{picture}
}

\vspace*{4mm}\begin{center}
{\footnotesize
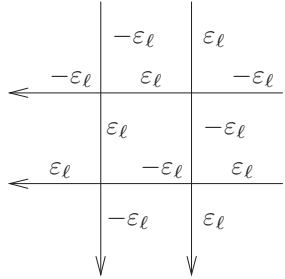
Figure 2.2: The ground state configuration of the fusion vertex
model labelled by $\ell$
}
\end{center}

\subsection{Fusion SOS models}

\subsubsection{The eight-vertex SOS model}
The eight-vertex SOS model, usually referred to as simply the SOS
model, is also a two-dimensional square lattice model\cite{ABF}. The dynamical variables $a_j$ are called local heights.
They take values in $\Z$. For each face, we associate a local height $a_j$ with each vertex $j$. 
We allow only the configurations satisfying the so-called admissibility condition $|a_j-a_k|=1$ for 
any two adjacent local heights $a_j$ and $a_k$. Then we have only the six possible 
configurations for each face depicted in Figure 2.3 $(b)$.
\be
{}\qquad \qquad \setlength{\unitlength}{0.01mm}%
\begingroup\makeatletter\ifx\SetFigFont\undefined%
\gdef\SetFigFont#1#2#3#4#5{%
  \reset@font\fontsize{#1}{#2pt}%
  \fontfamily{#3}\fontseries{#4}\fontshape{#5}%
  \selectfont}%
\fi\endgroup%
\begin{picture}(4500,3264)(526,-2040)
\put(2101,-2986){\makebox(0,0)[lb]{\smash{\SetFigFont{12}{14.4}{\rmdefault}{\mddefault}{\updefault}{\color[rgb]{0,0,0}$(a)$}%
}}}
\thinlines
{\color[rgb]{0,0,0}\put(3001,-61){\line( 1, 0){1800}}
\put(4801,-61){\line( 0,-1){1800}}
\put(4801,-1861){\line(-1, 0){1800}}
}%
{\color[rgb]{0,0,0}\put(3501,-61){\line(-1,-1){500}}
}%
\put(3726,-1050){\makebox(0,0)[lb]{\smash{\SetFigFont{12}{14.4}{\rmdefault}{\mddefault}{\updefault}{\color[rgb]{0,0,0}$u$}%
}}}
\put(5026,-2011){\makebox(0,0)[lb]{\smash{\SetFigFont{12}{14.4}{\rmdefault}{\mddefault}{\updefault}{\color[rgb]{0,0,0}$c$}%
}}}
\put(4951, 89){\makebox(0,0)[lb]{\smash{\SetFigFont{12}{14.4}{\rmdefault}{\mddefault}{\updefault}{\color[rgb]{0,0,0}b}%
}}}
\put(2701, 89){\makebox(0,0)[lb]{\smash{\SetFigFont{12}{14.4}{\rmdefault}{\mddefault}{\updefault}{\color[rgb]{0,0,0}$a$}%
}}}
\put(2701,-2086){\makebox(0,0)[lb]{\smash{\SetFigFont{12}{14.4}{\rmdefault}{\mddefault}{\updefault}{\color[rgb]{0,0,0}$d$}%
}}}
\put(-1026,-961){\makebox(0,0)[lb]{\smash{\SetFigFont{12}{14.4}{\rmdefault}{\mddefault}{\updefault}{\color[rgb]{0,0,0}$W\BW{a}{b}{d}{c}{u}$}%
}}}
\put(2101,-961){\makebox(0,0)[lb]{\smash{\SetFigFont{12}{14.4}{\rmdefault}{\mddefault}{\updefault}{\color[rgb]{0,0,0}=}%
}}}
{\color[rgb]{0,0,0}\put(3001,-61){\line( 0,-1){1800}}
}%
%\put(3001,-450){\line(1,1){400}}
\lb{SOS}\end{picture}
\en
\be
\setlength{\unitlength}{0.01mm}%
\begingroup\makeatletter\ifx\SetFigFont\undefined%
\gdef\SetFigFont#1#2#3#4#5{%
  \reset@font\fontsize{#1}{#2pt}%
  \fontfamily{#3}\fontseries{#4}\fontshape{#5}%
  \selectfont}%
\fi\endgroup%
\begin{picture}(11937,2739)(376,-1915)
\put(6301,-2461){\makebox(0,0)[lb]{\smash{\SetFigFont{12}{14.4}{\rmdefault}{\mddefault}{\updefault}{\color[rgb]{0,0,0}$(b)$}%
}}}
\thinlines
{\color[rgb]{0,0,0}\put(2701,-61){\line( 0,-1){1200}}
\put(2701,-1261){\line( 1, 0){1200}}
\put(3901,-1261){\line( 0, 1){1200}}
\put(3901,-61){\line(-1, 0){1200}}
}%
{\color[rgb]{0,0,0}\put(6901,-61){\line( 0,-1){1200}}
\put(6901,-1261){\line( 1, 0){1200}}
\put(8101,-1261){\line( 0, 1){1200}}
\put(8101,-61){\line(-1, 0){1200}}
}%
{\color[rgb]{0,0,0}\put(4801,-61){\line( 0,-1){1200}}
\put(4801,-1261){\line( 1, 0){1200}}
\put(6001,-1261){\line( 0, 1){1200}}
\put(6001,-61){\line(-1, 0){1200}}
}%
{\color[rgb]{0,0,0}\put(9001,-61){\line( 0,-1){1200}}
\put(9001,-1261){\line( 1, 0){1200}}
\put(10201,-1261){\line( 0, 1){1200}}
\put(10201,-61){\line(-1, 0){1200}}
}%
{\color[rgb]{0,0,0}\put(11101,-61){\line( 0,-1){1200}}
\put(11101,-1261){\line( 1, 0){1200}}
\put(12301,-1261){\line( 0, 1){1200}}
\put(12301,-61){\line(-1, 0){1200}}
}%
{\color[rgb]{0,0,0}\put(601,-436){\line( 1, 1){375}}
}%
{\color[rgb]{0,0,0}\put(2738,-473){\line( 1, 1){400}}
}%
{\color[rgb]{0,0,0}\put(4838,-473){\line( 1, 1){400}}
}%
{\color[rgb]{0,0,0}\put(6938,-473){\line( 1, 1){400}}
}%
{\color[rgb]{0,0,0}\put(9038,-473){\line( 1, 1){400}}
}%
{\color[rgb]{0,0,0}\put(11138,-473){\line( 1, 1){400}}
}%
\put(451, 89){\makebox(0,0)[lb]{\smash{\SetFigFont{10}{14.4}{\rmdefault}{\mddefault}{\updefault}{\color[rgb]{0,0,0}$a$}%
}}}
\put(1226, 89){\makebox(0,0)[lb]{\smash{\SetFigFont{10}{14.4}{\rmdefault}{\mddefault}{\updefault}{\color[rgb]{0,0,0}$a+1$}%
}}}
\put(76,-1636){\makebox(0,0)[lb]{\smash{\SetFigFont{10}{14.4}{\rmdefault}{\mddefault}{\updefault}{\color[rgb]{0,0,0}$a+1$}%
}}}
\put(1226,-1636){\makebox(0,0)[lb]{\smash{\SetFigFont{10}{14.4}{\rmdefault}{\mddefault}{\updefault}{\color[rgb]{0,0,0}$a+2$}%
}}}
\put(2551, 89){\makebox(0,0)[lb]{\smash{\SetFigFont{10}{14.4}{\rmdefault}{\mddefault}{\updefault}{\color[rgb]{0,0,0}$a$}%
}}}
\put(4651, 89){\makebox(0,0)[lb]{\smash{\SetFigFont{10}{14.4}{\rmdefault}{\mddefault}{\updefault}{\color[rgb]{0,0,0}$a$}%
}}}
\put(6751, 89){\makebox(0,0)[lb]{\smash{\SetFigFont{10}{14.4}{\rmdefault}{\mddefault}{\updefault}{\color[rgb]{0,0,0}$a$}%
}}}
\put(8851, 89){\makebox(0,0)[lb]{\smash{\SetFigFont{10}{14.4}{\rmdefault}{\mddefault}{\updefault}{\color[rgb]{0,0,0}$a$}%
}}}
\put(10951, 89){\makebox(0,0)[lb]{\smash{\SetFigFont{10}{14.4}{\rmdefault}{\mddefault}{\updefault}{\color[rgb]{0,0,0}$a$}%
}}}
\put(3901,-1636){\makebox(0,0)[lb]{\smash{\SetFigFont{10}{14.4}{\rmdefault}{\mddefault}{\updefault}{\color[rgb]{0,0,0}$a$}%
}}}
\put(6001,-1636){\makebox(0,0)[lb]{\smash{\SetFigFont{10}{14.4}{\rmdefault}{\mddefault}{\updefault}{\color[rgb]{0,0,0}$a$}%
}}}
\put(8101,-1636){\makebox(0,0)[lb]{\smash{\SetFigFont{10}{14.4}{\rmdefault}{\mddefault}{\updefault}{\color[rgb]{0,0,0}$a$}%
}}}
\put(10201,-1636){\makebox(0,0)[lb]{\smash{\SetFigFont{10}{14.4}{\rmdefault}{\mddefault}{\updefault}{\color[rgb]{0,0,0}$a$}%
}}}
\put(3351, 89){\makebox(0,0)[lb]{\smash{\SetFigFont{10}{14.4}{\rmdefault}{\mddefault}{\updefault}{\color[rgb]{0,0,0}$a+1$}%
}}}
\put(2476,-1636){\makebox(0,0)[lb]{\smash{\SetFigFont{10}{14.4}{\rmdefault}{\mddefault}{\updefault}{\color[rgb]{0,0,0}$a+1$}%
}}}
\put(5451, 89){\makebox(0,0)[lb]{\smash{\SetFigFont{10}{14.4}{\rmdefault}{\mddefault}{\updefault}{\color[rgb]{0,0,0}$a+1$}%
}}}
\put(6751,-1636){\makebox(0,0)[lb]{\smash{\SetFigFont{10}{14.4}{\rmdefault}{\mddefault}{\updefault}{\color[rgb]{0,0,0}$a+1$}%
}}}
\put(4651,-1636){\makebox(0,0)[lb]{\smash{\SetFigFont{10}{14.4}{\rmdefault}{\mddefault}{\updefault}{\color[rgb]{0,0,0}$a-1$}%
}}}
\put(7551, 89){\makebox(0,0)[lb]{\smash{\SetFigFont{10}{14.4}{\rmdefault}{\mddefault}{\updefault}{\color[rgb]{0,0,0}$a-1$}%
}}}
\put(9551, 89){\makebox(0,0)[lb]{\smash{\SetFigFont{10}{14.4}{\rmdefault}{\mddefault}{\updefault}{\color[rgb]{0,0,0}$a-1$}%
}}}
\put(8851,-1636){\makebox(0,0)[lb]{\smash{\SetFigFont{10}{14.4}{\rmdefault}{\mddefault}{\updefault}{\color[rgb]{0,0,0}$a-1$}%
}}}
\put(10726,-1636){\makebox(0,0)[lb]{\smash{\SetFigFont{10}{14.4}{\rmdefault}{\mddefault}{\updefault}{\color[rgb]{0,0,0}$a-1$}%
}}}
\put(11900, 89){\makebox(0,0)[lb]{\smash{\SetFigFont{10}{14.4}{\rmdefault}{\mddefault}{\updefault}{\color[rgb]{0,0,0}$a-1$}%
}}}
\put(11900,-1636){\makebox(0,0)[lb]{\smash{\SetFigFont{10}{14.4}{\rmdefault}{\mddefault}{\updefault}{\color[rgb]{0,0,0}$a-2$}%
}}}
{\color[rgb]{0,0,0}\put(601,-61){\line( 0,-1){1200}}
\put(601,-1261){\line( 1, 0){1200}}
\put(1801,-1261){\line( 0, 1){1200}}
\put(1801,-61){\line(-1, 0){1200}}
}%
\end{picture}
\en

\begin{center}
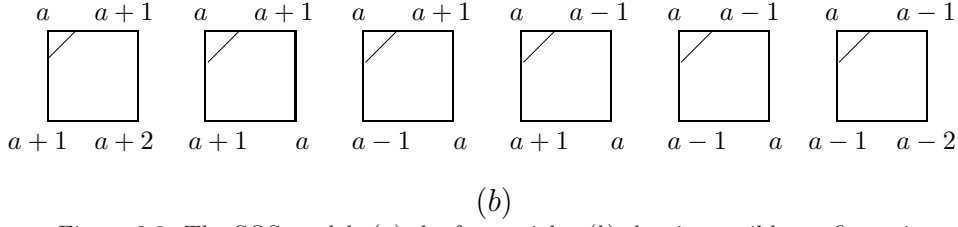
{\footnotesize Figure 2.3: The SOS model: $(a)$ the face weight; $(b)$ the six possible configurations.}
\end{center}
We assign the following     
Boltzmann weight (or face weight) $W\BW{a_1}{a_2}{a_4}{a_3}{u}$ to each configuration.  
  \begin{eqnarray}
&&{W}\left(\left.
\begin{array}{cc}
n&n\pm 1\\
n\pm1&n\pm2
\end{array}\right|u\right)={R}_0(u),\nn\\
&&{W}\left(\left.
\begin{array}{cc}
n&n\pm 1\\
n\pm1&n
\end{array}\right|u\right)={R}_0(u)\frac{[n\mp u][1]}{[ n][ 1+u]},\lb{face}\\
&&{W}\left(\left.
\begin{array}{cc}
n&n\pm 1\\
n\mp 1&n
\end{array}\right|u\right)=
{R}_0(u)
\frac{[ n\pm 1][u] }{[ n][1+u]}.\nn
\end{eqnarray}
The face weights satisfy the following face-type Yang-Baxter equation,
unitarity,
and 
crossing symmetry relations:  
\bea
&&\sum_{g}W\BW{a}{b}{f}{g}{u}W\BW{f}{g}{e}{d}{v}W\BW{b}{c}{g}{d}{u-v}\nn\\
&&\qquad =\sum_{g}W\BW{a}{g}{f}{e}{u-v}W\BW{a}{b}{g}{c}{v}W\BW{g}{c}{e}{d}{u},\lb{fyb}\\
&&\sum_{e}W\BW{a}{b}{e}{c}{u}W\BW{a}{e}{d}{c}{-u}=1,\lb{fcross}\\
&&W\BW{a}{b}{d}{c}{u}=(-)^{\frac{a+d-b-c}{2}}\frac{[b]}{[a]}W\BW{d}{a}{c}{b}{-u-1}.
\ena

\subsubsection{Fusion of the  SOS model}
The $k\times k$ fusion of face weights is obtained as follows:
define
\be
W^{(k,1)}\BW{a}{b}{d}{c}{u}
%\nn\\&&
&=&\sum_{d_1,..,d_{k-1}}  \!\!\!\! W\BW{a}{a_1}{d}{d_1}{u+k-1}W\BW{a_1}{a_2}{d_1}{d_2}{u+k-2}\\
&&\qquad\qquad  \times \, \cdots \, 
W\BW{a_{k-1}}{b}{d_{k-1}}{c}{u}.
\en
Then the RHS is independent of the choice of $a_1,..,a_{k-1}$ provided
$|a-a_1|=|a_1-a_2|=\cdots=|a_{k-1}-b|=1$ \cite{DJKMO88}. We then define
\bea
%\hspace*{-20mm}
W^{(k,k)}\BW{a}{b}{d}{c}{u}
%\nn\\&&\hspace*{-30mm}
&=&\sum_{a_1,..,a_{k-1}}W^{(k,1)}\BW{a}{b}{a_1}{b_1}{u-k+1}
\!\! W^{(k,1)}\BW{a_1}{b_1}{a_2}{b_2}{u-k+2}\nn\\
&&\qquad\qquad \times \, \cdots \, W^{(k,1)}\BW{a_{k-1}}{b_{k-1}}{d}{c}{u}.
% \nn\\
\lb{kkfusionW}
\ena
Here the RHS is independent of the choice of $b_1,..,b_{k-1}$ provided $|b-b_1|=|b_1-b_2|=\cdots=|b_{k-1}-c|=1$.
In $W^{(k,k)}$, the admissible condition for the dynamical variables is extended to $a_j-a_k\ \in \{-k, -k+2, .., k\}$ for any two adjacent local heights $a_j, a_k$. 
The fused face weight $W^{(k,k)}$ satisfies the face type YBE,
unitarity 
and crossing symmetry relations. 
The latter two relations are given by
\bea
&&\sli_{s}\wet{a}{s}{d}{c}{-u}{(k,k)}{} \wet{a}{b}{s}{c}{u}{(k,k)}{}=\delta_{b,d},\label{unitarity}\\
&& \wet{d}{c}{a}{b}{u}{(k,k)}{}=\frac{G^{(k)}_{a,d}}{G^{(k)}_{b,c}} \,\wet{a}{d}{b}{c}{-1-u}{(k,k)}{}.
\label{crossing}
\ena
where $g_a=s_a\sqrt{[a]}$\, $s_a=\pm1,\ s_a s_{a+1}=(-)^a$ and 
$\displaystyle{G^{(k)}_{a,b}
=\frac{g_a}{g_b (a,b)_k}}$.
In addition, we have a symmetry for  $k\in \Z_{>0}$ \cite{DJKMO88}
\bea
&&W^{(k,k)}\BW{d}{a}{c}{b}{u}=\frac{(a,b)_k(d,a)_k}{(d,c)_k(c,b)_k}W^{(k,k)}
\BW{d}{c}{a}{b}{u}.\lb{gaugek}
\ena

\subsubsection{The ground states}
We discuss  regime III specified by the region
\be
0<p<1,\qquad -1<u<0.
\en
In this regime, the ground states of the level $k$ fusion SOS weights are
of the form shown in Figure 2.4, where $m\in \Z$, $\ell\in\{0,1,\cdots,k\}$, and we define $\bar{\ell}= k-\ell$.

%\newpage 
%\input{FIG32.eepic}
\setlength{\unitlength}{0.0004in}
\begingroup\makeatletter\ifx\SetFigFont\undefined%
\gdef\SetFigFont#1#2#3#4#5{%
  \reset@font\fontsize{#1}{#2pt}%
  \fontfamily{#3}\fontseries{#4}\fontshape{#5}%
  \selectfont}%
\fi\endgroup%
{\renewcommand{\dashlinestretch}{30}
\begin{picture}(3882,3180)(-4000,600)
\thinlines
\path(675,2700)(3075,2700)(3075,300)
        (675,300)(675,2700)
\path(1875,2700)(1875,300)
\path(675,1500)(3075,1500)
\path(975,2700)(675,2400)
\path(2175,2700)(1875,2400)
\path(2175,1500)(1875,1200)
\path(975,1500)(675,1200)
\put(1450,-150){\makebox(0,0)[lb]{{\SetFigFont{10}{16.8}{\rmdefault}{\mddefault}{\updefault}$m+\bar{\ell}$}}}
\put(-500,-150){\makebox(0,0)[lb]{{\SetFigFont{10}{16.8}{\rmdefault}{\mddefault}{\updefault}$m+\ell$}}}
\put(3375,-150){\makebox(0,0)[lb]{{\SetFigFont{10}{16.8}{\rmdefault}{\mddefault}{\updefault}     $m+\ell$}}}
\put(3375,1425){\makebox(0,0)[lb]{{\SetFigFont{10}{16.8}{\rmdefault}{\mddefault}{\updefault}$m+\bar{\ell}$}}}
\put(3375,3000){\makebox(0,0)[lb]{{\SetFigFont{10}{16.8}{\rmdefault}{\mddefault}{\updefault}$m+\ell$}}}
\put(1450,3000){\makebox(0,0)[lb]{{\SetFigFont{10}{16.8}{\rmdefault}{\mddefault}{\updefault}$m+\bar{\ell}$}}}
\put(-500,3000){\makebox(0,0)[lb]{{\SetFigFont{10}{16.8}{\rmdefault}{\mddefault}{\updefault}$m+\ell$}}}
\put(-500,1500){\makebox(0,0)[lb]{{\SetFigFont{10}{16.8}{\rmdefault}{\mddefault}{\updefault}$m+\bar{\ell}$}}}
\put(1450,1575){\makebox(0,0)[lb]{{\SetFigFont{10}{16.8}{\rmdefault}{\mddefault}{\updefault}$m+\ell$}}}
\lb{fgs}\end{picture}
}
\vspace*{-20 mm}
%\begin{rightline}
{\footnotesize

\hfill Figure 2.4: A ground state  configuration\\[-7mm]

\hfill   of the fusion SOS model \hspace*{5mm}
}
%\end{rightline}
\newpage \noindent The ground state indicated in Figure 2.4 with height 
variable $m+\ell$ on a specified reference site is labelled 
by the pair $(m,\ell)$.

\subsection{The Vertex-Face correspondence }
In order to solve the eight-vertex model by the Bethe ansatz, Baxter discovered the celebrated 
identity referred to as the vertex-face correspondence \cite{Bax73aII}. This
correspondence was later 
generalised to higher fusion level $k$ in \cite{DJMO86}.

\subsubsection{The simple $k=1$ case}
Let us consider the following vector (Figure 2.5 $(a)$)
\bea
&&\psi(u)^a_b=\psi_+(u)^a_b\ v_+ + \psi_-(u)^a_b\ v_-,\qquad \\
&&{\psi}_+(u)^a_b=\vtf{0}{\frac{(a-b)u+a}{2r}}{\frac{\tau}{2}},\qquad \psi_-(u)^a_b=\vtf{3}{\frac{(a-b)u+a}{2r}}{\frac{\tau}{2}} \lb{intertwinvec}
\ena 
with $|a-b|=1$. Baxter showed the following identity (Figure 2.6 $(a)$).
\bea
\sum_{\vep_1',\vep_2'}
R(u-v)_{\vep_1 \vep_2}^{\vep_1' \vep_2'}\ 
\psi_{\vep_1'}(u)_{b}^{a}
\psi_{\vep_2'}(v)_{c}^{b}
=\sum_{b' \in {\mathbb{Z}}}
\psi_{\vep_2}(v)_{b'}^{a}
\psi_{\vep_1}(u)_{c}^{b'}
W\left(\left.
\begin{array}{cc}
a&b\\
b'&c
\end{array}\right|u-v\right).\lb{vertexface}
\ena
We hence call $\psi(u)^a_b$ the intertwining vector.
This identity is a key formula throughout this paper. 
 One should note that Baxter's original 
intertwining vector intertwines the eight-vertex model in the disordered regime with the SOS model 
in regime III \cite{Bax73aII,DJKMO88}.
In order to consider
the eight-vertex model in
 the principal regime, we have derived the eight-vertex $R$-matrix
 \eqref{evR},  the SOS face weight \eqref{face} and the intertwining
 vector \eqref{intertwinvec} from those of Baxter\cite{Bax73aII}
 (which are the same as those used in
\cite{DJKMO88}) by using Jacobi's
 imaginary transformation.  

\vspace{5mm}
\begin{center}
\setlength{\unitlength}{0.01mm}%
\begingroup\makeatletter\ifx\SetFigFont\undefined%
\gdef\SetFigFont#1#2#3#4#5{%
  \reset@font\fontsize{#1}{#2pt}%
  \fontfamily{#3}\fontseries{#4}\fontshape{#5}%
  \selectfont}%
\fi\endgroup%
\begin{picture}(9312,2889)(301,-2815)
\put(7651,-2761){\makebox(0,0)[lb]{\smash{\SetFigFont{12}{14.4}{\rmdefault}{\mddefault}{\updefault}{\color[rgb]{0,0,0}$(b)$}%
}}}
\thinlines
{\color[rgb]{0,0,0}\put(3001,-361){\vector( 0,-1){1200}}
}%
{\color[rgb]{0,0,0}\put(8401,-1561){\line( 1, 0){1200}}
}%
{\color[rgb]{0,0,0}\put(9001,-361){\vector( 0,-1){600}}
}%
{\color[rgb]{0,0,0}\put(9001,-961){\line( 0,-1){600}}
}%
\put(2101,-61){\makebox(0,0)[lb]{\smash{\SetFigFont{12}{14.4}{\rmdefault}{\mddefault}{\updefault}{\color[rgb]{0,0,0}$a$}%
}}}
\put(3451,-61){\makebox(0,0)[lb]{\smash{\SetFigFont{12}{14.4}{\rmdefault}{\mddefault}{\updefault}{\color[rgb]{0,0,0}$b$}%
}}}
\put(3226,-1561){\makebox(0,0)[lb]{\smash{\SetFigFont{12}{14.4}{\rmdefault}{\mddefault}{\updefault}{\color[rgb]{0,0,0}$\vep$}%
}}}
\put(2901,-1861){\makebox(0,0)[lb]{\smash{\SetFigFont{12}{14.4}{\rmdefault}{\mddefault}{\updefault}{\color[rgb]{0,0,0}$u$}%
}}}
\put(201,-961){\makebox(0,0)[lb]{\smash{\SetFigFont{12}{14.4}{\rmdefault}{\mddefault}{\updefault}{\color[rgb]{0,0,0}$\psi_\vep(u)^a_b$}%
}}}
\put(1501,-961){\makebox(0,0)[lb]{\smash{\SetFigFont{12}{14.4}{\rmdefault}{\mddefault}{\updefault}{\color[rgb]{0,0,0}$=$}%
}}}
\put(1801,-2761){\makebox(0,0)[lb]{\smash{\SetFigFont{12}{14.4}{\rmdefault}{\mddefault}{\updefault}{\color[rgb]{0,0,0}$(a)$}%
}}}
\put(7501,-961){\makebox(0,0)[lb]{\smash{\SetFigFont{12}{14.4}{\rmdefault}{\mddefault}{\updefault}{\color[rgb]{0,0,0}$=$}%
}}}
\put(8901,-61){\makebox(0,0)[lb]{\smash{\SetFigFont{12}{14.4}{\rmdefault}{\mddefault}{\updefault}{\color[rgb]{0,0,0}$u$}%
}}}
\put(9151,-661){\makebox(0,0)[lb]{\smash{\SetFigFont{12}{14.4}{\rmdefault}{\mddefault}{\updefault}{\color[rgb]{0,0,0}$\vep$}%
}}}
\put(8326,-1861){\makebox(0,0)[lb]{\smash{\SetFigFont{12}{14.4}{\rmdefault}{\mddefault}{\updefault}{\color[rgb]{0,0,0}$a$}%
}}}
\put(9601,-1861){\makebox(0,0)[lb]{\smash{\SetFigFont{12}{14.4}{\rmdefault}{\mddefault}{\updefault}{\color[rgb]{0,0,0}$b$}%
}}}
\put(6001,-961){\makebox(0,0)[lb]{\smash{\SetFigFont{12}{14.4}{\rmdefault}{\mddefault}{\updefault}{\color[rgb]{0,0,0}$\psi^*_\vep(u)_a^b$}%
}}}
{\color[rgb]{0,0,0}\put(2401,-361){\line( 1, 0){1200}}
}%
\end{picture}

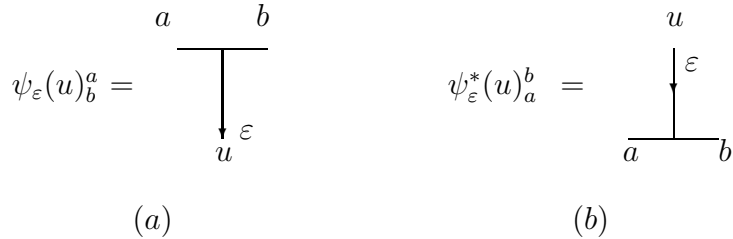
{\footnotesize Figure 2.5: \lb{ITV} $(a)$ The intertwining vector ; $(b)$ the dual intertwining vector}
\end{center}

\vspace{5mm}
\begin{center}
\setlength{\unitlength}{0.012mm}%
\begingroup\makeatletter\ifx\SetFigFont\undefined%
\gdef\SetFigFont#1#2#3#4#5{%
  \reset@font\fontsize{#1}{#2pt}%
  \fontfamily{#3}\fontseries{#4}\fontshape{#5}%
  \selectfont}%
\fi\endgroup%
\begin{picture}(11850,3369)(451,-3190)
\put(9001,-3136){\makebox(0,0)[lb]{\smash{\SetFigFont{12}{14.4}{\rmdefault}{\mddefault}{\updefault}{\color[rgb]{0,0,0}$(b)$}%
}}}
\thinlines
{\color[rgb]{0,0,0}\put(1801,-361){\vector( 0,-1){1500}}
}%
{\color[rgb]{0,0,0}\put(2401,-961){\vector(-1, 0){1500}}
}%
{\color[rgb]{0,0,0}\put(8401,-1561){\line(-1, 0){1200}}
\put(7201,-1561){\line( 0, 1){1200}}
}%
{\color[rgb]{0,0,0}\put(10501,-361){\line( 0,-1){1200}}
\put(10501,-1561){\line( 1, 0){1200}}
\put(11701,-1561){\line( 0, 1){1200}}
\put(11701,-361){\line(-1, 0){1200}}
}%
{\color[rgb]{0,0,0}\put(7801,-136){\vector( 0,-1){1500}}
}%
{\color[rgb]{0,0,0}\put(8701,-961){\vector(-1, 0){1500}}
}%
{\color[rgb]{0,0,0}\put(11101,-61){\vector( 0,-1){300}}
}%
{\color[rgb]{0,0,0}\put(4501,-361){\line( 0,-1){1200}}
\put(4501,-1561){\line( 1, 0){1200}}
\put(5701,-1561){\line( 0, 1){1200}}
\put(5701,-361){\line(-1, 0){1200}}
}%
{\color[rgb]{0,0,0}\put(4501,-961){\vector(-1, 0){300}}
}%
{\color[rgb]{0,0,0}\put(5101,-1561){\vector( 0,-1){300}}
}%
{\color[rgb]{0,0,0}\put(12001,-961){\vector(-1, 0){300}}
}%
{\color[rgb]{0,0,0}\put(4801,-361){\line(-1,-1){300}}
}%
{\color[rgb]{0,0,0}\put(10801,-361){\line(-1,-1){300}}
}%
\put(3301,-961){\makebox(0,0)[lb]{\smash{\SetFigFont{12}{14.4}{\rmdefault}{\mddefault}{\updefault}{\color[rgb]{0,0,0}$=$ }%
}}}
\put(1651,-2311){\makebox(0,0)[lb]{\smash{\SetFigFont{12}{14.4}{\rmdefault}{\mddefault}{\updefault}{\color[rgb]{0,0,0}$u$}%
}}}
\put(1951,-1861){\makebox(0,0)[lb]{\smash{\SetFigFont{12}{14.4}{\rmdefault}{\mddefault}{\updefault}{\color[rgb]{0,0,0}$\vep_1$}%
}}}
\put(901,-1261){\makebox(0,0)[lb]{\smash{\SetFigFont{12}{14.4}{\rmdefault}{\mddefault}{\updefault}{\color[rgb]{0,0,0}$\vep_2$}%
}}}
\put(451,-961){\makebox(0,0)[lb]{\smash{\SetFigFont{12}{14.4}{\rmdefault}{\mddefault}{\updefault}{\color[rgb]{0,0,0}$v$}%
}}}
\put(901,-211){\makebox(0,0)[lb]{\smash{\SetFigFont{12}{14.4}{\rmdefault}{\mddefault}{\updefault}{\color[rgb]{0,0,0}$a$}%
}}}
\put(2251,-211){\makebox(0,0)[lb]{\smash{\SetFigFont{12}{14.4}{\rmdefault}{\mddefault}{\updefault}{\color[rgb]{0,0,0}$b$}%
}}}
\put(2551,-1561){\makebox(0,0)[lb]{\smash{\SetFigFont{12}{14.4}{\rmdefault}{\mddefault}{\updefault}{\color[rgb]{0,0,0}$c$}%
}}}
\put(4351,-211){\makebox(0,0)[lb]{\smash{\SetFigFont{12}{14.4}{\rmdefault}{\mddefault}{\updefault}{\color[rgb]{0,0,0}$a$}%
}}}
\put(5551,-211){\makebox(0,0)[lb]{\smash{\SetFigFont{12}{14.4}{\rmdefault}{\mddefault}{\updefault}{\color[rgb]{0,0,0}$b$}%
}}}
\put(5851,-1561){\makebox(0,0)[lb]{\smash{\SetFigFont{12}{14.4}{\rmdefault}{\mddefault}{\updefault}{\color[rgb]{0,0,0}$c$}%
}}}
\put(5026,-2311){\makebox(0,0)[lb]{\smash{\SetFigFont{12}{14.4}{\rmdefault}{\mddefault}{\updefault}{\color[rgb]{0,0,0}$u$}%
}}}
\put(3751,-961){\makebox(0,0)[lb]{\smash{\SetFigFont{12}{14.4}{\rmdefault}{\mddefault}{\updefault}{\color[rgb]{0,0,0}$v$}%
}}}
\put(5251,-1936){\makebox(0,0)[lb]{\smash{\SetFigFont{12}{14.4}{\rmdefault}{\mddefault}{\updefault}{\color[rgb]{0,0,0}$\vep_1$}%
}}}
\put(3901,-1261){\makebox(0,0)[lb]{\smash{\SetFigFont{12}{14.4}{\rmdefault}{\mddefault}{\updefault}{\color[rgb]{0,0,0}$\vep_2$}%
}}}
\put(8326,-1936){\makebox(0,0)[lb]{\smash{\SetFigFont{12}{14.4}{\rmdefault}{\mddefault}{\updefault}{\color[rgb]{0,0,0}$a$}%
}}}
\put(6901,-1936){\makebox(0,0)[lb]{\smash{\SetFigFont{12}{14.4}{\rmdefault}{\mddefault}{\updefault}{\color[rgb]{0,0,0}$b$}%
}}}
\put(6901,-211){\makebox(0,0)[lb]{\smash{\SetFigFont{12}{14.4}{\rmdefault}{\mddefault}{\updefault}{\color[rgb]{0,0,0}$c$}%
}}}
\put(7951,-286){\makebox(0,0)[lb]{\smash{\SetFigFont{12}{14.4}{\rmdefault}{\mddefault}{\updefault}{\color[rgb]{0,0,0}$\vep_1$}%
}}}
\put(8251,-1186){\makebox(0,0)[lb]{\smash{\SetFigFont{12}{14.4}{\rmdefault}{\mddefault}{\updefault}{\color[rgb]{0,0,0}$\vep_2$}%
}}}
\put(7726, 89){\makebox(0,0)[lb]{\smash{\SetFigFont{12}{14.4}{\rmdefault}{\mddefault}{\updefault}{\color[rgb]{0,0,0}$u$}%
}}}
\put(8926,-1036){\makebox(0,0)[lb]{\smash{\SetFigFont{12}{14.4}{\rmdefault}{\mddefault}{\updefault}{\color[rgb]{0,0,0}$v$}%
}}}
\put(12301,-1036){\makebox(0,0)[lb]{\smash{\SetFigFont{12}{14.4}{\rmdefault}{\mddefault}{\updefault}{\color[rgb]{0,0,0}$v$}%
}}}
\put(11026, 89){\makebox(0,0)[lb]{\smash{\SetFigFont{12}{14.4}{\rmdefault}{\mddefault}{\updefault}{\color[rgb]{0,0,0}$u$}%
}}}
\put(11326,-136){\makebox(0,0)[lb]{\smash{\SetFigFont{12}{14.4}{\rmdefault}{\mddefault}{\updefault}{\color[rgb]{0,0,0}$\vep_1$}%
}}}
\put(11851,-1186){\makebox(0,0)[lb]{\smash{\SetFigFont{12}{14.4}{\rmdefault}{\mddefault}{\updefault}{\color[rgb]{0,0,0}$\vep_2$}%
}}}
\put(10201,-1936){\makebox(0,0)[lb]{\smash{\SetFigFont{12}{14.4}{\rmdefault}{\mddefault}{\updefault}{\color[rgb]{0,0,0}$b$}%
}}}
\put(10201,-136){\makebox(0,0)[lb]{\smash{\SetFigFont{12}{14.4}{\rmdefault}{\mddefault}{\updefault}{\color[rgb]{0,0,0}$c$}%
}}}
\put(11626,-1936){\makebox(0,0)[lb]{\smash{\SetFigFont{12}{14.4}{\rmdefault}{\mddefault}{\updefault}{\color[rgb]{0,0,0}$a$}%
}}}
\put(9601,-1036){\makebox(0,0)[lb]{\smash{\SetFigFont{12}{14.4}{\rmdefault}{\mddefault}{\updefault}{\color[rgb]{0,0,0}$=$}%
}}}
\put(3001,-3136){\makebox(0,0)[lb]{\smash{\SetFigFont{12}{14.4}{\rmdefault}{\mddefault}{\updefault}{\color[rgb]{0,0,0}$(a)$}%
}}}
{\color[rgb]{0,0,0}\put(1201,-361){\line( 1, 0){1200}}
\put(2401,-361){\line( 0,-1){1200}}
}%
\end{picture}

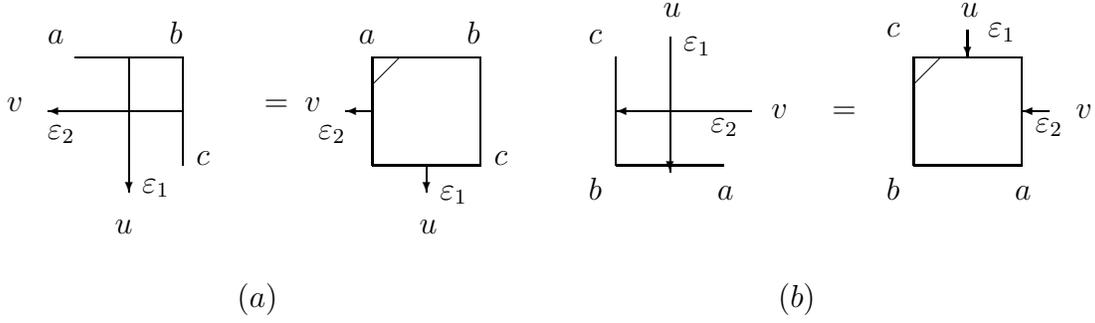
{\footnotesize Figure 2.6:\lb{VFC} The vertex-face correspondence: $(a)$ via the intertwining vector ; $(b)$ via the dual intertwining vector}
\end{center}

In addition to the intertwining vector, it is necessary to introduce its dual counterpart and a second intertwining vector.
The dual intertwining vector $\psi^*(u)^a_b$ (Figure 2.5 $(b)$) is defined by 
\bea
&&\psi^*(u)^a_b\ v_{\vep}= \psi^*_\vep(u)^a_b,\qquad 
\psi^*_\vep(u)^a_b=-\vep\frac{a-b}{2[b][u]}C^2\ \psi_{-\vep}(u-1)^a_b, \lb{dualintvec}
\ena
whereas the second intertwining vector $\psi'(u)^a_b\ (b=a\pm a)$ is given by
\be
&&\psi'(u)^a_b=\sum_{\vep=\pm}\psi'_\vep(u)^a_b\ v_\vep,\qquad \psi'_\vep(u)^a_b=\frac{[u][a]}{[u-1][b]}\psi_\vep(u-2)^a_b 
\en
with $|a-b|=1$ in both cases.
Then by direct calculation, one can verify the following inversion relations
(Figure 2.7)
\bea
&&\sum_{\vep=\pm}\psi_\vep^*(u)^a_b\psi_\vep(u)^b_c=\delta_{a,c},\lb{inversiona}\\
&&\sum_{a=b\pm1}\psi_{\vep'}^*(u)^a_b\psi_\vep(u)^b_a=\delta_{\vep',\vep},\lb{inversionb}\\
&&\sum_{\vep=\pm}\psi^*_\vep(u)^a_b\psi'_\vep(u)^c_a=\delta_{b,c},\\
&&\sum_{b=a\pm 1}\psi^*_\vep(u)^a_b\psi'_{\vep'}(u)^b_a=\delta_{\vep,\vep'}.
\lb{inversionbp}
\ena
These inversion properties are the reason that we 
call $\psi^*(u)^a_b$ the dual intertwining vector. 

\vspace{5mm}
\begin{center}
\setlength{\unitlength}{0.012mm}%
\begingroup\makeatletter\ifx\SetFigFont\undefined%
\gdef\SetFigFont#1#2#3#4#5{%
  \reset@font\fontsize{#1}{#2pt}%
  \fontfamily{#3}\fontseries{#4}\fontshape{#5}%
  \selectfont}%
\fi\endgroup%
\begin{picture}(8025,2739)(1051,-2815)
\put(2326,-1861){\makebox(0,0)[lb]{\smash{\SetFigFont{12}{14.4}{\rmdefault}{\mddefault}{\updefault}{\color[rgb]{0,0,0}$a$}%
}}}
\thinlines
{\color[rgb]{0,0,0}\put(1201,-361){\line( 1, 0){1200}}
}%
{\color[rgb]{0,0,0}\put(1801,-361){\vector( 0,-1){600}}
}%
{\color[rgb]{0,0,0}\put(1801,-961){\line( 0,-1){600}}
}%
\dottedline[$\circle*{15}$]{100}(1201,-361)(1201,-1561)
%{\color[rgb]{0,0,0}\multiput(1201,-361)(0.00000,-9.02256){134}{\makebox(1.6667,11.6667)
%{\SetFigFont{5}{6}{\rmdefault}{\mddefault}{\updefault}.}}}%
{\color[rgb]{0,0,0}\put(1201,-1561){\line( 1, 0){1200}}
}%
{\color[rgb]{0,0,0}\put(6301,-361){\line( 1, 0){1270}}
}%
{\color[rgb]{0,0,0}\put(6301,-1561){\line( 1, 0){1270}}
}%
{\color[rgb]{0,0,0}\put(6901,-361){\vector( 0,-1){450}}
}%
{\color[rgb]{0,0,0}\put(6901,-1111){\vector( 0,-1){450}}
}%
\dottedline[$\circle*{15}$]{100}(6301,-361)(6301,-1561)
%{\color[rgb]{0,0,0}\multiput(6301,-361)(0.00000,-9.02256){134}{\makebox(1.6667,11.6667)
%{\SetFigFont{5}{6}{\rmdefault}{\mddefault}{\updefault}.}}}%
\put(3001,-961){\makebox(0,0)[lb]{\smash{\SetFigFont{12}{14.4}{\rmdefault}{\mddefault}{\updefault}{\color[rgb]{0,0,0}$=$}%
}}}
\put(3601,-961){\makebox(0,0)[lb]{\smash{\SetFigFont{12}{14.4}{\rmdefault}{\mddefault}{\updefault}{\color[rgb]{0,0,0}$\delta_{a,c}$}%
}}}
\put(1051,-211){\makebox(0,0)[lb]{\smash{\SetFigFont{12}{14.4}{\rmdefault}{\mddefault}{\updefault}{\color[rgb]{0,0,0}$b$}%
}}}
\put(1051,-1861){\makebox(0,0)[lb]{\smash{\SetFigFont{12}{14.4}{\rmdefault}{\mddefault}{\updefault}{\color[rgb]{0,0,0}$b$}%
}}}
\put(2326,-211){\makebox(0,0)[lb]{\smash{\SetFigFont{12}{14.4}{\rmdefault}{\mddefault}{\updefault}{\color[rgb]{0,0,0}$c$}%
}}}
\put(7051,-736){\makebox(0,0)[lb]{\smash{\SetFigFont{12}{14.4}{\rmdefault}{\mddefault}{\updefault}{\color[rgb]{0,0,0}$\vep$}%
}}}
\put(7051,-1261){\makebox(0,0)[lb]{\smash{\SetFigFont{12}{14.4}{\rmdefault}{\mddefault}{\updefault}{\color[rgb]{0,0,0}$\vep'$}%
}}}
\put(2401,-2761){\makebox(0,0)[lb]{\smash{\SetFigFont{12}{14.4}{\rmdefault}{\mddefault}{\updefault}{\color[rgb]{0,0,0}$(a)$}%
}}}
\put(7201,-2761){\makebox(0,0)[lb]{\smash{\SetFigFont{12}{14.4}{\rmdefault}{\mddefault}{\updefault}{\color[rgb]{0,0,0}$(b)$}%
}}}
\put(8551,-1036){\makebox(0,0)[lb]{\smash{\SetFigFont{12}{14.4}{\rmdefault}{\mddefault}{\updefault}{\color[rgb]{0,0,0}$=$}%
}}}
\put(9076,-1036){\makebox(0,0)[lb]{\smash{\SetFigFont{12}{14.4}{\rmdefault}{\mddefault}{\updefault}{\color[rgb]{0,0,0}$\delta_{\vep,\vep'}$}%
}}}
\put(6226,-1861){\makebox(0,0)[lb]{\smash{\SetFigFont{12}{14.4}{\rmdefault}{\mddefault}{\updefault}{\color[rgb]{0,0,0}$b$}%
}}}
\put(6226,-211){\makebox(0,0)[lb]{\smash{\SetFigFont{12}{14.4}{\rmdefault}{\mddefault}{\updefault}{\color[rgb]{0,0,0}$b$}%
}}}
{\color[rgb]{0,0,0}\put(7501,-961){\oval(1200,1200)[br]}
\put(7501,-961){\oval(1200,1200)[tr]}
}%
\end{picture}

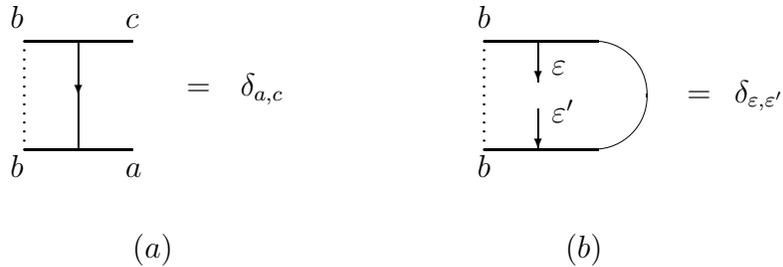
{\footnotesize Figure 2.7:\lb{ITVINV} The inversion relations between the intertwining vector and its dual.}
\end{center}
It then follows from the crossing symmetry properties of $R$ and $W$
that the following vertex-face correspondence holds:
\bea
\sum_{\vep_1',\vep_2'}
R(u-v)^{\vep_1 \vep_2}_{\vep_1' \vep_2'}\ 
\psi^*_{\vep_1'}(u)_{b}^{a}
\psi^*_{\vep_2'}(v)_{c}^{b}
=\sum_{s \in {\mathbb{Z}}}
\psi^*_{\vep_2}(v)_{b'}^{a}
\psi^*_{\vep_1}(u)_{c}^{b'}
W\left(\left.
\begin{array}{cc}
c&b'\\
b&a
\end{array}\right|u-v\right).\lb{vertexfacedual} 
\ena
This relation is represented by Figure 2.6 $(b)$. 

\subsubsection{The general $k$ case}
Let us now discuss the fusion of the vertex-face relationship \eqref{vertexface}. 
We define fused intertwining vectors by
 \bea
 &&\psi^{(k)}(u)^a_b=\Pi_{1\cdots k}\ \psi(u+k-1)^a_{c_1}\otimes \psi(u+k-2)^{c_1}_{c_2} \otimes \cdots 
\otimes \psi(u)^{c_{k-1}}_{b},\lb{kpsi}
 \ena
 The RHS is independent of the choice of $c_1,..,c_{k-1}$ provided
 $|a-c_1|=|c_1-c_2|=\cdots=|c_{k-1}-b|=1$. The components of $\psi^{(k)}(u)^a_b$ are given by the following formula.
\be
\psi^{(k)}(u)^a_b&=&\sum_{\ep\in\{-k,-k+2,..,k\}}\psi^{(k)}_\ep(u)^a_b\ v^{(k)}_\ep,\\
\psi^{(k)}_\ep(u)^a_b&=&
\sum_{\vep_1,\cdots,\vep_{k}\atop {\ep_1+\ep_2+\cdots+\ep_k=\ep}}
\psi_{\vep_1}(u+k-1)^a_{c_1}
\psi_{\vep_2}(u+k-2)^{c_1}_{c_2} \cdots  \psi_{\ep_k}(u)^{c_{k-1}}_{b}.\lb{compokpsi}
\en

From \eqref{kkfusionR},
 \eqref{kkfusionW},  and \eqref{vertexface}, it follows that
 $\psi^{(k)}(u)^a_b$ satisfies the $k\times k$ fusion vertex-face 
 correspondence 
relations with respect to $R^{(k,k)}$ and $W^{(k,k)}$. That is, we have
 \bea
 \hspace*{-20mm}&&\sum_{\ep_1',\ep_2'}
 R^{(k,k)}(u-v)_{\ep_1 \ep_2}^{\ep_1' \ep_2'}\ 
 \psi^{(k)}_{\ep_1'}(u)_{b}^{a}
 \psi^{(k)}_{\ep_2'}(v)_{c}^{b}
 =\sum_{b' \in {\mathbb{Z}}}
 \psi^{(k)}_{\ep_2}(v)_{b'}^{a}
 \psi^{(k)}_{\ep_1}(u)_{c}^{b'}
 W^{(k,k)}\left(\left.
 \begin{array}{cc}
 a&b\\
 b'&c
 \end{array}\right|u-v\right).
% \nn\\&&
 \lb{fusionvertexface}
 \ena
Similarly, we fuse the second intertwining vector $\psi'(u)^a_b$ as follows. 
\be
\psi^{'(k)}(u)^a_b&=&\Pi_{1\cdots k}\ \psi'(u+k-1)^a_{c_1}\otimes \psi'(u+k-2)^{c_1}_{c_2} \otimes \cdots 
\otimes \psi'(u)^{c_{k-1}}_{b}.\lb{kpsip}
\en
Then we find that the components of $\psi^{'(k)}(u)^a_b$ are given by
\be
\psi^{'(k)}(u)^a_b&=&\sum_{\ep\in\{-k,-k+2,..,k\}}\psi^{'(k)}_\ep(u)^a_b\ v^{(k)}_\ep,\\
\psi^{'(k)}_\ep(u)^a_b&=&\frac{[u+k-1][a]}{[u-1][b]}\psi^{(k)}_\ep(u-2)^a_b.\lb{compopsip}
\en

In addition, we define the fusion of the  dual intertwining vector in the 
following way.
 \begin{eqnarray}
\psi^{*(k)}(u)^b_a=\sum_{c_1,..,c_{k-1}}\psi^*(u+k-1)_a^{c_1}\otimes \psi^*(u+k-2)^{c_2}_{c_1} 
\otimes \cdots \otimes \psi^*(u)_{c_{k-1}}^{b}.
\lb{kpsis}
 \end{eqnarray}
 with the property 
 \be
 &&\Pi_{1 \cdots k}\ \psi^{*(k)}(u)^b_a = \psi^{*(k)}(u)^b_a\ \Pi_{1 \cdots k}. 
 \en
Written out in component form, the last relation indicates that the RHS of   
 \be
 \psi^{*(k)}_{\ep}(u)_a^b=
 \sum_{c_1,..,c_{k-1}} \psi^*_{\vep_1}(u+k-1)_a^{c_1}\psi^*_{\vep_2}(u+k-2)^{c_2}_{c_1} 
\cdots \psi^*_{\vep_k}(u)_{c_{k-1}}^{b}
 \en
 is independent of the choice of $\vep_1,\cdots,
 \vep_k$  provided $\ep=\vep_1+\cdots+\vep_k$.

As above, it follows immediately from \mref{vertexfacedual},
\mref{kkfusionR} and \mref{kkfusionW}, that we have
 \bea
\hspace*{-30mm} &&\sum_{\ep_1',\ep_2'}
 R^{(k,k)}(u-v)_{\ep_1 \ep_2}^{\ep_1' \ep_2'}\ 
 \psi^{*(k)}_{\ep_1'}(u)_{b}^{a}
 \psi^{*(k)}_{\ep_2'}(v)_{c}^{b}
 =\sum_{b' \in {\mathbb{Z}}}
 \psi^{*(k)}_{\ep_2}(v)_{b'}^{a}
 \psi^{*(k)}_{\ep_1}(u)_{c}^{b'}
 W^{(k,k)}\left(\left.
 \begin{array}{cc}
 c&b'\\
 b&a
 \end{array}\right|u-v\right).
% \nn\\ &&
\lb{fusionvertexfacedual}
 \ena
 It is worth noting that this formula is also obtained from \eqref{fusionvertexface} 
 by using the crossing symmetry properties \eqref{qconj} and \eqref{crossing}. 
 In fact, the fused dual intertwining vector is related to the intertwining vector 
 as follows.
 \bea
 &&\psi^{*(k)}_\vep(u)^a_b=C^{(k)}(u)G^{(k)}_{a,b}\sum_{\vep'}Q^{\vep'}_\vep\psi^{(k)}_{\vep'}(u-1)^a_b.\lb{crosspsi}
 \ena 
Here $C^{(k)}(u)$ is a certain normalisation function. The $k=1$ case is given in \eqref{dualintvec}, whereas the $k=2$ case in \cite{Konno04}. 

Finally,  using \eqref{inversiona} - \eqref{inversionbp}, it is easy to verify the following inversion relations.
% \begin{prop}
 \bea
 \sum_{\ep\in\{-k,-k+2,\cdots,k\}}\psi^{*(k)}_\ep(u)^a_b\psi^{(k)}_\ep(u)^b_c&=&\delta_{a,c},\lb{fusioninv1}
 \\
 \sum_{a\in b+\{-k,-k+2,..,k\}}\psi^{*(k)}_{\ep'}(u)^a_b\psi^{(k)}_\ep(u)^b_a&=&\delta_{\ep', \ep},\lb{fusioninv2}\\
\sum_{\ep\in\{-k,-k+2,..,k\}}\psi^{*(k)}_\ep(u)^a_b\psi^{'(k)}_\ep(u)^c_a
&=&\delta_{b,c},\lb{fusioninv1p}\\
\sum_{b\in a+\{-k,-k+2,..,k\} }\psi^{*(k)}_{\ep'}(u)^a_b
\psi^{'(k)}_{\ep}(u)^b_a
&=&\delta_{\ep',\ep}.\lb{fusioninv2p}
 \ena
% \end{prop}
\subsubsection{The $L$-Matrix}
In the next section, we will make use of the `$L$-matrix',
defined in terms of the intertwiner and dual intertwiner by
 \bea
 \lmatk{a}{b}{c}{d}{u}{\vep}= \psi^{*(k)}_{\vep}(u)^d_c \, \psi^{(k)}_{\vep}(u)^a_b.
 \lb{Lmate}
 \ena
The graphical representation is given in Figure 2.8.

\setlength{\unitlength}{0.0004in}
\begingroup\makeatletter\ifx\SetFigFont\undefined%
\gdef\SetFigFont#1#2#3#4#5{%
  \reset@font\fontsize{#1}{#2pt}%
  \fontfamily{#3}\fontseries{#4}\fontshape{#5}%
  \selectfont}%
\fi\endgroup%
{\renewcommand{\dashlinestretch}{30}
\begin{picture}(2073,2430)(-5000,600)
\thinlines
\texture{55888888 88555555 5522a222 a2555555 55888888 88555555 552a2a2a 2a555555 
        55888888 88555555 55a222a2 22555555 55888888 88555555 552a2a2a 2a555555 
        55888888 88555555 5522a222 a2555555 55888888 88555555 552a2a2a 2a555555 
        55888888 88555555 55a222a2 22555555 55888888 88555555 552a2a2a 2a555555 }
\path(150,2100)(1950,2100)
\path(150,2100)(1950,2100)
\path(150,300)(1950,300)
\path(150,300)(1950,300)
\path(1050,2100)(1050,1050)
\path(1050,2100)(1050,1050)
\path(990.000,1290.000)(1050.000,1050.000)(1110.000,1290.000)
\path(1050,1050)(1050,300)
\path(1050,1050)(1050,300)
\put(0,2250){\makebox(0,0)[lb]{{\SetFigFont{12}{16.8}{\rmdefault}{\mddefault}{\updefault}$a$}}}
\put(1950,2250){\makebox(0,0)[lb]{{\SetFigFont{12}{16.8}{\rmdefault}{\mddefault}{\updefault}$b$}}}
\put(0,0){\makebox(0,0)[lb]{{\SetFigFont{12}{16.8}{\rmdefault}{\mddefault}{\updefault}$c$}}}
\put(1950,0){\makebox(0,0)[lb]{{\SetFigFont{12}{16.8}{\rmdefault}{\mddefault}{\updefault}$d$}}}
\put(1200,1500){\makebox(0,0)[lb]{{\SetFigFont{12}{16.8}{\rmdefault}{\mddefault}{\updefault}$u$}}}
\put(750,900){\makebox(0,0)[lb]{{\SetFigFont{12}{16.8}{\rmdefault}{\mddefault}{\updefault}$\vep$}}}
\end{picture}
}

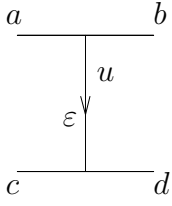
\vspace*{4mm}\begin{center}{\ft Figure 2.8: The graphical representation of the $L-$matrix}
\end{center}

\nin It is also useful to define the matrix
\bea 
\lmatsk{a}{b}{c}{d}{u}=\sli_{\vep\in\{-k,-k+2,\cdots,k\}} \lmatk{a}{b}{c}{d}{u}{\vep}.\lb{Lmat}
\ena

If we restrict $-1<u+\frac{k-1}{2}<0$, $m\geq 1+\frac{k}{2}$, and
choose $\ell\in\{0,1,\cdots,k\}$, we
find that the $L-$matrix with $a,b,c,d$ specified as follows
\ben 
 \lmatk{m+\ell}{m+\bar{\ell}}{m+\ell}{m+\bar{\ell}}{u}{\vep},
 \een
has a maximum absolute value for the choice $\ep_\ell=k-2\ell$. Thus, maximal
weight $L-$matrix configurations are of the form shown in Figure 2.9.

\setlength{\unitlength}{0.0004in}
\begingroup\makeatletter\ifx\SetFigFont\undefined%
\gdef\SetFigFont#1#2#3#4#5{%
  \reset@font\fontsize{#1}{#2pt}%
  \fontfamily{#3}\fontseries{#4}\fontshape{#5}%
  \selectfont}%
\fi\endgroup%
{\renewcommand{\dashlinestretch}{30}
\begin{picture}(2073,2430)(-5000,600)
\thinlines
\texture{55888888 88555555 5522a222 a2555555 55888888 88555555 552a2a2a 2a555555 
        55888888 88555555 55a222a2 22555555 55888888 88555555 552a2a2a 2a555555 
        55888888 88555555 5522a222 a2555555 55888888 88555555 552a2a2a 2a555555 
        55888888 88555555 55a222a2 22555555 55888888 88555555 552a2a2a 2a555555 }
\path(150,2100)(1950,2100)
\path(150,2100)(1950,2100)
\path(150,300)(1950,300)
\path(150,300)(1950,300)
\path(1050,2100)(1050,1050)
\path(1050,2100)(1050,1050)
\path(990.000,1290.000)(1050.000,1050.000)(1110.000,1290.000)
\path(1050,1050)(1050,300)
\path(1050,1050)(1050,300)
\put(-400,2250){\makebox(0,0)[lb]{{\SetFigFont{12}{16.8}{\rmdefault}{\mddefault}{\updefault}$m+\ell$}}}
\put(1950,2250){\makebox(0,0)[lb]{{\SetFigFont{12}{16.8}{\rmdefault}{\mddefault}{\updefault}$m+\bar{\ell}$}}}
\put(-400,-200){\makebox(0,0)[lb]{{\SetFigFont{12}{16.8}{\rmdefault}{\mddefault}{\updefault}$m+\ell$}}}
\put(1950,-200){\makebox(0,0)[lb]{{\SetFigFont{12}{16.8}{\rmdefault}{\mddefault}{\updefault}$m+\bar{\ell}$}}}
\put(1200,1500){\makebox(0,0)[lb]{{\SetFigFont{12}{16.8}{\rmdefault}{\mddefault}{\updefault}$u$}}}
\put(1150,900){\makebox(0,0)[lb]{{\SetFigFont{12}{16.8}{\rmdefault}{\mddefault}{\updefault}$\vep_\ell$}}}
\end{picture}
}

\vspace*{4mm}\begin{center}{\ft Figure 2.9: A maximal weight $L-$matrix configuration}
\end{center}

%%%%%%%%%%%%%%%%%%%%%%%%%%%%%%%%%%%%%%%%
\section{The Corner Transfer Matrix and Half Transfer Matrices}
%%%%%%%%%%%%%%%%%%%%%%%%%%%%%%%%%%%%%%%%
In this and the next section, we consider correlation functions of the fusion
eight-vertex models introduced above. We express and manipulate these
correlation functions using two ideas: the expression for correlation 
functions in terms of the corner transfer matrix (CTM) and half-transfer matrices 
(HTMs)\cite{JM,Fodal93}; and the vertex-face correspondence\cite{LaP98}. 
We here review the first part, i.e., the algebraic analysis approach 
to both the fusion eight-vertex and the SOS models.

\ssect{Fusion Eight-Vertex Models}
Correlation functions of the fusion vertex models correspond 
to the probabilities of the $\ep$ edge variables taking certain
values on some specified set of edges of a lattice. 
More concretely, consider the following dimension $(2L+N)\times 2L$ 
lattice on which the edge variables at the indicated sides have the values
 $\ep_1,\cdots\ep_N$ (where $N$ in Figure 3.1 as shown is actually $3$, and for simplicity
we assume that $L$ is even). 

\setlength{\unitlength}{0.0003in}
\begingroup\makeatletter\ifx\SetFigFont\undefined%
\gdef\SetFigFont#1#2#3#4#5{%
  \reset@font\fontsize{#1}{#2pt}%
  \fontfamily{#3}\fontseries{#4}\fontshape{#5}%
  \selectfont}%
\fi\endgroup%
{\renewcommand{\dashlinestretch}{30}
\begin{picture}(9944,7859)(-4000,1000)
\thinlines
\texture{55888888 88555555 5522a222 a2555555 55888888 88555555 552a2a2a 2a555555 
        55888888 88555555 55a222a2 22555555 55888888 88555555 552a2a2a 2a555555 
        55888888 88555555 5522a222 a2555555 55888888 88555555 552a2a2a 2a555555 
        55888888 88555555 55a222a2 22555555 55888888 88555555 552a2a2a 2a555555 }
\put(4972,4222){\shade\ellipse{100}{100}}
\put(4972,4222){\ellipse{100}{100}}
\put(4072,4222){\shade\ellipse{100}{100}}
\put(4072,4222){\ellipse{100}{100}}
\put(5872,4222){\shade\ellipse{100}{100}}
\put(5872,4222){\ellipse{100}{100}}
\path(1372,7822)(1372,1522)
\path(922,6472)(8122,6472)
\path(922,5572)(8122,5572)
\path(922,4672)(8122,4672)
\path(922,2872)(8122,2872)
\path(922,7372)(8122,7372)
\path(922,3772)(8122,3772)
\path(2272,7822)(2272,1522)
\path(3172,7822)(3172,1522)
\path(4072,7822)(4072,1522)
\path(4972,7822)(4972,1522)
\path(5872,7822)(5872,1522)
\path(6772,7822)(6772,1522)
\path(7672,7822)(7672,1522)
\path(922,1972)(8122,1972)
\path(8572,7822)(8572,1522)
\path(8572,7822)(8572,1522)
\path(8122,7372)(9022,7372)
\path(8122,7372)(9022,7372)
\path(8122,6472)(9022,6472)
\path(8122,6472)(9022,6472)
\path(8122,5572)(9022,5572)
\path(8122,5572)(9022,5572)
\path(8122,4672)(9022,4672)
\path(8122,4672)(9022,4672)
\path(8122,3772)(9022,3772)
\path(8122,3772)(9022,3772)
\path(8122,2872)(9022,2872)
\path(8122,2872)(9022,2872)
\path(8122,1972)(9022,1972)
\path(8122,1972)(9022,1972)
\path(922,1072)(9022,1072)
\path(922,1072)(9022,1072)
\path(1372,1522)(1372,622)
\path(1372,1522)(1372,622)
\path(1312.000,862.000)(1372.000,622.000)(1432.000,862.000)
\path(2272,1597)(2272,622)
\path(2272,1597)(2272,622)
\path(2212.000,862.000)(2272.000,622.000)(2332.000,862.000)
\path(3172,1522)(3172,622)
\path(3172,1522)(3172,622)
\path(3112.000,862.000)(3172.000,622.000)(3232.000,862.000)
\path(4072,1597)(4072,622)
\path(4072,1597)(4072,622)
\path(4012.000,862.000)(4072.000,622.000)(4132.000,862.000)
\path(4972,1522)(4972,622)
\path(4972,1522)(4972,622)
\path(4912.000,862.000)(4972.000,622.000)(5032.000,862.000)
\path(5872,1522)(5872,622)
\path(5872,1522)(5872,622)
\path(5812.000,862.000)(5872.000,622.000)(5932.000,862.000)
\path(6772,1522)(6772,622)
\path(6772,1522)(6772,622)
\path(6712.000,862.000)(6772.000,622.000)(6832.000,862.000)
\path(7672,1522)(7672,622)
\path(7672,1522)(7672,622)
\path(7612.000,862.000)(7672.000,622.000)(7732.000,862.000)
\path(8572,1522)(8572,622)
\path(8572,1522)(8572,622)
\path(8512.000,862.000)(8572.000,622.000)(8632.000,862.000)
\drawline(9472,22)(9472,22)
\path(9472,7822)(9472,622)
\path(472,7822)(472,622)
\path(412.000,862.000)(472.000,622.000)(532.000,862.000)
\path(9022,7372)(9922,7372)
\path(9022,6472)(9922,6472)
\path(9022,5572)(9922,5572)
\path(9022,4672)(9922,4672)
\path(9022,3772)(9922,3772)
\path(9022,2872)(9922,2872)
\path(9022,1972)(9922,1972)
\path(9022,1072)(9922,1072)
\path(922,7372)(22,7372)
\path(262.000,7432.000)(22.000,7372.000)(262.000,7312.000)
\path(922,6472)(22,6472)
\path(262.000,6532.000)(22.000,6472.000)(262.000,6412.000)
\path(922,5572)(22,5572)
\path(262.000,5632.000)(22.000,5572.000)(262.000,5512.000)
\path(922,4672)(22,4672)
\path(262.000,4732.000)(22.000,4672.000)(262.000,4612.000)
\path(922,3772)(22,3772)
\path(262.000,3832.000)(22.000,3772.000)(262.000,3712.000)
\path(922,2872)(22,2872)
\path(262.000,2932.000)(22.000,2872.000)(262.000,2812.000)
\path(922,1972)(97,1972)
\path(337.000,2032.000)(97.000,1972.000)(337.000,1912.000)
\path(922,1072)(22,1072)
\path(262.000,1132.000)(22.000,1072.000)(262.000,1012.000)
\put(3547,4147){\makebox(0,0)[lb]{{\SetFigFont{10}{16.8}{\rmdefault}{\mddefault}{\updefault}$\ep_1$}}}
\put(4522,4147){\makebox(0,0)[lb]{{\SetFigFont{10}{16.8}{\rmdefault}{\mddefault}{\updefault}$\ep_2$}}}
\put(5422,4147){\makebox(0,0)[lb]{{\SetFigFont{10}{16.8}{\rmdefault}{\mddefault}{\updefault}$\ep_3$}}}
\end{picture}
}

\vspace*{4mm}
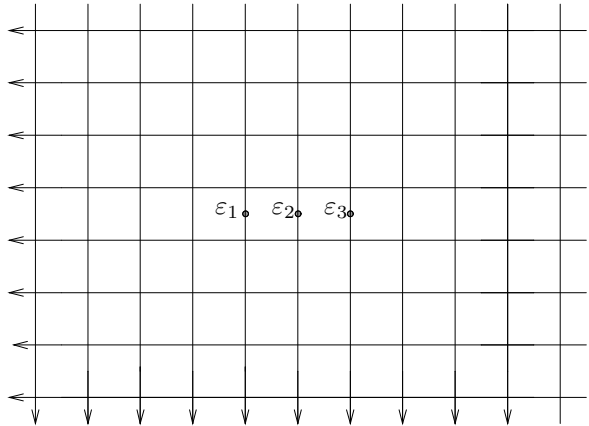
\begin{center}{\ft Figure 3.1: The restriction of edge variables associated with our correlation function}\end{center}

\nin The correlation function we consider specifies the probability
of such a configuration. It is the ratio of the weighted sum over such 
restricted configurations to the weighted sum over all configurations
(the latter sum being
the partition function). The total weight of any configuration is the
product of the local vertex Boltzmann weights.

The algebraic analysis approach of \cite{JM,Fodal93} gives a
way of computing such sums
in the infinite $L$ limit. To be more specific, it allows the computation of
this correlation function for the infinite-volume lattice in which 
sums are taken over edge variable configurations which are fixed 
to one of the ground state configurations of
Figure 2.2 beyond a finite, but arbitrarily large, 
distance from the centre of the lattice. We denote this 
correlation function by $P^{(\ell)}(\ep_1,\ep_2,\cdots,\ep_N)$,
where $\ell\in\{0,1,\cdots,k\}$ labels the chosen ground state configuration.

\subsubsection{The space of states}
The starting point is to replace the weighted sum by a trace
over a vector space $\cH^{(\ell)}$ representing a line running from the 
centre to the boundary of the infinite lattice. $\cH^{(\ell)}$ is called the
 space of states and  
defined in terms of basis vectors $v^{(k)}_{\vep}\ (\vep=-k,-k+2,..,k)$ by
\ben \cH^{(\ell)}
&=&\hb{Span}_{\C}\left\{ \cdots \ot  v^{(k)}_{\ep(1)} \ot
  v^{(k)}_{\ep(0)}|\ep(i)\in \{-k,-k+2,\cdots,k\},
\ep(i)=\bar{\ep}^{\,(\ell)}(i)\hb{ for } i\gg 0
\right\}\\
\bar{\ep}^{\,(\ell)}(i)&=& 
\begin{cases} 
2\ell-k & \hb{ for } i=0 \mod 2\\
 k-2\ell & \hb{ for }i=1 \mod 2
\end{cases}
\een
The correlation function $P^{(\ell)}(\ep_1,\ep_2,\cdots,\ep_N)$
is then represented in terms of a ratio of traces over
$\cH^{(\ell)}$ of CTMs and HTMs. 

\subsubsection{CTMs and the partition function}
These operators are best defined graphically. The North-West corner
transfer matrix $\cA(u):\cH^{(\ell)}\ra \cH^{(\ell)}$ is represented in
Figure 3.2, where $\cdots$ represent the infinite directions.

\setlength{\unitlength}{0.0003in}
\begingroup\makeatletter\ifx\SetFigFont\undefined%
\gdef\SetFigFont#1#2#3#4#5{%
  \reset@font\fontsize{#1}{#2pt}%
  \fontfamily{#3}\fontseries{#4}\fontshape{#5}%
  \selectfont}%
\fi\endgroup%
{\renewcommand{\dashlinestretch}{30}
\begin{picture}(5347,5182)(-4000,1000)
\thinlines
\path(975.000,982.000)(825.000,922.000)(975.000,862.000)
\path(825,922)(5325,922)
\path(3525,3622)(3525,22)
\path(3465.000,202.000)(3525.000,22.000)(3585.000,202.000)
\path(2625,2722)(2625,22)
\path(2565.000,202.000)(2625.000,22.000)(2685.000,202.000)
\path(1725,1822)(1725,22)
\path(1665.000,202.000)(1725.000,22.000)(1785.000,202.000)
\path(4425,4522)(4425,22)
\path(4365.000,202.000)(4425.000,22.000)(4485.000,202.000)
\path(1905.000,1882.000)(1725.000,1822.000)(1905.000,1762.000)
\path(1725,1822)(5325,1822)
\path(2805.000,2782.000)(2625.000,2722.000)(2805.000,2662.000)
\path(2625,2722)(5325,2722)
\path(3705.000,3682.000)(3525.000,3622.000)(3705.000,3562.000)
\path(3525,3622)(5325,3622)
\put(4350,5047){\makebox(0,0)[lb]{\smash{{{\SetFigFont{14}{16.8}{\rmdefault}{\mddefault}{\updefault}$\vdots$}}}}}
\put(0,847){\makebox(0,0)[lb]{\smash{{{\SetFigFont{14}{16.8}{\rmdefault}{\mddefault}{\updefault}$\hdots$}}}}}
\end{picture}
}

\vspace*{4mm}
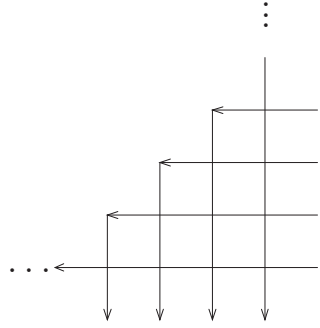
\begin{center}{\ft Figure 3.2: The North-West vertex model corner transfer matrix $\cA(u)$}\end{center}

\nin This and other graphical representations should 
be read is as follows: the matrix element
$\cA(u)_{\cdots,\ep'_2,\ep'_1}^{\cdots,\ep_2,\ep_1}$
is obtained by computing the weighted sum associated with the lattice 
in Figure 3.2, with all internal edge variables summed over, and with the
West external horizontal edge variables and  North vertical
edge variables fixed to the values
\ben 
\cdots \,\ep'_2 \,\ep'_1\quad \quad  \hb{        and
} \quad \quad {{\vdots\atop{\displaystyle{\ep_2}}}\atop{\ep_1}}
\een 
respectively. Clearly, one can then define South-West, South-East and
North-East corner transfer matrices in an analogous manner. A
simple consideration of the boundary conditions 
establishes that these operators act as 
\ben \cA_{SW}(u): \cH^{(\ell)}\ra \cH^{(\bar{\ell})},\quad \cA_{SE}(u):
\cH^{(\bar{\ell})}\ra \cH^{(\bar{\ell})},
\quad \cA_{NE}(u):\cH^{(\bar{\ell})}\ra \cH^{(\ell)}\quad \hb{where}\ws\bar{\ell}=
k-\ell.\een
It is an easy exercise to show that the crossing symmetry relation
\mref{qconj} implies that these new operators are related to $\cA(u)$ by
\bea \cA_{SW}(u) =\cQ \,\cA(-1-u),\quad \cA_{SE}(u) =\cQ\,\cA(u)\,\cQ^{-1} ,\quad \cA_{NE}(u)
=\cA(-1-u)\,\cQ^{-1},\lb{vctmreln}\ee
where
 $\cQ:\oplus_{\ell} \cH^{(\ell)} \ra \oplus_{\ell} \cH^{(\ell)}$ is the operator \ben \cQ=\cdots \ot Q \ot Q \ot Q.\een

Baxter's key observation about the corner transfer matrix is
that in the infinite volume limit we have $\cA(u)\sim x^{-2u H^{(\ell)}}$,
where $2 H^{(\ell)}$, the  corner Hamiltonian, has discrete and equidistant 
eigenvalues bounded from below,  
%only has  eigenvalues in $\N$, 
and $\sim$ means equal up to a scalar.
Thus, from \eqref{vctmreln}, we have
\bea
&&\cA_{NE}(u)\cA_{SE}(u)\cA_{SW}(u) \cA(u)\sim 
%\rho^{(\ell)}:=
\ x^{4H^{(\ell)}}.\lb{vctm2h}
\ena  

In terms of CTMs, the partition function $Z^{(\ell)}$ is expressed by 
\bea
Z^{(\ell)}&=& 
\Tr_{\cH^{(\ell)}}\Big( \cA_{NE}(u)\cA_{SE}(u)\cA_{SW}(u) \cA(u)\Big)
\sim \Tr_{\cH^{(\ell)}}\ x^{4H^{(\ell)}}. \lb{vpart}
\ena
It is remarkable that this partition function is known to coincides with 
 the principally 
specialised character of the level 
$k$ irreducible highest-weight $\slth$-module $V(\la_\ell)$ with 
highest weight $\la_\ell=(k-\ell)\Lambda_0+\ell\Lambda_1\ (\ell=0,1,..,k)$
\cite{DateLMP89}. 
Here $\Lambda_i\ (i=0,1)$ denotes the fundamental weight of $\slth$.
Namely we have
\bea
Z^{(\ell)}&\sim&\Tr_{\cH^{(\ell)}}\ x^{4H^(\ell)}=\chi^{(k)}_\ell(\bar{\tau}),\lb{part2ch}\\
\chi^{(k)}_\ell(\bar{\tau})&=&
\frac{x^{\frac{1}{2}}}{(x^2;x^2)_\infty(x^2;x^4)_\infty}[\ell+1]^{(k+2)}.\lb{chprincipal}
\ena
Here we set $e^{2\pi i \bar{\tau}}=x^4$.
Later we will use another expression of $\chi^{(k)}_\ell(\bar{\tau})$ in terms of
 the string function \cite{KaPe}:
\bea
&&\chi^{(k)}_{\ell}(\bar{\tau})=\sum_{n\in\Z}\sum_{M\equiv 
0 \atop {\rm mod} 2k}^{2k-1}\ 
c^{\la_\ell}_{\la_{M}}(\bar{\tau})\ x^{4k(n+\frac{M}{2k})^2-2k(n+\frac{M}{2k})},
\lb{ch1}\ena
where $c^{\la_\ell}_{\la_M}(\bar{\tau})$ 
denotes the string function defined  by 
\bea
c^{\la_\ell}_{\la_M}(\bar{\tau})&=&(x^4)^{\frac{\ell(\ell+2)}{4(k+2)}-\frac{M^2}{4k}-\frac{k}{8(k+2)}}\sum_{n\geq 0 }{\rm dim} V({\la_\ell})_{\la_M-n\delta}\ (x^4)^n\lb{stfn}
\ena
and $c^{\la_\ell}_{\la_M}(\bar{\tau})=0$ for $M\not\equiv \ell$ mod 2.
Here $V(\la_\ell)_\la$ denotes the  weight space 
of  $V(\la_\ell)$; $\delta$ denotes the null root of $\slth$ satisfying  
$(\delta,\delta)=0=(\delta,\al),\ (\delta,d)=1$ with a standard symmetric bilinear 
form $(\ ,\ ):P\times P \to \frac{1}{2}\Z$, $P=\Z\Lambda_0\oplus\Z\Lambda_1\oplus \Z\delta$. 
Note that the string function satisfies the following relations.
\be
&&c^{\la_\ell}_{\la_M}(\bar{\tau})=c^{\la_\ell}_{\la_{-M}}(\bar{\tau})=c^{\la_{k-\ell}}_{\la_{k-M}}(\bar{\tau})
=c^{\la_\ell}_{\la_{M+2k}}(\bar{\tau}).
\en

\subsubsection{HTMs}
To express the correlation functions in a similar way to \eqref{vpart}, we 
need to introduce HTMs, that is, 
North and South half-transfer matrices, denoted by
$\phi^{(\bar{\ell},\ell)}_\ep(u):\cH^{(\ell)}\ra \cH^{(\bar{\ell})}$
and $\phi^{(\bar{\ell},\ell)}_{S;\,\ep}(u):\cH^{(\ell)}\ra \cH^{(\bar{\ell})}$ 
respectively, and
defined graphically by Figure 3.3.

\setlength{\unitlength}{0.0003in}
\begingroup\makeatletter\ifx\SetFigFont\undefined%
\gdef\SetFigFont#1#2#3#4#5{%
  \reset@font\fontsize{#1}{#2pt}%
  \fontfamily{#3}\fontseries{#4}\fontshape{#5}%
  \selectfont}%
\fi\endgroup%
{\renewcommand{\dashlinestretch}{30}
\begin{picture}(9172,5475)(-2000,800)
\thinlines
\path(1950,4950)(1950,450)
\path(1890.000,600.000)(1950.000,450.000)(2010.000,600.000)
\path(1200.000,4110.000)(1050.000,4050.000)(1200.000,3990.000)
\path(1050,4050)(2850,4050)
\path(1200.000,3210.000)(1050.000,3150.000)(1200.000,3090.000)
\path(1050,3150)(2850,3150)
\path(1200.000,2310.000)(1050.000,2250.000)(1200.000,2190.000)
\path(1050,2250)(2850,2250)
\path(1200.000,1410.000)(1050.000,1350.000)(1200.000,1290.000)
\path(1050,1350)(2850,1350)
\path(8250,4950)(8250,450)
\path(8190.000,600.000)(8250.000,450.000)(8310.000,600.000)
\path(7500.000,4110.000)(7350.000,4050.000)(7500.000,3990.000)
\path(7350,4050)(9150,4050)
\path(7500.000,3210.000)(7350.000,3150.000)(7500.000,3090.000)
\path(7350,3150)(9150,3150)
\path(7500.000,2310.000)(7350.000,2250.000)(7500.000,2190.000)
\path(7350,2250)(9150,2250)
\path(7500.000,1410.000)(7350.000,1350.000)(7500.000,1290.000)
\path(7350,1350)(9150,1350)
\put(1875,5325){\makebox(0,0)[lb]{{\SetFigFont{12}{16.8}{\rmdefault}{\mddefault}{\updefault}$\vdots$}}}
\put(0,3000){\makebox(0,0)[lb]{{\SetFigFont{12}{16.8}{\rmdefault}{\mddefault}{\updefault}$a)$}}}
\put(2175,675){\makebox(0,0)[lb]{{\SetFigFont{12}{16.8}{\rmdefault}{\mddefault}{\updefault}$\ep$}}}
\put(8550,4875){\makebox(0,0)[lb]{{\SetFigFont{12}{16.8}{\rmdefault}{\mddefault}{\updefault}$\ep$}}}
\put(8175,0){\makebox(0,0)[lb]{{\SetFigFont{12}{16.8}{\rmdefault}{\mddefault}{\updefault}$\vdots$}}}
\put(5925,3000){\makebox(0,0)[lb]{{\SetFigFont{12}{16.8}{\rmdefault}{\mddefault}{\updefault}$b)$}}}
\end{picture}
}

\vspace*{4mm}
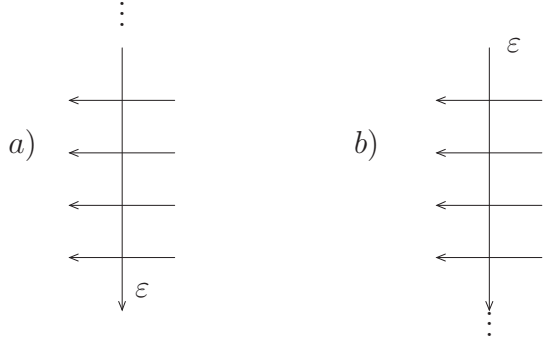
\begin{center}{\ft Figure 3.3: a) The North half-transfer matrix $\phi_\ep(u)$  \\[-1mm]
$ $\hspace*{23mm} b) The South half-transfer matrix $\phi_{S;\,\ep}(u)$}\end{center}

\nin These operators are viewed as acting
in
an anti-clockwise direction about the finite end, i.e., the end whose edge variable is 
fixed to the value $\ep$. Again, crossing symmetry implies that
we have the relation
\bea 
\phi_{S;\,\ep}^{(\bell,\ell)}(u)=
\cQ\ \phi^{*(\bell,\ell)}_{\vep}(u)\ \cQ^{-1}
\lb{vphireln}\ee
where we define the dual operator $\phi^{*(\bell,\ell)} _{\vep}(u)$ by
\be
&&\phi^{*(\bell,\ell)}_{\vep}(u)=\sum_{\vep'}Q^{\vep'}_{\vep}
\phi^{(\bell,\ell)}_{\vep'}(u-1).
\en 
We will often suppress the $(\bell,\ell)$
superscripts on these various operators.

The heuristic graphical arguments of \cite{Fodal93} then lead to the
following relations for half transfer and 
corner transfer matrices:
\bea
\phi^{(\ell,\bell)}_{\vep_2}(u_2)\phi^{(\bell,\ell)}_{\vep_1}(u_1)
&=&\hspace*{-5mm}\sum_{\vep_1',\,\vep_2'\in\{-k,-k+2,..,k\}} 
\hspace*{-10mm} R^{(k,k)}(u_1-u_2)^{\vep'_1\vep_2'}_{\vep_1\vep_2}
\phi^{(\ell,\bell)}_{\vep_1'}(u_1)\phi^{(\bell,\ell)}_{\vep_2'}(u_2),\lb{Rphiphi}\\
\sum_{\vep\in\{-k,-k+2,..,k\}}\hspace*{-10mm}\phi^{*(\ell,\bell)}_{\vep}(u)
\phi^{(\bell,\ell)}_{\vep}(u)&=&\id,\qquad\qquad\qquad\qquad\qquad\qquad\qquad\qquad\qquad \lb{8vinversion}\\
\cA(u) \phi_\ep(v)&=&\phi_\ep(v-u) \cA(u).\lb{vctmphireln}
\ena
Furthermore \eqref{vctm2h} and \eqref{vctmphireln} yield 
\bea
\phi^{(\bell,\ell)}_{\vep}(u)x^{4H^{(\ell)}}&=&x^{4H^{(\bell)}}\phi^{(\bell,\ell)}_{\vep}(u-2).\lb{rhophi}
\ena

\subsubsection{Correlation functions}
Now let us divide the lattice depicted in Figure 3.1 into the pieces corresponding to 
CTMs and HTMs. We obtain Figure 3.4. 

\setlength{\unitlength}{0.00035in}
\begingroup\makeatletter\ifx\SetFigFont\undefined%
\gdef\SetFigFont#1#2#3#4#5{%
  \reset@font\fontsize{#1}{#2pt}%
  \fontfamily{#3}\fontseries{#4}\fontshape{#5}%
  \selectfont}%
\fi\endgroup%
{\renewcommand{\dashlinestretch}{30}
\begin{picture}(11507,8715)(-1500,500)
\put(5550,6825){\makebox(0,0)[lb]{\smash{{{\SetFigFont{8}{9.6}{\rmdefault}{\mddefault}{\updefault}$\hdots$}}}}}
\put(11400,3525){\makebox(0,0)[lb]{\smash{{{\SetFigFont{8}{9.6}{\rmdefault}{\mddefault}{\updefault}$\hdots$}}}}}
\put(11400,5325){\makebox(0,0)[lb]{\smash{{{\SetFigFont{8}{9.6}{\rmdefault}{\mddefault}{\updefault}$\hdots$}}}}}
\put(-200,5325){\makebox(0,0)[lb]{\smash{{{\SetFigFont{8}{9.6}{\rmdefault}{\mddefault}{\updefault}$\hdots$}}}}}
\put(-200,3525){\makebox(0,0)[lb]{\smash{{{\SetFigFont{8}{9.6}{\rmdefault}{\mddefault}{\updefault}$\hdots$}}}}}
\put(7125,3825){\makebox(0,0)[lb]{\smash{{{\SetFigFont{8}{9.6}{\rmdefault}{\mddefault}{\updefault}$\ep_N$}}}}}
\put(7125,4950){\makebox(0,0)[lb]{\smash{{{\SetFigFont{8}{9.6}{\rmdefault}{\mddefault}{\updefault}$\ep_N$}}}}}
\put(4875,3825){\makebox(0,0)[lb]{\smash{{{\SetFigFont{8}{9.6}{\rmdefault}{\mddefault}{\updefault}$\ep_1$}}}}}
\put(4875,4950){\makebox(0,0)[lb]{\smash{{{\SetFigFont{8}{9.6}{\rmdefault}{\mddefault}{\updefault}$\ep_1$}}}}}
%\thicklines
\path(8250.000,2445.000)(8100.000,2400.000)(8250.000,2355.000)
\path(8100,2400)(9900,2400)
\path(8250.000,3045.000)(8100.000,3000.000)(8250.000,2955.000)
\path(8100,3000)(10500,3000)
\path(9300,1800)(8100,1800)
\path(8250.000,1845.000)(8100.000,1800.000)(8250.000,1755.000)
\path(1050,4200)(1050,3000)
\path(1005.000,3150.000)(1050.000,3000.000)(1095.000,3150.000)
\path(1650,4200)(1650,2400)
\path(1605.000,2550.000)(1650.000,2400.000)(1695.000,2550.000)
\path(2250,4200)(2250,1800)
\path(2205.000,1950.000)(2250.000,1800.000)(2295.000,1950.000)
\path(8250.000,3645.000)(8100.000,3600.000)(8250.000,3555.000)
\path(8100,3600)(11100,3600)
\path(10500,4200)(10500,3000)
\path(10455.000,3150.000)(10500.000,3000.000)(10545.000,3150.000)
\path(9900,4200)(9900,2400)
\path(9855.000,2550.000)(9900.000,2400.000)(9945.000,2550.000)
\path(9300,4200)(9300,1800)
\path(9255.000,1950.000)(9300.000,1800.000)(9345.000,1950.000)
\put(5550,6150){\makebox(0,0)[lb]{\smash{{{\SetFigFont{8}{9.6}{\rmdefault}{\mddefault}{\updefault}$\hdots$}}}}}
\put(6750,8850){\makebox(0,0)[lb]{\smash{{{\SetFigFont{8}{9.6}{\rmdefault}{\mddefault}{\updefault}$\phi_{\ep_N}(u)$}}}}}
\put(4375,8850){\makebox(0,0)[lb]{\smash{{{\SetFigFont{8}{9.6}{\rmdefault}{\mddefault}{\updefault}$\phi_{\ep_1}(u)$}}}}}
\put(6750,100){\makebox(0,0)[lb]{\smash{{{\SetFigFont{8}{9.6}{\rmdefault}{\mddefault}{\updefault}$\phi_{S\ep_N}(u)$}}}}}
\put(4375,100){\makebox(0,0)[lb]{\smash{{{\SetFigFont{8}{9.6}{\rmdefault}{\mddefault}{\updefault}$\phi_{S\ep_1}(u)$}}}}}
\put(5650,100){\makebox(0,0)[lb]{\smash{{{\SetFigFont{8}{9.6}{\rmdefault}{\mddefault}{\updefault}$\hdots$}}}}}
\put(5550,8950){\makebox(0,0)[lb]{\smash{{{\SetFigFont{8}{9.6}{\rmdefault}{\mddefault}{\updefault}$\hdots$}}}}}
\put(10050,1875){\makebox(0,0)[lb]{\smash{{{\SetFigFont{8}{9.6}{\rmdefault}{\mddefault}{\updefault}$\cA_{SW}(u)$}}}}}
\put(10125,6900){\makebox(0,0)[lb]{\smash{{{\SetFigFont{8}{9.6}{\rmdefault}{\mddefault}{\updefault}$\cA_{NE}(u)$}}}}}
\put(400,1875){\makebox(0,0)[lb]{\smash{{{\SetFigFont{8}{9.6}{\rmdefault}{\mddefault}{\updefault}$\cA_{SW}(u)$}}}}}
\put(400,6900){\makebox(0,0)[lb]{\smash{{{\SetFigFont{8}{9.6}{\rmdefault}{\mddefault}{\updefault}$\cA(u)$}}}}}
\put(8625,8100){\makebox(0,0)[lb]{\smash{{{\SetFigFont{8}{9.6}{\rmdefault}{\mddefault}{\updefault}$\vdots$}}}}}
\put(8625,625){\makebox(0,0)[lb]{\smash{{{\SetFigFont{8}{9.6}{\rmdefault}{\mddefault}{\updefault}$\vdots$}}}}}
\put(6825,625){\makebox(0,0)[lb]{\smash{{{\SetFigFont{8}{9.6}{\rmdefault}{\mddefault}{\updefault}$\vdots$}}}}}
\put(4575,625){\makebox(0,0)[lb]{\smash{{{\SetFigFont{8}{9.6}{\rmdefault}{\mddefault}{\updefault}$\vdots$}}}}}
\put(6825,8100){\makebox(0,0)[lb]{\smash{{{\SetFigFont{8}{9.6}{\rmdefault}{\mddefault}{\updefault}$\vdots$}}}}}
\put(4575,8100){\makebox(0,0)[lb]{\smash{{{\SetFigFont{8}{9.6}{\rmdefault}{\mddefault}{\updefault}$\vdots$}}}}}
\put(2775,625){\makebox(0,0)[lb]{\smash{{{\SetFigFont{8}{9.6}{\rmdefault}{\mddefault}{\updefault}$\vdots$}}}}}
\put(2775,8100){\makebox(0,0)[lb]{\smash{{{\SetFigFont{8}{9.6}{\rmdefault}{\mddefault}{\updefault}$\vdots$}}}}}
\put(5550,2625){\makebox(0,0)[lb]{\smash{{{\SetFigFont{8}{9.6}{\rmdefault}{\mddefault}{\updefault}$\hdots$}}}}}
\put(5550,1950){\makebox(0,0)[lb]{\smash{{{\SetFigFont{8}{9.6}{\rmdefault}{\mddefault}{\updefault}$\hdots$}}}}}
\put(5550,3225){\makebox(0,0)[lb]{\smash{{{\SetFigFont{8}{9.6}{\rmdefault}{\mddefault}{\updefault}$\hdots$}}}}}
\put(5550,5550){\makebox(0,0)[lb]{\smash{{{\SetFigFont{8}{9.6}{\rmdefault}{\mddefault}{\updefault}$\hdots$}}}}}
\path(8700,4200)(8700,1200)
\path(8655.000,1350.000)(8700.000,1200.000)(8745.000,1350.000)
\path(6900,4200)(6900,1200)
\path(6855.000,1350.000)(6900.000,1200.000)(6945.000,1350.000)
\path(6450.000,5445.000)(6300.000,5400.000)(6450.000,5355.000)
\path(6300,5400)(7500,5400)
\path(6450.000,6045.000)(6300.000,6000.000)(6450.000,5955.000)
\path(6300,6000)(7500,6000)
\path(6450.000,6645.000)(6300.000,6600.000)(6450.000,6555.000)
\path(6300,6600)(7500,6600)
\path(6450.000,7245.000)(6300.000,7200.000)(6450.000,7155.000)
\path(6300,7200)(7500,7200)
\path(6900,7800)(6900,4800)
\path(6855.000,4950.000)(6900.000,4800.000)(6945.000,4950.000)
\path(4650,7800)(4650,4800)
\path(4605.000,4950.000)(4650.000,4800.000)(4695.000,4950.000)
\path(4650,4200)(4650,1200)
\path(4605.000,1350.000)(4650.000,1200.000)(4695.000,1350.000)
\path(4200.000,5445.000)(4050.000,5400.000)(4200.000,5355.000)
\path(4050,5400)(5250,5400)
\path(4200.000,6045.000)(4050.000,6000.000)(4200.000,5955.000)
\path(4050,6000)(5250,6000)
\path(4200.000,6645.000)(4050.000,6600.000)(4200.000,6555.000)
\path(4050,6600)(5250,6600)
\path(4200.000,7245.000)(4050.000,7200.000)(4200.000,7155.000)
\path(4050,7200)(5250,7200)
\path(2850,4200)(2850,1200)
\path(2805.000,1350.000)(2850.000,1200.000)(2895.000,1350.000)
\path(3450,1800)(2250,1800)
\path(2400.000,1845.000)(2250.000,1800.000)(2400.000,1755.000)
\path(3450,2400)(1650,2400)
\path(1800.000,2445.000)(1650.000,2400.000)(1800.000,2355.000)
\path(3450,3000)(1050,3000)
\path(1200.000,3045.000)(1050.000,3000.000)(1200.000,2955.000)
\path(3450,3600)(450,3600)
\path(600.000,3645.000)(450.000,3600.000)(600.000,3555.000)
\path(1050,6000)(1050,4800)
\path(1005.000,4950.000)(1050.000,4800.000)(1095.000,4950.000)
\path(1650,6600)(1650,4800)
\path(1605.000,4950.000)(1650.000,4800.000)(1695.000,4950.000)
\path(2250,7200)(2250,4800)
\path(2205.000,4950.000)(2250.000,4800.000)(2295.000,4950.000)
\path(6450.000,3645.000)(6300.000,3600.000)(6450.000,3555.000)
\path(6300,3600)(7500,3600)
\path(8250.000,5445.000)(8100.000,5400.000)(8250.000,5355.000)
\path(8100,5400)(11100,5400)
\path(10500,6000)(10500,4725)
\path(10455.000,4875.000)(10500.000,4725.000)(10545.000,4875.000)
\path(8250.000,6045.000)(8100.000,6000.000)(8250.000,5955.000)
\path(8100,6000)(10500,6000)
\path(9900,6600)(9900,4800)
\path(9855.000,4950.000)(9900.000,4800.000)(9945.000,4950.000)
\path(8250.000,6645.000)(8100.000,6600.000)(8250.000,6555.000)
\path(8100,6600)(9900,6600)
\path(9300,7200)(9300,4800)
\path(9255.000,4950.000)(9300.000,4800.000)(9345.000,4950.000)
\path(8700,7800)(8700,4800)
\path(8655.000,4950.000)(8700.000,4800.000)(8745.000,4950.000)
\path(8250.000,7245.000)(8100.000,7200.000)(8250.000,7155.000)
\path(8100,7200)(9300,7200)
\path(600.000,5445.000)(450.000,5400.000)(600.000,5355.000)
\path(450,5400)(3450,5400)
\path(1200.000,6045.000)(1050.000,6000.000)(1200.000,5955.000)
\path(1050,6000)(3450,6000)
\path(1800.000,6645.000)(1650.000,6600.000)(1800.000,6555.000)
\path(1650,6600)(3450,6600)
\path(2400.000,7245.000)(2250.000,7200.000)(2400.000,7155.000)
\path(2250,7200)(3450,7200)
\path(2805.000,4950.000)(2850.000,4800.000)(2895.000,4950.000)
\path(2850,4800)(2850,7800)
\path(4200.000,1845.000)(4050.000,1800.000)(4200.000,1755.000)
\path(4050,1800)(5250,1800)
\path(4200.000,2445.000)(4050.000,2400.000)(4200.000,2355.000)
\path(4050,2400)(5250,2400)
\path(4200.000,3045.000)(4050.000,3000.000)(4200.000,2955.000)
\path(4050,3000)(5250,3000)
\path(4200.000,3645.000)(4050.000,3600.000)(4200.000,3555.000)
\path(4050,3600)(5250,3600)
\path(6450.000,1845.000)(6300.000,1800.000)(6450.000,1755.000)
\path(6300,1800)(7500,1800)
\path(6450.000,2445.000)(6300.000,2400.000)(6450.000,2355.000)
\path(6300,2400)(7500,2400)
\path(6450.000,3045.000)(6300.000,3000.000)(6450.000,2955.000)
\path(6300,3000)(7500,3000)
\end{picture}
}

\begin{center}
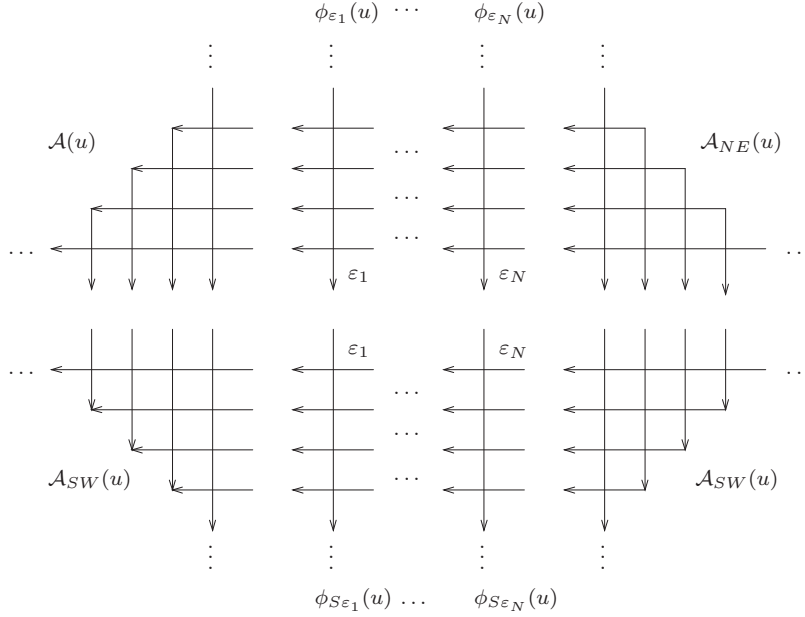
{\ft Figure 3.4: The vertex model correlation function trace}
\end{center}

According to this picture, we can express   
the correlation function $P^{(\ell)}(\ep_1,\ep_2,\cdots,\ep_N)$ as follows.
\bea P^{(\ell)}(\ep_1,\ep_2,\cdots,\ep_N)=
\frac{1}{Z^{(\ell)}}
F^{(\ell)}(\ep_1,\ep_2,\cdots,\ep_N)\lb{vcorr1}\ee 
where 
\bea 
&&
\hspace*{-18mm}F^{(\ell)}(\ep_1,\ep_2,\cdots,\ep_N)\nn \\
&&\hspace*{-15mm}=
\Tr_{\cH^{\sigma^N\!(\ell)}} \Big( 
\cA_{NE}(u)\cA_{SE}(u) \phi_{S\,\ep_N}(u)\cdots \phi_{S\,\ep_1}(u) \cA_{SW}(u) \cA(u)
\phi_{\ep_1}(u)\cdots \phi_{\ep_N}(u)
\Big)
\ee
with $\sigma(\ell)=k-\ell$.
One can then use the relations \mref{vctmreln}, \mref{vphireln} and \mref{vctmphireln} to write all operators
in terms of $\phi_{\ep}(u)$ and $\cA(u)$ and to re-order them.
We thus obtain the following simplified expression for the correlation 
function:
\bea  
&&P^{(\ell)}(\ep_1,\ep_2,\cdots,\ep_N)\nn\\
&&=\ws\frac{1}{\chi^{(k)}_{\ell}(\bar{\tau})}
 {\Tr_{\cH^{\sigma^N\!(\ell)}} \Big( x^{4H^{(\ell)}}\,
 \phi^{*(\ell,\bell)}_{\vep_N}(u)\cdots \phi^{*(\sigma^{N+1}(\ell),\sigma^{N}(\ell))}_{\vep_1}(u) \,
\phi^{(\sigma^{N}(\ell),\sigma^{N-1}(\ell))}_{\vep_1}(u)\cdots 
\phi^{(\bell,\ell)}_{\vep_N}(u)
\Big)
}.\nn\\ &&\lb{vtrace}\ee

\ssect{Fusion SOS Models}\lb{fusionSOS}
We next recall the analogous technology
developed in \cite{Fodal93} to write the infinite-volume limit of 
SOS correlation functions as 
traces of CTMs and HTMs. 

\subsubsection{The space of states}
The first step is to define the space of states $\cH^{(\ell)}_{m,a}$
on which our various SOS operators act. We define $\cH^{(\ell)}_{m,a}$ ($\ell\in \{0,1,\cdots,k\}$, $m\in \Z$, $a\in 2\Z+m+\ell$) by
\ben 
\cH^{(\ell)}_{m,a}=&&\hb{Span}_{\C}\left\{ \cdots \ot v_{s(2)}\ot v_{s(1)} \ot
  v_{s(0)}|\right. \\ && \left. s(i)\in \Z,s(i+1)-s(i)\in \{-k,-k+2,\cdots,k\},
s(i)=s^{(\ell)}_{m}(i)\hb{ for } i\gg 0, s(0)=a  
\right\}\\
s^{(\ell)}_m(i)&=& 
\begin{cases} 
  m+\ell & i=0 \mod 2\\
m+\bar{\ell} & i=1 \mod 2\end{cases}
\een
 $\cH^{(\ell)}_{m,a}$ is the vector space associated with the height
 variables along a line running from the centre of a lattice to the boundary, for
which the central height is fixed to $a$, and far from the centre the
boundary heights are fixed to the ground state configuration 
$\cdots ,m+\bell,m+\ell,m+\bell,m+\ell,m+\bell ,\cdots$.

\subsubsection{CTMs and the partition function}
The infinite volume North-West corner transfer matrix $A_{a}(u)$ is
now defined graphically by

\setlength{\unitlength}{0.0003in}
\begingroup\makeatletter\ifx\SetFigFont\undefined%
\gdef\SetFigFont#1#2#3#4#5{%
  \reset@font\fontsize{#1}{#2pt}%
  \fontfamily{#3}\fontseries{#4}\fontshape{#5}%
  \selectfont}%
\fi\endgroup%
{\renewcommand{\dashlinestretch}{30}
\begin{picture}(4848,4200)(-6000,600)
\put(4725,15){\makebox(0,0)[lb]{\smash{{{\SetFigFont{12}{16.8}{\rmdefault}{\mddefault}{\updefault}$a$}}}}}
\put(0,540){\makebox(0,0)[lb]{\smash{{{\SetFigFont{12}{16.8}{\rmdefault}{\mddefault}{\updefault}$\hdots$}}}}}
\put(3825,4065){\makebox(0,0)[lb]{\smash{{{\SetFigFont{12}{16.8}{\rmdefault}{\mddefault}{\updefault}$\vdots$}}}}}
\thinlines
\path(1650,1890)(4350,1890)
\path(2550,2790)(4350,2790)
\path(750,990)(750,90)
\path(1650,1890)(1650,90)
\path(2550,2790)(2550,90)
\path(750,90)(4350,90)
\path(4350,3690)(4350,90)
\path(3450,3690)(4350,3690)
\path(3675,990)(3450,765)
\path(2775,990)(2550,765)
\path(1875,990)(1650,765)
\path(975,990)(750,765)
\path(3675,1890)(3450,1665)
\path(2775,1890)(2550,1665)
\path(1875,1890)(1650,1665)
\path(3675,2790)(3450,2565)
\path(2775,2790)(2550,2565)
\path(3675,3690)(3450,3465)
\path(750,990)(4350,990)
\path(3450,3690)(3450,90)
\texture{55888888 88555555 5522a222 a2555555 55888888 88555555 552a2a2a 2a555555 
        55888888 88555555 55a222a2 22555555 55888888 88555555 552a2a2a 2a555555 
        55888888 88555555 5522a222 a2555555 55888888 88555555 552a2a2a 2a555555 
        55888888 88555555 55a222a2 22555555 55888888 88555555 552a2a2a 2a555555 }
\put(4350,90){\shade\ellipse{150}{150}}
\put(4350,90){\ellipse{150}{150}}
\end{picture}
}

\vspace*{4mm}\begin{center}{\ft Figure 3.5: The North-West corner transfer matrix $A_{a}(u)$}
\end{center}

\nin where we suppress the appearance of the spectral parameter $u$ associated
with each SOS face weight.
Note that the centre height is fixed to $a$ and all internal height
variables are summed over. $A_{a}(u)$ preserves the boundary
conditions, and can be viewed as an operator that acts in an
anti-clockwise direction about the centre of the lattice as  
$A_a(u):\cH^{(\ell)}_{m,a}\ra \cH^{(\ell)}_{m,a}$. One can
define the South-West, South-East and North-East corner transfer matrices 
with fixed central heights in an analogous manner. 
Using the crossing symmetry relations \mref{crossing}, allows us to write each
of these in terms $A_{a}(u)$. As for vertex models, it is a simple, but 
illuminating exercise, to show that we have
\bea
A_{SW;a}(u)= g_a \Gamma A_{a}(-1-u),\ws A_{SE;a}(u)= \Gamma\, A_{a}(u)
\Gamma^{-1},\ws A_{NE;a}(u)=g_a A_{a}(-1-u)\, \Gamma^{-1}\nn\\[-2mm]
\lb{idctm2}\ee
where $\Gamma:\cH^{(\ell)}_{m,a}\ra \cH^{(\ell)}_{m,a}$ is defined by 
\ben
\\[-5mm]
\cdots \ot v_{s(2)}\ot v_{s(1)} \ot
  v_{s(0)} &\goto&  \cdots \,(s(2),s(1))_k\, (s(1),s(0))_k \,
     \left(\cdots \ot v_{s(2)}\ot v_{s(1)} \ot
  v_{s(0)}\right) 
  \een
and $g_a$ and $(a,b)_k$ are as previously defined. 

Again, in parallel to the vertex case, it is known that in the infinite volume limit 
we have
$A_a(u)\sim x^{-2u H_{m,a}^{(\ell)}}$, where the corner Hamiltonian  
$H^{(\ell)}_{m,a}$ has discrete and equidistant 
eigenvalues bounded from below.  
%a spectrum in $\N$. 
From \eqref{idctm2}, we therefore have
\bea
A_{NE;a}(u) A_{SE;a}(u) A_{SW;a}(u)A_a(u)\sim
\ [a]\ 
x^{4H^{(\ell)}_{m,a}}.\lb{fctm2h}
\ena

Then the partition function $Z^{(\ell)}_m$ is expressed by CTMs as follows.
\be
Z^{(\ell)}_m&=&\sum_{a\in m+\ell+2\Z} \Tr_{\cH^{(\ell)}_{m,a}}
\big(A_{NE;a_N}(u)A_{SE;a'_N}(u)A_{SW;a}(u)A_{a}(u)\big),\nn\\
&\sim&\sum_{a\in m+\ell+2\Z}[a] \, \Tr_{\H^{(\ell)}_{m,a}} x^{4H^{(\ell)}_{m,a}}.
\en
It is known that the trace part is given by the string function
\cite{DateLMP89,DJKMO89} 
\bea
\Tr_{\H^{(\ell)}_{m,a}} x^{4H^{(\ell)}_{m,a}}&=&
c^{\la_\ell}_{\la_{M}}(\bar{\tau})\  
x^{\frac{(mr-ar^*)^2}{krr^*}},\lb{foneptf}
\ena
where $M\equiv a-m\ \mod 2$. 
Then by the calculation given in Appendix \ref{ChPartition}, 
we obtain the following expression of the partition function.  
\begin{thm}\lb{charpartSOS}
\bea
Z^{(\ell)}_m&\sim&\sum_{a\in m+\ell+2\Z}[a] \, \Tr_{\H^{(\ell)}_{m,a}} x^{4H^{(\ell)}_{m,a}}=[m]^*
\chi^{(k)}_\ell(\bar{\tau}),\lb{charpart}
\ena
were $\chi^{(k)}_\ell(\bar{\tau})$ is the principally specialised character 
given by \eqref{ch1}.
\end{thm}
\noi The case $k=1$ was obtained by Lashkevich and Pugai\cite{LaP98}. 
Note that this represents the vertex-face correspondence 
between the spaces of states. 

It is worth noting the resemblance of \eqref{charpart} to the branching formula 
for the product of characters of irreducible integrable representations of $\slth$:
\bea
&&\chi^{(k)}_\ell(\bar{\tau})\chi^{(r-k-2)}_{m-1}(\bar{\tau})=
\sum_{1\leq a\leq r-1}b^{(\ell)}_{m,a}(\bar{\tau}) \chi^{(r-2)}_{a-1}(\bar{\tau}),  
\lb{branch}
\ena
where the principally specialised character 
$\chi^{(s)}_{a}(\bar{\tau})$ is given in \eqref{chprincipal}. The branching function 
$b^{(\ell)}_{m,a}(\bar{\tau})$ is known to be the 
character of the irreducible Virasoro module $Vir_{m,a}$ associated with the 
coset $(\slth)_k\oplus (\slth)_{r-k-2}/(\slth)_{r-2}$, and with the 
highest weight $h_{m,a}=\frac{\ell(k-\ell)}{2k(k+2)}
+\frac{(mr-ar^*)^2-k^2}{4krr^*}$ and central charge 
$c_{Vir}=\frac{3k}{k+2}\left(1-\frac{2(k+2)}{rr^*}\right)$. The main difference 
between the two formulae are:

\begin{itemize}
\item[1)]  \eqref{branch} corresponds to 
the direct sum decomposition of the tensor product representation  
%$V(\la^{(k)}_\ell)\otimes V(\la^{(r-k-2)}_{m-1})$. 
\bea
&&V(\la^{(k)}_\ell)\otimes V(\la^{(r-k-2)}_{m-1})
=\bigoplus_{1\leq a\leq r-1}\Omega^{(\ell)}_{m,a}\otimes V(\la^{(r-2)}_{a-1}).
\lb{VVdecomp}
\ena
Here $V(\la^{(s)}_a)$ denotes the level $s$ irreducible integrable representation 
with the highest weight $\la^{(s)}_a=(s-a)\Lambda_0+a\Lambda_1$ being 
dominant integral. $\Omega^{(\ell)}_{m,a}$ denotes 
the corresponding irreducible coset Virasoro module $Vir_{m,a}$. 
For generic $r$, the complete reducibility of the tensor product
representation is unknown, and
one can not expect a formula like (3.21). However one can define the coset Virasoro
algebra even in this case by the Goddard-Kent-Olive construction
associated with the tensor
product representation. Its irreducible representation will be
realised in terms of the representation theory
of the elliptic algebra $U_{x,p}(\slth)$ and identified with $\H^{(\ell)}_{m,a}$ in Section 5.2.1.

% For generic $r$, one can not expect a formula like \eqref{VVdecomp}.
\item[2)] In terms of lattice models, \eqref{branch} appears in a consideration of the  
fusion RSOS model with  $r>k+2 \in\Z$ and restricted heights $1\leq a\leq r-1$ 
and $1\leq m\leq r-k-1$\cite{DJMO86b
%,JMO93,Ko98
}. On the other hand,
\mref{charpart} is associated with the fusion 
SOS model with $r (>k+2)$ generic  and with no restriction 
on the local heights. 
\end{itemize}

\subsubsection{HTMs}
The next object to consider is the North half-transfer matrix
$\Phi_{b,a}(u):\cH^{(\ell)}_{m,a}\ra \cH^{(\bell)}_{m,b}$,
which we define graphically by Figure 3.6 a).

\setlength{\unitlength}{0.0003in}
\begingroup\makeatletter\ifx\SetFigFont\undefined%
\gdef\SetFigFont#1#2#3#4#5{%
  \reset@font\fontsize{#1}{#2pt}%
  \fontfamily{#3}\fontseries{#4}\fontshape{#5}%
  \selectfont}%
\fi\endgroup%
{\renewcommand{\dashlinestretch}{30}
\begin{picture}(1173,4575)(-300,900)
\thinlines
\texture{55888888 88555555 5522a222 a2555555 55888888 88555555 552a2a2a 2a555555 
        55888888 88555555 55a222a2 22555555 55888888 88555555 552a2a2a 2a555555 
        55888888 88555555 5522a222 a2555555 55888888 88555555 552a2a2a 2a555555 
        55888888 88555555 55a222a2 22555555 55888888 88555555 552a2a2a 2a555555 }
\put(1050,375){\shade\ellipse{100}{100}}
\put(1050,375){\ellipse{100}{100}}
\put(150,375){\shade\ellipse{100}{100}}
\put(150,375){\ellipse{100}{100}}
\path(150,3975)(1050,3975)(1050,375)
        (150,375)(150,3975)
\path(150,3075)(1050,3075)
\path(150,2175)(1050,2175)
\path(150,1275)(1050,1275)
\path(375,3975)(150,3750)
\path(375,3075)(150,2850)
\path(375,2175)(150,1950)
\path(375,1275)(150,1050)
\put(525,4425){\makebox(0,0)[lb]{{\SetFigFont{12}{16.8}{\rmdefault}{\mddefault}{\updefault}$\vdots$}}}
\put(1050,0){\makebox(0,0)[lb]{{\SetFigFont{12}{16.8}{\rmdefault}{\mddefault}{\updefault}$a$}}}
\put(0,0){\makebox(0,0)[lb]{{\SetFigFont{12}{16.8}{\rmdefault}{\mddefault}{\updefault}$b$}}}
\put(-1000,2200){\makebox(0,0)[lb]{{\SetFigFont{12}{16.8}{\rmdefault}{\mddefault}{\updefault}$a)$}}}
\end{picture}
}

\setlength{\unitlength}{0.0003in}
\begingroup\makeatletter\ifx\SetFigFont\undefined%
\gdef\SetFigFont#1#2#3#4#5{%
  \reset@font\fontsize{#1}{#2pt}%
  \fontfamily{#3}\fontseries{#4}\fontshape{#5}%
  \selectfont}%
\fi\endgroup%
{\renewcommand{\dashlinestretch}{30}
\begin{picture}(9244,4833)(-4000,-3820)
\thinlines
\texture{55888888 88555555 5522a222 a2555555 55888888 88555555 552a2a2a 2a555555 
        55888888 88555555 55a222a2 22555555 55888888 88555555 552a2a2a 2a555555 
        55888888 88555555 5522a222 a2555555 55888888 88555555 552a2a2a 2a555555 
        55888888 88555555 55a222a2 22555555 55888888 88555555 552a2a2a 2a555555 }
\put(4350,4125){\shade\ellipse{100}{100}}
\put(4350,4125){\ellipse{100}{100}}
\put(4350,3225){\shade\ellipse{100}{100}}
\put(4350,3225){\ellipse{100}{100}}
\put(750,1425){\shade\ellipse{100}{100}}
\put(750,1425){\ellipse{100}{100}}
\put(750,525){\shade\ellipse{100}{100}}
\put(750,525){\ellipse{100}{100}}
\put(8250,4125){\shade\ellipse{100}{100}}
\put(8250,4125){\ellipse{100}{100}}
\put(9150,4125){\shade\ellipse{100}{100}}
\put(9150,4125){\ellipse{100}{100}}
\path(750,4125)(4350,4125)(4350,3225)
        (750,3225)(750,4125)
\path(750,1425)(4350,1425)(4350,525)
        (750,525)(750,1425)
\path(8250,4125)(9150,4125)(9150,525)
        (8250,525)(8250,4125)
\path(1650,4125)(1650,3225)
\path(2550,4125)(2550,3225)
\path(3450,4125)(3450,3225)
\path(975,4125)(750,3900)
\path(1875,4125)(1650,3900)
\path(2775,4125)(2550,3900)
\path(1650,1425)(1650,525)
\path(2550,1425)(2550,525)
\path(3450,1425)(3450,525)
\path(975,1425)(750,1200)
\path(1875,1425)(1650,1200)
\path(2775,1425)(2550,1200)
\path(3675,1425)(3450,1200)
\path(8475,4125)(8250,3900)
\path(8250,3225)(9150,3225)
\path(8250,2325)(9150,2325)
\path(8250,1425)(9150,1425)
\path(8475,3225)(8250,3000)
\path(8475,2325)(8250,2100)
\path(8475,1425)(8250,1200)
\path(3675,4125)(3450,3900)
\put(0,3675){\makebox(0,0)[lb]{{\SetFigFont{12}{16.8}{\rmdefault}{\mddefault}{\updefault}$\hdots$}}}
\put(4800,900){\makebox(0,0)[lb]{{\SetFigFont{12}{16.8}{\rmdefault}{\mddefault}{\updefault}$\hdots$}}}
\put(8625,0){\makebox(0,0)[lb]{{\SetFigFont{12}{16.8}{\rmdefault}{\mddefault}{\updefault}$\vdots$}}}
\put(750,4650){\makebox(0,0)[lb]{{\SetFigFont{12}{16.8}{\rmdefault}{\mddefault}{\updefault}$b$)}}}
\put(750,2025){\makebox(0,0)[lb]{{\SetFigFont{12}{16.8}{\rmdefault}{\mddefault}{\updefault}$d$)}}}
\put(7200,2550){\makebox(0,0)[lb]{{\SetFigFont{12}{16.8}{\rmdefault}{\mddefault}{\updefault}$c$)}}}
\put(4725,4050){\makebox(0,0)[lb]{{\SetFigFont{12}{16.8}{\rmdefault}{\mddefault}{\updefault}$a$}}}
\put(4725,3225){\makebox(0,0)[lb]{{\SetFigFont{12}{16.8}{\rmdefault}{\mddefault}{\updefault}$b$}}}
\put(225,450){\makebox(0,0)[lb]{{\SetFigFont{12}{16.8}{\rmdefault}{\mddefault}{\updefault}$a$}}}
\put(225,1425){\makebox(0,0)[lb]{{\SetFigFont{12}{16.8}{\rmdefault}{\mddefault}{\updefault}$a$}}}
\put(8100,4500){\makebox(0,0)[lb]{{\SetFigFont{12}{16.8}{\rmdefault}{\mddefault}{\updefault}$a$}}}
\put(9075,4500){\makebox(0,0)[lb]{{\SetFigFont{12}{16.8}{\rmdefault}{\mddefault}{\updefault}$b$}}}
\end{picture}
}

\vspace*{-35mm}
%\vspace{3mm}
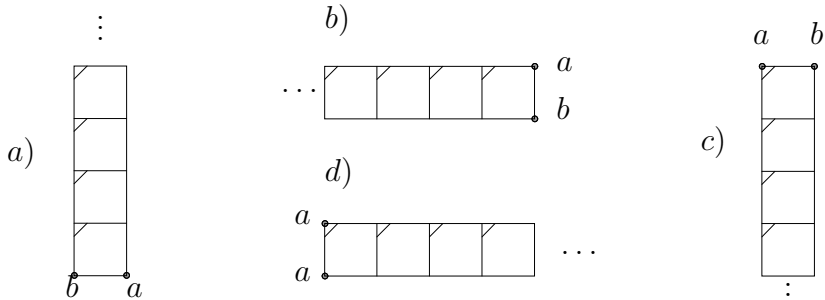
\begin{center}{\ft Figures 3.6 a), b), c) and d) : The North, West, South and East half-transfer matrices}
\end{center}

\nin Again the centre weights are fixed, and the boundary conditions change as
indicated. We view the operator as  acting anti-clockwise about the
centre of the lattice. There are 3 other half-transfer matrices: $\Phi_{W;b,a}(u)$,
$\Phi_{S;b,a}(u)$, and $\Phi_{E;b,a}(u)$,
associated with our SOS model that we might consider, and these are
defined graphically by Figures 3.6 b), c) and d).
\nin These West, South and East half-transfer matrices are again related to
the North half-transfer matrix by crossing symmetry. We find
\ben \Phi_{W;b,a}(u)=(a,b)_k \frac{g_a}{g_{b}}\Phi_{b,a}(-1-u),\\
\Phi_{S;b,a}(u)=(b,a)_k \, \Gamma\, \Phi_{b,a}(u)\,\Gamma^{-1},\\
\Phi_{E;b,a}(u)=\frac{g_a}{g_b}\,\Gamma\, \Phi_{b,a}(-1-u)\,  \Gamma^{-1}.\een

As for vertex models, the graphical arguments 
of \cite{Fodal93}, that rely on the Yang-Baxter equation and unitarity,
lead to the following relations:

\bea
\Phi_{c,b}(u_2)\Phi_{b,a}(u_1)&=&\sum_{d}W^{(k,k)}\left(\left.
\begin{array}{cc}
c&d\\
b&a
\end{array}
\right|u_1-u_2\right)
\Phi_{c,d}(u_1)\Phi_{d,a}(u_2),\lb{commrel}\\
\sum_{a}\Phi^*_{b,a}(u)\Phi_{a,b}(u)&=&
\id,\lb{Eqn:Inversion}\\
A_{b}(u) \Phi_{b,a}(v)&=& \Phi_{b,a}(v-u) A_a(u),\label{Aphi}
\ee
where we have introduced the dual HTM defined by
\be
 \Phi^*_{b,a}(u)&=&\Phi_{W;b,a}(-u)=\frac{1}{G^{(k)}_{b,a}}\,\Phi_{b,a}(u-1).
\en
From \eqref{fctm2h} and \eqref{Aphi}, we also have 
\bea
&&\Phi_{b,a}(u)\ x^{4H^{(\ell)}_{m,a}}=
x^{4H^{(\bell)}_{m,b}}\ \Phi_{b,a}(u-2).\lb{cornerPhi}
\ena

\subsubsection{Correlation functions}
The $N+1$-point correlation function $P^{SOS(\ell)}_m(a,a_1,..,a_N)$ of 
the fusion SOS model is the probability of the local 
height variables 
taking certain values $a,a_1,..,a_N$ on 
specified set of vertices $0,1,,..,N$ of a lattice. Figure 3.7
represents the $N+1$-point function 
$F_m^{SOS(\ell)}(a,a_1,..,a_N)$ divided 
into the parts corresponding to CTMs and HTMs.

\setlength{\unitlength}{0.0003in}
\begingroup\makeatletter\ifx\SetFigFont\undefined%
\gdef\SetFigFont#1#2#3#4#5{%
  \reset@font\fontsize{#1}{#2pt}%
  \fontfamily{#3}\fontseries{#4}\fontshape{#5}%
  \selectfont}%
\fi\endgroup%
{\renewcommand{\dashlinestretch}{30}
\begin{picture}(12573,10530)(-1500,1500)
\thinlines
\path(4725,9525)(4725,5925)
\path(1125,5925)(4650,5925)
\path(2925,8625)(4725,8625)
\path(2025,7725)(4725,7725)
\path(1125,6825)(1125,5925)
\path(2925,8625)(2925,5925)
\path(2025,7725)(2025,5925)
\path(1125,4725)(4650,4725)
\path(4725,4725)(4725,1125)
\path(1125,4725)(1125,3825)
\path(3825,1125)(4725,1125)
\path(2025,2925)(4725,2925)
\path(2925,2025)(4725,2025)
\path(2925,4725)(2925,2025)
\path(2025,4725)(2025,2925)
\path(8625,9525)(8625,5925)
\path(8625,5925)(12225,5925)
\path(8625,9525)(9525,9525)
\path(12225,6825)(12225,5925)
\path(8625,8625)(10425,8625)
\path(8625,7725)(11325,7725)
\path(3825,9525)(4725,9525)
\path(2250,4725)(2025,4500)
\path(3150,4725)(2925,4500)
\path(4050,4725)(3825,4500)
\path(2250,3825)(2025,3600)
\path(3150,3825)(2925,3600)
\path(4050,3825)(3825,3600)
\path(3150,2925)(2925,2700)
\path(4050,2925)(3825,2700)
\path(4050,2025)(3825,1800)
\path(4050,9525)(3825,9300)
\path(3150,8625)(2925,8400)
\path(4050,8625)(3825,8400)
\path(2250,7725)(2025,7500)
\path(3150,7725)(2925,7500)
\path(4050,7725)(3825,7500)
\path(1350,6825)(1125,6600)
\path(2250,6825)(2025,6600)
\path(3150,6825)(2925,6600)
\path(4050,6825)(3825,6600)

\dottedline{68}(5325,5925)(4725,5925)(4725,4725)(5325,4725)
\dottedline{68}(8025,5925)(8625,5925)(8625,
%5625
4725)
\dottedline{68}(8625,4725)(8025,4725)
%(8625,5025)
\path(10425,8625)(10425,5925)
\put(12450,4275){\makebox(0,0)[lb]{\smash{{{\SetFigFont{8}{16.8}{\rmdefault}{\mddefault}{\updefault}$\hdots$}}}}}
\put(4925,5250){\makebox(0,0)[lb]{\smash{{{\SetFigFont{8}{16.8}{\rmdefault}{\mddefault}{\updefault}$a$}}}}}
\put(6075,5250){\makebox(0,0)[lb]{\smash{{{\SetFigFont{8}{16.8}{\rmdefault}{\mddefault}{\updefault}$a_1$}}}}}
\put(6675,5250){\makebox(0,0)[lb]{\smash{{{\SetFigFont{8}{16.8}{\rmdefault}{\mddefault}{\updefault}$a_{N-1}$}}}}}
\put(8075,5250){\makebox(0,0)[lb]{\smash{{{\SetFigFont{8}{16.8}{\rmdefault}{\mddefault}{\updefault}$a_N$}}}}}
%\put(6075,4875){\makebox(0,0)[lb]{\smash{{{\SetFigFont{8}{16.8}{\rmdefault}{\mddefault}{\updefault}$a_1$}}}}}
%\put(6675,4875){\makebox(0,0)[lb]{\smash{{{\SetFigFont{8}{16.8}{\rmdefault}{\mddefault}{\updefault}$a_{N-1}$}}}}}
%\put(7875,4875){\makebox(0,0)[lb]{\smash{{{\SetFigFont{8}{16.8}{\rmdefault}{\mddefault}{\updefault}$a_N$}}}}}

%%%%%%%%%%%%%%%%%%%%%%%%%%%%%%%%%%%%%%%%%%%%%
%\put(5375,5250){\makebox(0,0)[lb]{\smash{{{\SetFigFont{8}{16.8}{\rmdefault}{\mddefault}{\updefault}$\vep_1$}}}}}
%\put(6300,5250){\makebox(0,0)[lb]{\smash{{{\SetFigFont{8}{16.8}{\rmdefault}{\mddefault}{\updefault}$\vep_2$}}}}}
%\put(7075,5250){\makebox(0,0)[lb]{\smash{{{\SetFigFont{8}{16.8}{\rmdefault}{\mddefault}{\updefault}$\vep_N$}}}}}
%\put(-800,5400){\makebox(0,0)[lb]{\smash{{{\SetFigFont{16}{16.8}{\rmdefault}{\mddefault}{\updefault}$\sum\limits_{a, a_j,a'_j}$}}}}}
%\put(13400,5200){\makebox(0,0)[lb]{\smash{{{\SetFigFont{12}{16.8}{\rmdefault}{\mddefault}{\updefault}$\widetilde{\Lambda}_{a_N,a'_N}(u)$}}}}}

%%%%%%%%%%%%%%%%%%%%%%%%%%%%%%%%%%%%%%%%%%%%%%%%%%%%%%%%%
\put(1050,8700){\makebox(0,0)[lb]{\smash{{{\SetFigFont{12}{16.8}{\rmdefault}{\mddefault}{\updefault}$A_a(u)$}}}}}
\put(11250,8700){\makebox(0,0)[lb]{\smash{{{\SetFigFont{12}{16.8}{\rmdefault}{\mddefault}{\updefault}$A_{NE;\,a_N}(u)$}}}}}
\put(550,1900){\makebox(0,0)[lb]{\smash{{{\SetFigFont{12}{16.8}{\rmdefault}{\mddefault}{\updefault}$A_{SW;\,a}(u)$}}}}}
\put(11550,1900){\makebox(0,0)[lb]{\smash{{{\SetFigFont{12}{16.8}{\rmdefault}{\mddefault}{\updefault}$A_{SE;\,a_N}(u)$}}}}}
\put(4825,10850){\makebox(0,0)[lb]{\smash{{{\SetFigFont{8}{16.8}{\rmdefault}{\mddefault}{\updefault}$\Phi_{a,\,a_1}(u)$}}}}}
\put(6425,10850){\makebox(0,0)[lb]{\smash{{{\SetFigFont{8}{16.8}{\rmdefault}{\mddefault}{\updefault}$\hdots$}}}}}
\put(7025,10850){\makebox(0,0)[lb]{\smash{{{\SetFigFont{8}{16.8}{\rmdefault}{\mddefault}{\updefault}$\Phi_{a_{N-1},a_N}(u)$}}}}}
\put(4625,-200){\makebox(0,0)[lb]{\smash{{{\SetFigFont{8}{16.8}{\rmdefault}{\mddefault}{\updefault}$\Phi_{S;a_1,a}(u)$}}}}}
\put(6425,9000){\makebox(0,0)[lb]{\smash{{{\SetFigFont{8}{16.8}{\rmdefault}{\mddefault}{\updefault}$\hdots$}}}}}
\put(6425,8100){\makebox(0,0)[lb]{\smash{{{\SetFigFont{8}{16.8}{\rmdefault}{\mddefault}{\updefault}$\hdots$}}}}}
\put(6425,7200){\makebox(0,0)[lb]{\smash{{{\SetFigFont{8}{16.8}{\rmdefault}{\mddefault}{\updefault}$\hdots$}}}}}
\put(6425,6300){\makebox(0,0)[lb]{\smash{{{\SetFigFont{8}{16.8}{\rmdefault}{\mddefault}{\updefault}$\hdots$}}}}}
\put(6425,4300){\makebox(0,0)[lb]{\smash{{{\SetFigFont{8}{16.8}{\rmdefault}{\mddefault}{\updefault}$\hdots$}}}}}
\put(6425,3400){\makebox(0,0)[lb]{\smash{{{\SetFigFont{8}{16.8}{\rmdefault}{\mddefault}{\updefault}$\hdots$}}}}}
\put(6425,2500){\makebox(0,0)[lb]{\smash{{{\SetFigFont{8}{16.8}{\rmdefault}{\mddefault}{\updefault}$\hdots$}}}}}
\put(6425,1600){\makebox(0,0)[lb]{\smash{{{\SetFigFont{8}{16.8}{\rmdefault}{\mddefault}{\updefault}$\hdots$}}}}}
\put(6425,-200){\makebox(0,0)[lb]{\smash{{{\SetFigFont{8}{16.8}{\rmdefault}{\mddefault}{\updefault}$\hdots$}}}}}

\put(7025,-200){\makebox(0,0)[lb]{\smash{{{\SetFigFont{8}{16.8}{\rmdefault}{\mddefault}{\updefault}$\Phi_{S;a_N,a_{N-1}}(u)$}}}}}
%\put(12450,5325){\makebox(0,0)[lb]{\smash{{{\SetFigFont{8}{16.8}{\rmdefault}{\mddefault}{\updefault}$\hdots$}}}}}
\path(11325,7725)(11325,5925)
\path(8700,4725)(12225,4725)
\path(8625,4725)(8625,1125)
\path(8625,1125)(9525,1125)
\path(12225,4725)(12225,3825)
\path(11325,4725)(11325,2925)
\path(8625,2025)(10425,2025)
\path(8625,2925)(11325,2925)
\path(10425,4725)(10425,2025)
\put(4125,9900){\makebox(0,0)[lb]{\smash{{{\SetFigFont{8}{16.8}{\rmdefault}{\mddefault}{\updefault}$\vdots$}}}}}
\put(5700,9900){\makebox(0,0)[lb]{\smash{{{\SetFigFont{8}{16.8}{\rmdefault}{\mddefault}{\updefault}$\vdots$}}}}}
%\put(6525,9900){\makebox(0,0)[lb]{\smash{{{\SetFigFont{8}{16.8}{\rmdefault}{\mddefault}{\updefault}$\vdots$}}}}}
\put(7500,9900){\makebox(0,0)[lb]{\smash{{{\SetFigFont{8}{16.8}{\rmdefault}{\mddefault}{\updefault}$\vdots$}}}}}
\put(8925,9900){\makebox(0,0)[lb]{\smash{{{\SetFigFont{8}{16.8}{\rmdefault}{\mddefault}{\updefault}$\vdots$}}}}}
\put(4125,650){\makebox(0,0)[lb]{\smash{{{\SetFigFont{8}{16.8}{\rmdefault}{\mddefault}{\updefault}$\vdots$}}}}}
\put(5700,650){\makebox(0,0)[lb]{\smash{{{\SetFigFont{8}{16.8}{\rmdefault}{\mddefault}{\updefault}$\vdots$}}}}}
%\put(6600,650){\makebox(0,0)[lb]{\smash{{{\SetFigFont{8}{16.8}{\rmdefault}{\mddefault}{\updefault}$\vdots$}}}}}
\put(7500,650){\makebox(0,0)[lb]{\smash{{{\SetFigFont{8}{16.8}{\rmdefault}{\mddefault}{\updefault}$\vdots$}}}}}
\put(8925,650){\makebox(0,0)[lb]{\smash{{{\SetFigFont{8}{16.8}{\rmdefault}{\mddefault}{\updefault}$\vdots$}}}}}
\put(525,6450){\makebox(0,0)[lb]{\smash{{{\SetFigFont{8}{16.8}{\rmdefault}{\mddefault}{\updefault}$\hdots$}}}}}
\put(525,4275){\makebox(0,0)[lb]{\smash{{{\SetFigFont{8}{16.8}{\rmdefault}{\mddefault}{\updefault}$\hdots$}}}}}
\put(12450,6375){\makebox(0,0)[lb]{\smash{{{\SetFigFont{8}{16.8}{\rmdefault}{\mddefault}{\updefault}$\hdots$}}}}}
\path(1350,4725)(1125,4500)
\path(7125,9525)(7125,5925)
\path(9525,9525)(9525,5925)
\path(1125,6825)(4725,6825)
%\path(5325,8625)(8025,8625)
%\path(5325,7725)(8025,7725)
%\path(5325,6825)(8025,6825)
\path(8625,6825)(12225,6825)
\path(1125,3825)(4725,3825)
%\path(5325,3825)(8025,3825)
%\path(5325,2925)(8025,2925)
%\path(5325,2025)(8025,2025)
%%%%% new %%%%
\path(5325,4725)(6225,4725)
\path(5325,3825)(6225,3825)
\path(5325,2925)(6225,2925)
\path(5325,2025)(6225,2025)
\path(5325,1125)(6225,1125)

\path(7125,4725)(8025,4725)
\path(7125,3825)(8025,3825)
\path(7125,2925)(8025,2925)
\path(7125,2025)(8025,2025)
\path(7125,1125)(8025,1125)

\path(5325,9525)(6225,9525)
\path(5325,8625)(6225,8625)
\path(5325,7725)(6225,7725)
\path(5325,6825)(6225,6825)
\path(5325,5925)(6225,5925)

\path(7125,9525)(8025,9525)
\path(7125,8625)(8025,8625)
\path(7125,7725)(8025,7725)
\path(7125,6825)(8025,6825)
\path(7125,5925)(8025,5925)

\path(5325,5925)(5325,9525)
\path(5325,1125)(5325,4725)

\path(8025,5925)(8025,9525)
\path(8025,1125)(8025,4725)

%%%%%%%%%%%%%%
\path(8625,3825)(12225,3825)
\path(3825,4725)(3825,1125)
\path(6225,4725)(6225,1125)
\path(7125,4725)(7125,1125)
%%%%%%%%%%%%%%%%%%%%%%%%%%%%%%%%%%
%L_vep lines  
%\path(5775,5925)(5775,4725)
%\path(6675,5925)(6675,4725)
%\path(7575,5925)(7575,4725)

%%%%%%%%%%%%%%%%%%%%%%%%%%%%%%%%%%%%
%Tail verttical lines
%\path(9075,5625)(9075,5025)
%\path(9975,5625)(9975,5025)
%\path(10875,5625)(10875,5025)
%\path(11775,5625)(11775,5025)

%Tail horizontal line
%\path(8625,5625)(12225,5625)
%\path(8625,5025)(12225,5025)

%\put(8625,5025){\shade\ellipse{100}{100}}
%\put(8625,5025){\ellipse{100}{100}}

%\put(8625,5625){\shade\ellipse{100}{100}}
%\put(8625,5625){\ellipse{100}{100}}

%%%%%%%%%%%%%%%%%%%%%%%%%%%%%%
\path(6225,9525)(6225,5925)
\texture{55888888 88555555 5522a222 a2555555 55888888 88555555 552a2a2a 2a555555 
        55888888 88555555 55a222a2 22555555 55888888 88555555 552a2a2a 2a555555 
        55888888 88555555 5522a222 a2555555 55888888 88555555 552a2a2a 2a555555 
        55888888 88555555 55a222a2 22555555 55888888 88555555 552a2a2a 2a555555 }
%\put(5775,5325){\shade\ellipse{100}{100}}
%\put(5775,5325){\ellipse{100}{100}}

%\put(6675,5325){\shade\ellipse{100}{100}}
%\put(6675,5325){\ellipse{100}{100}}

%\put(7575,5325){\shade\ellipse{100}{100}}
%\put(7575,5325){\ellipse{100}{100}}

\put(5325,5925){\shade\ellipse{100}{100}}
\put(5325,5925){\ellipse{100}{100}}

\put(4725,5925){\shade\ellipse{100}{100}}
\put(4725,5925){\ellipse{100}{100}}
\put(4725,4725){\shade\ellipse{100}{100}}
\put(4725,4725){\ellipse{100}{100}}
\put(5325,4725){\shade\ellipse{100}{100}}
\put(5325,4725){\ellipse{100}{100}}
\put(8025,5925){\shade\ellipse{100}{100}}
\put(8025,5925){\ellipse{100}{100}}
\put(7125,5925){\shade\ellipse{100}{100}}
\put(7125,5925){\ellipse{100}{100}}
\put(6225,5925){\shade\ellipse{100}{100}}
\put(6225,5925){\ellipse{100}{100}}
\put(6225,4725){\shade\ellipse{100}{100}}
\put(6225,4725){\ellipse{100}{100}}
\put(7125,4725){\shade\ellipse{100}{100}}
\put(7125,4725){\ellipse{100}{100}}
\put(8025,4725){\shade\ellipse{100}{100}}
\put(8025,4725){\ellipse{100}{100}}

\put(8625,4725){\shade\ellipse{100}{100}}
\put(8625,4725){\ellipse{100}{100}}

\put(8625,5925){\shade\ellipse{100}{100}}
\put(8625,5925){\ellipse{100}{100}}

%\path(8025,9525)(5325,9525)(5325,5925)
%        (8025,5925)(8025,9525)
%\path(5325,4725)(8025,4725)(8025,1125)
%        (5325,1125)(5325,4725)

\path(3825,9525)(3825,5925)

\path(8850,4725)(8625,4500)
\path(9750,4725)(9525,4500)
\path(10650,4725)(10425,4500)
\path(11550,4725)(11325,4500)
\path(8850,3825)(8625,3600)
\path(9750,3825)(9525,3600)
\path(10650,3825)(10425,3600)
\path(8850,2925)(8625,2700)
\path(9750,2925)(9525,2700)
\path(8850,2025)(8625,1800)
\path(5550,4725)(5325,4500)
%\path(6450,4725)(6225,4500)
\path(7350,4725)(7125,4500)
\path(5550,3825)(5325,3600)
%\path(6450,3825)(6225,3600)
\path(7350,3825)(7125,3600)
\path(5550,2925)(5325,2700)
%\path(6450,2925)(6225,2700)
\path(7350,2925)(7125,2700)
\path(5550,2025)(5325,1800)
%\path(6450,2025)(6225,1800)
\path(7350,2025)(7125,1800)
\path(11550,6825)(11325,6600)
\path(9525,4725)(9525,1125)
\path(9525,4725)(9525,1125)
\path(5550,9525)(5325,9300)
\path(5550,9525)(5325,9300)
%\path(6450,9525)(6225,9300)
%\path(6450,9525)(6225,9300)
\path(7350,9525)(7125,9300)
\path(7350,9525)(7125,9300)
\path(5550,8625)(5325,8400)
%\path(6450,8625)(6225,8400)
\path(7350,8625)(7125,8400)
\path(5550,7725)(5325,7500)
%\path(6450,7725)(6225,7500)
\path(7350,7725)(7125,7500)
\path(5550,6825)(5325,6600)
%\path(6450,6825)(6225,6600)
\path(7350,6825)(7125,6600)
\path(8850,9525)(8625,9300)
\path(8850,8625)(8625,8400)
\path(9750,8625)(9525,8400)
\path(8850,7725)(8625,7500)
\path(9750,7725)(9525,7500)
\path(10650,7725)(10425,7500)
\path(8850,6825)(8625,6600)
\path(9750,6825)(9525,6600)
\path(10650,6825)(10425,6600)

% vert arrows
%\path(5725,5025)(5775,4875)(5825,5025)
%\path(6625,5025)(6675,4875)(6725,5025)
%\path(7525,5025)(7575,4875)(7625,5025)

%\path(9025,5425)(9075,5275)(9125,5425)
%\path(9925,5425)(9975,5275)(10025,5425)
%\path(10825,5425)(10875,5275)(10925,5425)
%\path(11725,5425)(11775,5275)(11825,5425)
\end{picture}
}

\vspace*{10mm}\begin{center}{\ft Figure 3.7: The $N+1$-point function $F_m^{SOS(\ell)}(a,a_1,\cdots,a_N)$}
\end{center}
\noi 
According to this picture we have 
\be
P^{SOS(\ell)}_m(a,a_1,..,a_N)&=&\frac{1}{Z^{(\ell)}_m}
{F_m^{SOS(\ell)}(a,a_1,..,a_N)},\\
F_m^{SOS(\ell)}(a,a_1,..,a_N)
&=&\Tr_{\cH^{(\ell)}_{m,a}}\big(
\Phi_{a,a_1}(u)\cdots \Phi_{a_{N-1},a_N}(u) A_{NE;a_N}(u) A_{SE;a_N}(u)\nn\\ 
&&\qquad \times
\Phi_{S;a_N,a_{N-1}}(u) \cdots \Phi_{S;a_1,a}(u) A_{SW;a}(u)A_{a}(u) \big).
\en
Using the relations \eqref{fctm2h}, \eqref{Aphi} and \eqref{cornerPhi}, we obtain the simplified 
expression. 
\bea
&&\hspace*{-5mm} P^{SOS(\ell)}_m(a,a_1,..,a_N)\nn\\
&&\hspace*{-5mm}=\frac{[a_N]}{
%Z^{(\ell)}_m
[m]^*\chi^{(k)}_\ell(\bar{\tau})}\,
Tr_{\cH^{(\sigma^N(\ell))}_{m,a_N}}\big(x^{4H^{(\ell)}_{m,a}}
\Phi^*_{a_N,a_{N-1}}(u)\cdots \Phi^*_{a_1,a}(u)\Phi_{a,a_1}(u) \cdots
\Phi_{a_{N-1},a_N}(u)\big).
\ena

%%%%%%%%%%%%%%%%%%%%%%%%%%%%%%%%%%%%%%%%%%%%%%%%%%%%%%%%%%%%%%%%%%%%%%%%%%%%%
\section{The Vertex-Face Correspondence}
%%%%%%%%%%%%%%%%%%%%%%%%%%%%%%%%%%%%%%%%%%%%%%%%%%%%%%%%%%%%%%%%%%%%%%%%%%%%%%
While the algebraic analysis approach works well for many models, including fusion
six-vertex and fusion SOS models \cite{idzal93,Idz94,BoWe94b,JMO93,Ko98}, it runs into a technical obstacle
for the case of the fusion eight-vertex model. The difficulty is the lack
of a suitable free-field realization to evaluate the trace 
occurring in \mref{vtrace}.
In order to overcome this problem, we shall follow Baxter \cite{Bax73aII} 
and Lashkevich
and Pugai \cite{LaP98} and relate the fusion eight-vertex 
model to the fusion SOS model. This latter model
is more tractable from the point of view of the algebraic analysis 
approach\cite{LP96,Ko98,JKOS}.

\subsection{Dressed Vertex Models}
The starting point in establishing the connection of the expression \mref{vtrace} 
with SOS models is to dress the
 boundary of the vertex model defined on the finite $(2L+N)\times 2L$ lattice 
shown in Figure 3.1 with the 
intertwining and dual intertwining vectors expressed 
by 3-vertices in Figure 2.5 $(a)$ and $(b)$.
The procedure is shown in Figure 4.1. In this section, we
identify the vertex model spectral parameter $u$ as 
$u=u_0-v$, where $u_0$ and $v$ are the spectral parameters
associated with vertical and horizontal lines respectively.

\begin{scriptsize}
\setlength{\unitlength}{0.0004in}
\begingroup\makeatletter\ifx\SetFigFont\undefined%
\gdef\SetFigFont#1#2#3#4#5{%
  \reset@font\fontsize{#1}{#2pt}%
  \fontfamily{#3}\fontseries{#4}\fontshape{#5}%
  \selectfont}%
\fi\endgroup%
{\renewcommand{\dashlinestretch}{30}
\begin{picture}(11607,8055)(-2000,600)
\thinlines
\texture{55888888 88555555 5522a222 a2555555 55888888 88555555 552a2a2a 2a555555 
        55888888 88555555 55a222a2 22555555 55888888 88555555 552a2a2a 2a555555 
        55888888 88555555 5522a222 a2555555 55888888 88555555 552a2a2a 2a555555 
        55888888 88555555 55a222a2 22555555 55888888 88555555 552a2a2a 2a555555 }
\put(5850,3975){\shade\ellipse{100}{100}}
\put(5850,3975){\ellipse{100}{100}}
\put(4950,3975){\shade\ellipse{100}{100}}
\put(4950,3975){\ellipse{100}{100}}
\put(1800,7575){\shade\ellipse{100}{100}}
\put(1800,7575){\ellipse{100}{100}}
\put(2700,7575){\shade\ellipse{100}{100}}
\put(2700,7575){\ellipse{100}{100}}
\put(3600,7575){\shade\ellipse{100}{100}}
\put(3600,7575){\ellipse{100}{100}}
\put(4500,7575){\shade\ellipse{100}{100}}
\put(4500,7575){\ellipse{100}{100}}
\put(5400,7575){\shade\ellipse{100}{100}}
\put(5400,7575){\ellipse{100}{100}}
\put(6300,7575){\shade\ellipse{100}{100}}
\put(6300,7575){\ellipse{100}{100}}
\put(7200,7575){\shade\ellipse{100}{100}}
\put(7200,7575){\ellipse{100}{100}}
\put(8100,7575){\shade\ellipse{100}{100}}
\put(8100,7575){\ellipse{100}{100}}
\put(6750,3975){\shade\ellipse{100}{100}}
\put(6750,3975){\ellipse{100}{100}}
\put(9000,7575){\shade\ellipse{100}{100}}
\put(9000,7575){\ellipse{100}{100}}
\put(9900,7575){\shade\ellipse{100}{100}}
\put(9900,7575){\ellipse{100}{100}}
\put(9900,375){\shade\ellipse{100}{100}}
\put(9900,375){\ellipse{100}{100}}
\put(1800,375){\shade\ellipse{100}{100}}
\put(1800,375){\ellipse{100}{100}}
\put(2700,375){\shade\ellipse{100}{100}}
\put(2700,375){\ellipse{100}{100}}
\put(3600,375){\shade\ellipse{100}{100}}
\put(3600,375){\ellipse{100}{100}}
\put(4500,375){\shade\ellipse{100}{100}}
\put(4500,375){\ellipse{100}{100}}
\put(5400,375){\shade\ellipse{100}{100}}
\put(5400,375){\ellipse{100}{100}}
\put(6300,375){\shade\ellipse{100}{100}}
\put(6300,375){\ellipse{100}{100}}
\put(7200,375){\shade\ellipse{100}{100}}
\put(7200,375){\ellipse{100}{100}}
\put(8100,375){\shade\ellipse{100}{100}}
\put(8100,375){\ellipse{100}{100}}
\put(9000,375){\shade\ellipse{100}{100}}
\put(9000,375){\ellipse{100}{100}}
\put(900,7575){\shade\ellipse{100}{100}}
\put(900,7575){\ellipse{100}{100}}
\put(900,6675){\shade\ellipse{100}{100}}
\put(900,6675){\ellipse{100}{100}}
\put(900,5775){\shade\ellipse{100}{100}}
\put(900,5775){\ellipse{100}{100}}
\put(900,4875){\shade\ellipse{100}{100}}
\put(900,4875){\ellipse{100}{100}}
\put(900,3975){\shade\ellipse{100}{100}}
\put(900,3975){\ellipse{100}{100}}
\put(900,3075){\shade\ellipse{100}{100}}
\put(900,3075){\ellipse{100}{100}}
\put(900,2175){\shade\ellipse{100}{100}}
\put(900,2175){\ellipse{100}{100}}
\put(900,1275){\shade\ellipse{100}{100}}
\put(900,1275){\ellipse{100}{100}}
\put(900,375){\shade\ellipse{100}{100}}
\put(900,375){\ellipse{100}{100}}
\put(10800,7575){\shade\ellipse{100}{100}}
\put(10800,7575){\ellipse{100}{100}}
\put(10800,6675){\shade\ellipse{100}{100}}
\put(10800,6675){\ellipse{100}{100}}
\put(10800,5775){\shade\ellipse{100}{100}}
\put(10800,5775){\ellipse{100}{100}}
\put(10800,4875){\shade\ellipse{100}{100}}
\put(10800,4875){\ellipse{100}{100}}
\put(10800,3975){\shade\ellipse{100}{100}}
\put(10800,3975){\ellipse{100}{100}}
\put(10800,3075){\shade\ellipse{100}{100}}
\put(10800,3075){\ellipse{100}{100}}
\put(10800,2175){\shade\ellipse{100}{100}}
\put(10800,2175){\ellipse{100}{100}}
\put(10800,1275){\shade\ellipse{100}{100}}
\put(10800,1275){\ellipse{100}{100}}
\put(10800,375){\shade\ellipse{100}{100}}
\put(10800,375){\ellipse{100}{100}}
\path(2250,7575)(2250,1275)
\path(1800,7575)(9000,7575)
\path(1800,6225)(9000,6225)
\path(1800,5325)(9000,5325)
\path(1800,4425)(9000,4425)
\path(1800,2625)(9000,2625)
\path(1800,7125)(9000,7125)
\path(1800,3525)(9000,3525)
\path(3150,7575)(3150,1275)
\path(4050,7575)(4050,1275)
\path(4950,7575)(4950,1275)
\path(5850,7575)(5850,1275)
\path(6750,7575)(6750,1275)
\path(7650,7575)(7650,1275)
\path(8550,7575)(8550,1275)
\path(1800,1725)(9000,1725)
\path(9900,7575)(9000,7575)
\path(9900,7575)(9000,7575)
\path(9450,7575)(9450,1275)
\path(9450,7575)(9450,1275)
\path(9000,7125)(9900,7125)
\path(9000,7125)(9900,7125)
\path(9000,6225)(9900,6225)
\path(9000,6225)(9900,6225)
\path(9000,5325)(9900,5325)
\path(9000,5325)(9900,5325)
\path(9000,4425)(9900,4425)
\path(9000,4425)(9900,4425)
\path(9000,3525)(9900,3525)
\path(9000,3525)(9900,3525)
\path(9000,2625)(9900,2625)
\path(9000,2625)(9900,2625)
\path(9000,1725)(9900,1725)
\path(9000,1725)(9900,1725)
\path(1800,825)(9900,825)
\path(1800,825)(9900,825)
\path(1800,375)(9900,375)
\path(1800,375)(9900,375)
\path(2250,1275)(2250,375)
\path(2250,1275)(2250,375)
\path(3150,1350)(3150,375)
\path(3150,1350)(3150,375)
\path(4050,1275)(4050,375)
\path(4050,1275)(4050,375)
\path(4950,1350)(4950,375)
\path(4950,1350)(4950,375)
\path(5850,1275)(5850,375)
\path(5850,1275)(5850,375)
\path(6750,1275)(6750,375)
\path(6750,1275)(6750,375)
\path(7650,1275)(7650,375)
\path(7650,1275)(7650,375)
\path(8550,1275)(8550,375)
\path(8550,1275)(8550,375)
\path(9450,1275)(9450,375)
\path(9450,1275)(9450,375)
\path(10800,7575)(10800,375)
\path(9900,7125)(10800,7125)
\path(9900,6225)(10800,6225)
\path(9900,5325)(10800,5325)
\path(9900,4425)(10800,4425)
\path(9900,3525)(10800,3525)
\path(9900,2625)(10800,2625)
\path(9900,1725)(10800,1725)
\path(9900,825)(10800,825)
\path(9900,375)(10800,375)
\path(9900,7575)(10800,7575)
\path(10350,7575)(10350,375)
\path(900,7575)(1800,7575)
\path(900,7125)(1800,7125)
\path(900,6225)(1950,6225)
\path(900,5325)(1800,5325)
\path(900,4425)(1800,4425)
\path(900,3525)(1875,3525)
\path(900,2625)(1800,2625)
\path(900,1725)(1800,1725)
\path(900,825)(1875,825)
\path(900,375)(1800,375)
\path(900,7575)(900,375)
\path(1350,7575)(1350,375)
%%% horiz arrows
\path(1300,875)(1100,825)(1300,775)
\path(1300,1775)(1100,1725)(1300,1675)
\path(1300,2675)(1100,2625)(1300,2575)
\path(1300,3575)(1100,3525)(1300,3475)
\path(1300,4475)(1100,4425)(1300,4375)
\path(1300,5375)(1100,5325)(1300,5275)
\path(1300,6275)(1100,6225)(1300,6175)
\path(1300,7175)(1100,7125)(1300,7075)
%%% vert arrows
\path(1300,650)(1350,450)(1400,650)
\path(2200,650)(2250,450)(2300,650)
\path(3100,650)(3150,450)(3200,650)
\path(4000,650)(4050,450)(4100,650)
\path(4900,650)(4950,450)(5000,650)
\path(5800,650)(5850,450)(5900,650)
\path(6700,650)(6750,450)(6800,650)
\path(7600,650)(7650,450)(7700,650)
\path(8500,650)(8550,450)(8600,650)
\path(9400,650)(9450,450)(9500,650)
\path(10300,650)(10350,450)(10400,650)

%%%
\put(4500,3900){\makebox(0,0)[lb]{{\SetFigFont{10}{16.8}{\rmdefault}{\mddefault}{\updefault}$\vep_1$}}}
\put(5400,3900){\makebox(0,0)[lb]{{\SetFigFont{10}{16.8}{\rmdefault}{\mddefault}{\updefault}$\vep_2$}}}
\put(6300,3900){\makebox(0,0)[lb]{{\SetFigFont{10}{16.8}{\rmdefault}{\mddefault}{\updefault}$\vep_3$}}}
\put(4200,7875){\makebox(0,0)[lb]{{\SetFigFont{8}{16.8}{\rmdefault}{\mddefault}{\updefault}$m+\ell$}}}
\put(5175,7875){\makebox(0,0)[lb]{{\SetFigFont{8}{16.8}{\rmdefault}{\mddefault}{\updefault}$m+\bell$}}}
\put(6150,7875){\makebox(0,0)[lb]{{\SetFigFont{8}{16.8}{\rmdefault}{\mddefault}{\updefault}$m+\ell$}}}
\put(6975,7875){\makebox(0,0)[lb]{{\SetFigFont{8}{16.8}{\rmdefault}{\mddefault}{\updefault}$m+\bell$}}}
\put(7950,7875){\makebox(0,0)[lb]{{\SetFigFont{8}{16.8}{\rmdefault}{\mddefault}{\updefault}$m+\ell$}}}
\put(8850,7875){\makebox(0,0)[lb]{{\SetFigFont{8}{16.8}{\rmdefault}{\mddefault}{\updefault}$m+\bell$}}}
\put(9750,7875){\makebox(0,0)[lb]{{\SetFigFont{8}{16.8}{\rmdefault}{\mddefault}{\updefault}$m+\ell$}}}
\put(3375,7875){\makebox(0,0)[lb]{{\SetFigFont{8}{16.8}{\rmdefault}{\mddefault}{\updefault}$m+\bell$}}}
\put(2550,7875){\makebox(0,0)[lb]{{\SetFigFont{8}{16.8}{\rmdefault}{\mddefault}{\updefault}$m+\ell$}}}
\put(1650,7875){\makebox(0,0)[lb]{{\SetFigFont{8}{16.8}{\rmdefault}{\mddefault}{\updefault}$m+\bell$}}}
\put(750,7875){\makebox(0,0)[lb]{{\SetFigFont{8}{16.8}{\rmdefault}{\mddefault}{\updefault}$m+\ell$}}}
\put(10650,7875){\makebox(0,0)[lb]{{\SetFigFont{8}{16.8}{\rmdefault}{\mddefault}{\updefault}$m+\bell$}}}
\put(11100,6675){\makebox(0,0)[lb]{{\SetFigFont{8}{16.8}{\rmdefault}{\mddefault}{\updefault}$m+\ell$}}}
\put(11100,5775){\makebox(0,0)[lb]{{\SetFigFont{8}{16.8}{\rmdefault}{\mddefault}{\updefault}$m+\bell$}}}
\put(11100,4875){\makebox(0,0)[lb]{{\SetFigFont{8}{16.8}{\rmdefault}{\mddefault}{\updefault}$m+\ell$}}}
\put(11100,3975){\makebox(0,0)[lb]{{\SetFigFont{8}{16.8}{\rmdefault}{\mddefault}{\updefault}$m+\bell$}}}
\put(11100,3075){\makebox(0,0)[lb]{{\SetFigFont{8}{16.8}{\rmdefault}{\mddefault}{\updefault}$m+\ell$}}}
\put(11100,2175){\makebox(0,0)[lb]{{\SetFigFont{8}{16.8}{\rmdefault}{\mddefault}{\updefault}$m+\bell$}}}
%\put(11100,1575){\makebox(0,0)[lb]{{\SetFigFont{8}{16.8}{\rmdefault}{\mddefault}{\updefault}$m+\ell$}}}
\put(11100,1275){\makebox(0,0)[lb]{{\SetFigFont{8}{16.8}{\rmdefault}{\mddefault}{\updefault}$m+\ell$}}}
\put(11100,375){\makebox(0,0)[lb]{{\SetFigFont{8}{16.8}{\rmdefault}{\mddefault}{\updefault}$m+\bell$}}}
\put(0,6675){\makebox(0,0)[lb]{{\SetFigFont{8}{16.8}{\rmdefault}{\mddefault}{\updefault}$m+\bell$}}}
\put(0,5775){\makebox(0,0)[lb]{{\SetFigFont{8}{16.8}{\rmdefault}{\mddefault}{\updefault}$m+\ell$}}}
\put(0,4875){\makebox(0,0)[lb]{{\SetFigFont{8}{16.8}{\rmdefault}{\mddefault}{\updefault}$m+\bell$}}}
\put(0,3975){\makebox(0,0)[lb]{{\SetFigFont{8}{16.8}{\rmdefault}{\mddefault}{\updefault}$m+\ell$}}}
\put(0,3075){\makebox(0,0)[lb]{{\SetFigFont{8}{16.8}{\rmdefault}{\mddefault}{\updefault}$m+\bell$}}}
\put(0,2175){\makebox(0,0)[lb]{{\SetFigFont{8}{16.8}{\rmdefault}{\mddefault}{\updefault}$m+\ell$}}}
\put(0,1275){\makebox(0,0)[lb]{{\SetFigFont{8}{16.8}{\rmdefault}{\mddefault}{\updefault}$m+\bell$}}}
\put(600,0){\makebox(0,0)[lb]{{\SetFigFont{8}{16.8}{\rmdefault}{\mddefault}{\updefault}$m+\ell$}}}
\put(1575,0){\makebox(0,0)[lb]{{\SetFigFont{8}{16.8}{\rmdefault}{\mddefault}{\updefault}$m+\bell$}}}
\put(2475,0){\makebox(0,0)[lb]{{\SetFigFont{8}{16.8}{\rmdefault}{\mddefault}{\updefault}$m+\ell$}}}
\put(3375,0){\makebox(0,0)[lb]{{\SetFigFont{8}{16.8}{\rmdefault}{\mddefault}{\updefault}$m+\bell$}}}
\put(4275,0){\makebox(0,0)[lb]{{\SetFigFont{8}{16.8}{\rmdefault}{\mddefault}{\updefault}$m+\ell$}}}
\put(5175,0){\makebox(0,0)[lb]{{\SetFigFont{8}{16.8}{\rmdefault}{\mddefault}{\updefault}$m+\bell$}}}
\put(6000,0){\makebox(0,0)[lb]{{\SetFigFont{8}{16.8}{\rmdefault}{\mddefault}{\updefault}$m+\ell$}}}
\put(6975,0){\makebox(0,0)[lb]{{\SetFigFont{8}{16.8}{\rmdefault}{\mddefault}{\updefault}$m+\bell$}}}
\put(7875,0){\makebox(0,0)[lb]{{\SetFigFont{8}{16.8}{\rmdefault}{\mddefault}{\updefault}$m+\ell$}}}
\put(8775,0){\makebox(0,0)[lb]{{\SetFigFont{8}{16.8}{\rmdefault}{\mddefault}{\updefault}$m+\bell$}}}
\put(9675,0){\makebox(0,0)[lb]{{\SetFigFont{8}{16.8}{\rmdefault}{\mddefault}{\updefault}$m+\ell$}}}
\end{picture}
}
\end{scriptsize}

\vspace*{4mm}
\begin{center}{\ft Figure 4.1: The dressed vertex model}\end{center}

\nin As well as the fixed interior edge variables
$\vep_1,\vep_2,\cdots\vep_N$ ($N=3$ is shown in the Figure), we 
fix the boundary edge variables, also marked by bullets, to take the values
shown in Figure 4.1, where $m\geq 1+\frac{k-1}{2}$ and $\ell\in \{0,1,\cdots,k\}$. The total Boltzmann weight associated
with a configuration of	edge variables is given by the product
of $R$ matrix values around 4-vertices and $\psi$ and $\psi^*$ values 
around 3-vertices. 
The rational for fixing the $\cdots
,m+\bell,m+\ell,m+\bell,m+\ell,m+\bell , \cdots$ boundary condition is that for suitably large 
lattices it imposes the vertex model boundary condition corresponding to 
the ground state shown in
Figure 2.2. This follows from the observations concerning Figure 2.9.

Let us denote the weighted sum over all edge variably
configurations of Figure 4.1, with $\ep_1,\cdots,\ep_N$ fixed, by $F^{(\ell)}_{L;m}(\vep_1,\vep_2,\cdots,\vep_N)$. The correlation
function we are interested in is the  ratios
\ben P^{(\ell)}_{L;m}(\vep_1,\vep_2,\cdots,\vep_N)\equiv \frac{1}{Z^{(\ell)}_{L;m}}
    F^{(\ell)}_{L;m}(\vep_1,\vep_2,\cdots,\vep_N)\een
where the partition function $Z^{(\ell)}_{L;m}$ is the corresponding unrestricted sum,
 i.e.,
\ben 
Z^{(\ell)}_{L;m}=\sum\limits_{\vep_1,\vep_2,\cdots,\vep_N}  F^{(\ell)}_{L;m}(\vep_1,\vep_2,\cdots,\vep_N).
\een

Ultimately, we will consider the infinite $L$ limit of
$P^{(\ell)}_{L;m}(\vep_1,\vep_2,\cdots,\vep_N)$.
The conjecture is that in this limit the $m$ dependence associated with the boundary will disappear, and that we can identify the vertex model correlation function as 
\bea
P^{(\ell)}(\vep_1,\vep_2,\cdots,\vep_N)=\lim\limits_{L\ra \infty}
P^{(\ell)}_{L;m}(\vep_1,\vep_2,\cdots,\vep_N)\lb{ci}.\ee

\ssect{The Relationship with Fusion SOS Correlation Functions}

We will next show how to relate the dressed function   
$F^{(\ell)}_{L;m}(\vep_1,\vep_2,\cdots,\vep_N)$ corresponding to Figure 4.1 
to one associated with SOS models. The argument precedes via a number
of diagrammatic equivalences.
The first step is to successively use the vertex-face correspondence relations 
depicted in Figures 2.6 $(a)$ and $(b)$ to 
turn vertex weights into SOS weights. We use 2.6 $(a)$ 
starting from the NE corner and 2.6 $(b)$ starting the SW corner of Figure 4.1. 
A little thought and diagram drawing
will convince the reader that you can carrying this procedure through until
you end up with a dislocation that extends  from the NW to the SE
corner of
the diagram passing through the $N$ fixed edges. There
are many ways to draw this dislocation - one is shown in Figure 4.2.

\setlength{\unitlength}{0.0004in}
\begingroup\makeatletter\ifx\SetFigFont\undefined%
\gdef\SetFigFont#1#2#3#4#5{%
  \reset@font\fontsize{#1}{#2pt}%
  \fontfamily{#3}\fontseries{#4}\fontshape{#5}%
  \selectfont}%
\fi\endgroup%
{\renewcommand{\dashlinestretch}{30}
\begin{picture}(11607,8070)(-2000,600)
\thinlines
\texture{55888888 88555555 5522a222 a2555555 55888888 88555555 552a2a2a 2a555555 
        55888888 88555555 55a222a2 22555555 55888888 88555555 552a2a2a 2a555555 
        55888888 88555555 5522a222 a2555555 55888888 88555555 552a2a2a 2a555555 
        55888888 88555555 55a222a2 22555555 55888888 88555555 552a2a2a 2a555555 }
\put(5850,3975){\shade\ellipse{100}{100}}
\put(5850,3975){\ellipse{100}{100}}
\put(4950,3975){\shade\ellipse{100}{100}}
\put(4950,3975){\ellipse{100}{100}}
\put(1800,7575){\shade\ellipse{100}{100}}
\put(1800,7575){\ellipse{100}{100}}
\put(2700,7575){\shade\ellipse{100}{100}}
\put(2700,7575){\ellipse{100}{100}}
\put(3600,7575){\shade\ellipse{100}{100}}
\put(3600,7575){\ellipse{100}{100}}
\put(5400,7575){\shade\ellipse{100}{100}}
\put(5400,7575){\ellipse{100}{100}}
\put(6300,7575){\shade\ellipse{100}{100}}
\put(6300,7575){\ellipse{100}{100}}
\put(7200,7575){\shade\ellipse{100}{100}}
\put(7200,7575){\ellipse{100}{100}}
\put(8100,7575){\shade\ellipse{100}{100}}
\put(8100,7575){\ellipse{100}{100}}
\put(6750,3975){\shade\ellipse{100}{100}}
\put(6750,3975){\ellipse{100}{100}}
\put(9000,7575){\shade\ellipse{100}{100}}
\put(9000,7575){\ellipse{100}{100}}
\put(9900,7575){\shade\ellipse{100}{100}}
\put(9900,7575){\ellipse{100}{100}}
\put(9900,375){\shade\ellipse{100}{100}}
\put(9900,375){\ellipse{100}{100}}
\put(1800,375){\shade\ellipse{100}{100}}
\put(1800,375){\ellipse{100}{100}}
\put(2700,375){\shade\ellipse{100}{100}}
\put(2700,375){\ellipse{100}{100}}
\put(3600,375){\shade\ellipse{100}{100}}
\put(3600,375){\ellipse{100}{100}}
\put(5400,375){\shade\ellipse{100}{100}}
\put(5400,375){\ellipse{100}{100}}
\put(6300,375){\shade\ellipse{100}{100}}
\put(6300,375){\ellipse{100}{100}}
\put(7200,375){\shade\ellipse{100}{100}}
\put(7200,375){\ellipse{100}{100}}
\put(8100,375){\shade\ellipse{100}{100}}
\put(8100,375){\ellipse{100}{100}}
\put(9000,375){\shade\ellipse{100}{100}}
\put(9000,375){\ellipse{100}{100}}
\put(900,6675){\shade\ellipse{100}{100}}
\put(900,6675){\ellipse{100}{100}}
\put(900,5775){\shade\ellipse{100}{100}}
\put(900,5775){\ellipse{100}{100}}
\put(900,4875){\shade\ellipse{100}{100}}
\put(900,4875){\ellipse{100}{100}}
\put(900,3075){\shade\ellipse{100}{100}}
\put(900,3075){\ellipse{100}{100}}
\put(900,2175){\shade\ellipse{100}{100}}
\put(900,2175){\ellipse{100}{100}}
\put(900,1275){\shade\ellipse{100}{100}}
\put(900,1275){\ellipse{100}{100}}
\put(900,375){\shade\ellipse{100}{100}}
\put(900,375){\ellipse{100}{100}}
\put(10800,7575){\shade\ellipse{100}{100}}
\put(10800,7575){\ellipse{100}{100}}
\put(10800,6675){\shade\ellipse{100}{100}}
\put(10800,6675){\ellipse{100}{100}}
\put(10800,5775){\shade\ellipse{100}{100}}
\put(10800,5775){\ellipse{100}{100}}
\put(10800,4875){\shade\ellipse{100}{100}}
\put(10800,4875){\ellipse{100}{100}}
\put(10800,3075){\shade\ellipse{100}{100}}
\put(10800,3075){\ellipse{100}{100}}
\put(10800,2175){\shade\ellipse{100}{100}}
\put(10800,2175){\ellipse{100}{100}}
\put(10800,1275){\shade\ellipse{100}{100}}
\put(10800,1275){\ellipse{100}{100}}
\put(885,7365){\shade\ellipse{100}{100}}
\put(885,7365){\ellipse{100}{100}}
\put(4650,7575){\shade\ellipse{100}{100}}
\put(4650,7575){\ellipse{100}{100}}
\put(900,3825){\shade\ellipse{100}{100}}
\put(900,3825){\ellipse{100}{100}}
\put(4350,375){\shade\ellipse{100}{100}}
\put(4350,375){\ellipse{100}{100}}
\put(10800,4125){\shade\ellipse{100}{100}}
\put(10800,4125){\ellipse{100}{100}}
\put(10575,375){\shade\ellipse{100}{100}}
\put(10575,375){\ellipse{100}{100}}
\path(1800,7575)(9000,7575)
\path(9900,7575)(9000,7575)
\path(9900,7575)(9000,7575)
\path(1800,375)(9900,375)
\path(1800,375)(9900,375)
\path(9900,7575)(10800,7575)
\path(900,375)(1800,375)
%\path(885,7365)(1800,7575)
%%%
\path(885,7575)(1800,7575)
%%%
\path(1800,7365)(1800,420)
\path(885,7365)(885,360)
\path(2265,7575)(2265,7365)
\path(2700,7350)(2700,405)
\path(3165,7560)(3165,7350)
\path(4080,7560)(4080,7365)
\path(885,3810)(10575,3810)
\path(885,3810)(10575,3810)
\path(10590,3810)(10590,375)
\path(10590,3810)(10590,375)
\path(885,5775)(4290,5775)
\path(900,4875)(4275,4875)
\path(4650,6675)(10800,6675)
\path(4650,5775)(10800,5775)
\path(4650,4875)(10800,4875)
\path(4650,7575)(4650,4200)
\path(4350,7350)(4350,3825)
\path(4350,7350)(4350,3825)
\path(900,7350)(4350,7350)
\path(900,7350)(4350,7350)
\path(3600,7350)(3600,375)
\path(3600,7350)(3600,375)
\path(4350,3825)(4350,375)
\path(4350,3825)(4350,375)
\path(4350,4875)(4200,4875)
\path(4350,4875)(4200,4875)
\path(4350,5775)(4125,5775)
\path(4350,5775)(4125,5775)
\path(900,6675)(4350,6675)
\path(900,6675)(4350,6675)
\path(4650,4125)(10800,4125)
\path(4650,4125)(10800,4125)
\path(4650,4275)(4650,4125)
\path(4650,4275)(4650,4125)
\path(5400,7575)(5400,4125)
\path(5400,7575)(5400,4125)
\path(6300,7575)(6300,4125)
\path(6300,7575)(6300,4125)
\path(7200,7575)(7200,4125)
\path(7200,7575)(7200,4125)
\path(8100,7575)(8100,4125)
\path(8100,7575)(8100,4125)
\path(9000,7575)(9000,4125)
\path(9000,7575)(9000,4125)
\path(9900,7575)(9900,4125)
\path(9900,7575)(9900,4125)
\path(900,3075)(10575,3075)
\path(900,3075)(10575,3075)
\path(900,2175)(10575,2175)
\path(900,2175)(10575,2175)
\path(900,1275)(10575,1275)
\path(900,1275)(10575,1275)
%\path(1425,7500)(1425,7350)
%\path(1425,7500)(1425,7350)
%\path(10800,1275)(10575,375)
%\path(10800,1275)(10575,375)
\path(10800,7575)(10800,1275)
\path(10800,7575)(10800,1275)
\path(9900,375)(10575,375)
\path(9900,375)(10575,375)
\path(4350,7050)(4650,7050)
\path(4350,7050)(4650,7050)
\path(4350,6225)(4650,6225)
\path(4350,6225)(4650,6225)
\path(4350,5325)(4650,5325)
\path(4350,5325)(4650,5325)
\path(4350,4425)(4650,4425)
\path(4350,4425)(4650,4425)
\path(4950,4125)(4950,3825)
\path(4950,4125)(4950,3825)
\path(5850,4125)(5850,3825)
\path(5850,4125)(5850,3825)
\path(6750,4125)(6750,3825)
\path(6750,4125)(6750,3825)
\path(7650,4125)(7650,3825)
\path(7650,4125)(7650,3825)
\path(8550,4125)(8550,3825)
\path(8550,4125)(8550,3825)
\path(9450,4125)(9450,3825)
\path(9450,4125)(9450,3825)
\path(10350,4125)(10350,3825)
\path(10350,4125)(10350,3825)
\path(10575,3450)(10800,3450)
\path(10575,3450)(10800,3450)
\path(10575,2625)(10800,2625)
\path(10575,2625)(10800,2625)
\path(10575,1725)(10800,1725)
\path(10575,1725)(10800,1725)
%\path(10575,975)(10725,975)
%\path(10575,975)(10725,975)
\path(5400,3825)(5400,375)
\path(5400,3825)(5400,375)
\path(6300,3825)(6300,375)
\path(6300,3825)(6300,375)
\path(7200,3825)(7200,375)
\path(7200,3825)(7200,375)
\path(8100,3825)(8100,375)
\path(8100,3825)(8100,375)
\path(9000,3825)(9000,375)
\path(9000,3825)(9000,375)
\path(9900,3825)(9900,375)
\path(9900,3825)(9900,375)
\path(900,6450)(1125,6675)
\path(900,6450)(1125,6675)
\path(2025,6675)(1800,6450)
\path(2025,6675)(1800,6450)
\path(2925,6675)(2700,6450)
\path(2925,6675)(2700,6450)
\path(2025,7350)(1800,7125)
\path(2025,7350)(1800,7125)
\path(1125,7350)(900,7125)
\path(1125,7350)(900,7125)
\path(2925,7350)(2700,7125)
\path(2925,7350)(2700,7125)
\path(3825,7350)(3600,7125)
\path(3825,7350)(3600,7125)
\path(3825,6675)(3600,6450)
\path(3825,6675)(3600,6450)
\path(1125,5775)(900,5550)
\path(1125,5775)(900,5550)
\path(2025,5775)(1800,5550)
\path(2025,5775)(1800,5550)
\path(2925,5775)(2700,5550)
\path(2925,5775)(2700,5550)
\path(3825,5775)(3600,5550)
\path(3825,5775)(3600,5550)
\path(1125,4875)(900,4650)
\path(1125,4875)(900,4650)
\path(2025,4875)(1800,4650)
\path(2025,4875)(1800,4650)
\path(2925,4875)(2700,4650)
\path(2925,4875)(2700,4650)
\path(3825,4875)(3600,4650)
\path(3825,4875)(3600,4650)
\path(5625,7575)(5400,7350)
\path(5625,7575)(5400,7350)
\path(4875,7575)(4650,7350)
\path(4875,7575)(4650,7350)
\path(6525,7575)(6300,7350)
\path(6525,7575)(6300,7350)
\path(7425,7575)(7200,7350)
\path(7425,7575)(7200,7350)
\path(8325,7575)(8100,7350)
\path(8325,7575)(8100,7350)
\path(9225,7575)(9000,7350)
\path(9225,7575)(9000,7350)
\path(10125,7575)(9900,7350)
\path(10125,7575)(9900,7350)
\path(4875,6675)(4650,6450)
\path(4875,6675)(4650,6450)
\path(5625,6675)(5400,6450)
\path(5625,6675)(5400,6450)
\path(6525,6675)(6300,6450)
\path(6525,6675)(6300,6450)
\path(7425,6675)(7200,6450)
\path(7425,6675)(7200,6450)
\path(8325,6675)(8100,6450)
\path(8325,6675)(8100,6450)
\path(9225,6675)(9000,6450)
\path(9225,6675)(9000,6450)
\path(10125,6675)(9900,6450)
\path(10125,6675)(9900,6450)
\path(4875,5775)(4650,5550)
\path(4875,5775)(4650,5550)
\path(5625,5775)(5400,5550)
\path(5625,5775)(5400,5550)
\path(6525,5775)(6300,5550)
\path(6525,5775)(6300,5550)
\path(7425,5775)(7200,5550)
\path(7425,5775)(7200,5550)
\path(8325,5775)(8100,5550)
\path(8325,5775)(8100,5550)
\path(9225,5775)(9000,5550)
\path(9225,5775)(9000,5550)
\path(10125,5775)(9900,5550)
\path(10125,5775)(9900,5550)
\path(900,3600)(1125,3825)
\path(900,3600)(1125,3825)
\path(2025,3825)(1800,3600)
\path(2025,3825)(1800,3600)
\path(2925,3825)(2700,3600)
\path(2925,3825)(2700,3600)
\path(3825,3825)(3600,3600)
\path(3825,3825)(3600,3600)
\path(1125,3075)(900,2850)
\path(1125,3075)(900,2850)
\path(2025,3075)(1800,2850)
\path(2025,3075)(1800,2850)
\path(2925,3075)(2700,2850)
\path(2925,3075)(2700,2850)
\path(3825,3075)(3600,2850)
\path(3825,3075)(3600,2850)
\path(4575,3075)(4350,2850)
\path(4575,3075)(4350,2850)
\path(5625,3075)(5400,2850)
\path(5625,3075)(5400,2850)
\path(6525,3075)(6300,2850)
\path(6525,3075)(6300,2850)
\path(7425,3075)(7200,2850)
\path(7425,3075)(7200,2850)
\path(8325,3075)(8100,2850)
\path(8325,3075)(8100,2850)
\path(4350,3600)(4575,3825)
\path(4350,3600)(4575,3825)
\path(5400,3600)(5625,3825)
\path(5400,3600)(5625,3825)
\path(6300,3600)(6525,3825)
\path(6300,3600)(6525,3825)
\path(7200,3600)(7425,3825)
\path(7200,3600)(7425,3825)
\path(8100,3600)(8325,3825)
\path(8100,3600)(8325,3825)
\path(9000,3600)(9225,3825)
\path(9000,3600)(9225,3825)
\path(9900,3600)(10125,3825)
\path(9900,3600)(10125,3825)
\path(9900,2850)(10125,3075)
\path(9900,2850)(10125,3075)
\path(9000,2850)(9225,3075)
\path(9000,2850)(9225,3075)
\path(900,1950)(1125,2175)
\path(900,1950)(1125,2175)
\path(2025,2175)(1800,1950)
\path(2025,2175)(1800,1950)
\path(2925,2175)(2700,1950)
\path(2925,2175)(2700,1950)
\path(3825,2175)(3600,1950)
\path(3825,2175)(3600,1950)
\path(4575,2175)(4350,1950)
\path(4575,2175)(4350,1950)
\path(5625,2175)(5400,1950)
\path(5625,2175)(5400,1950)
\path(6525,2175)(6300,1950)
\path(6525,2175)(6300,1950)
\path(7425,2175)(7200,1950)
\path(7425,2175)(7200,1950)
\path(8325,2175)(8100,1950)
\path(8325,2175)(8100,1950)
\path(9225,2175)(9000,1950)
\path(9225,2175)(9000,1950)
\path(10125,2175)(9900,1950)
\path(10125,2175)(9900,1950)
\path(4875,4875)(4650,4650)
\path(4875,4875)(4650,4650)
\path(5625,4875)(5400,4650)
\path(5625,4875)(5400,4650)
\path(6525,4875)(6300,4650)
\path(6525,4875)(6300,4650)
\path(7425,4875)(7200,4650)
\path(7425,4875)(7200,4650)
\path(8325,4875)(8100,4650)
\path(8325,4875)(8100,4650)
\path(9225,4875)(9000,4650)
\path(9225,4875)(9000,4650)
\path(10125,4875)(9900,4650)
\path(10125,4875)(9900,4650)
\path(1125,1275)(900,1050)
\path(1125,1275)(900,1050)
\path(2025,1275)(1800,1050)
\path(2025,1275)(1800,1050)
\path(2925,1275)(2700,1050)
\path(2925,1275)(2700,1050)
\path(3825,1275)(3600,1050)
\path(3825,1275)(3600,1050)
\path(4575,1275)(4350,1050)
\path(4575,1275)(4350,1050)
\path(5625,1275)(5400,1050)
\path(5625,1275)(5400,1050)
\path(6525,1275)(6300,1050)
\path(6525,1275)(6300,1050)
\path(7425,1275)(7200,1050)
\path(7425,1275)(7200,1050)
\path(8325,1275)(8100,1050)
\path(8325,1275)(8100,1050)
\path(9225,1275)(9000,1050)
\path(9225,1275)(9000,1050)
\path(10125,1275)(9900,1050)
\path(10125,1275)(9900,1050)
\put(5175,7875){\makebox(0,0)[lb]{{\SetFigFont{8}{16.8}{\rmdefault}{\mddefault}{\updefault}$m+\bell$}}}
\put(6100,7875){\makebox(0,0)[lb]{{\SetFigFont{8}{16.8}{\rmdefault}{\mddefault}{\updefault}$m+\ell$}}}
\put(6975,7875){\makebox(0,0)[lb]{{\SetFigFont{8}{16.8}{\rmdefault}{\mddefault}{\updefault}$m+\bell$}}}
\put(7950,7875){\makebox(0,0)[lb]{{\SetFigFont{8}{16.8}{\rmdefault}{\mddefault}{\updefault}$m+\ell$}}}
\put(8850,7875){\makebox(0,0)[lb]{{\SetFigFont{8}{16.8}{\rmdefault}{\mddefault}{\updefault}$m+\bell$}}}
\put(9750,7875){\makebox(0,0)[lb]{{\SetFigFont{8}{16.8}{\rmdefault}{\mddefault}{\updefault}$m+\ell$}}}
\put(3375,7875){\makebox(0,0)[lb]{{\SetFigFont{8}{16.8}{\rmdefault}{\mddefault}{\updefault}$m+\bell$}}}
\put(2550,7875){\makebox(0,0)[lb]{{\SetFigFont{8}{16.8}{\rmdefault}{\mddefault}{\updefault}$m+\ell$}}}
\put(1650,7875){\makebox(0,0)[lb]{{\SetFigFont{8}{16.8}{\rmdefault}{\mddefault}{\updefault}$m+\bell$}}}
\put(750,7875){\makebox(0,0)[lb]{{\SetFigFont{8}{16.8}{\rmdefault}{\mddefault}{\updefault}$m+\ell$}}}
\put(10650,7875){\makebox(0,0)[lb]{{\SetFigFont{8}{16.8}{\rmdefault}{\mddefault}{\updefault}$m+\bell$}}}
\put(11100,6675){\makebox(0,0)[lb]{{\SetFigFont{8}{16.8}{\rmdefault}{\mddefault}{\updefault}$m+\ell$}}}
\put(11100,5775){\makebox(0,0)[lb]{{\SetFigFont{8}{16.8}{\rmdefault}{\mddefault}{\updefault}$m+\bell$}}}
\put(11100,4875){\makebox(0,0)[lb]{{\SetFigFont{8}{16.8}{\rmdefault}{\mddefault}{\updefault}$m+\ell$}}}
\put(11100,3075){\makebox(0,0)[lb]{{\SetFigFont{8}{16.8}{\rmdefault}{\mddefault}{\updefault}$m+\ell$}}}
\put(11100,2175){\makebox(0,0)[lb]{{\SetFigFont{8}{16.8}{\rmdefault}{\mddefault}{\updefault}$m+\bell$}}}
\put(11100,1275){\makebox(0,0)[lb]{{\SetFigFont{8}{16.8}{\rmdefault}{\mddefault}{\updefault}$m+\ell$}}}
\put(11100,375){\makebox(0,0)[lb]{{\SetFigFont{8}{16.8}{\rmdefault}{\mddefault}{\updefault}$m+\bell$}}}
\put(0,6675){\makebox(0,0)[lb]{{\SetFigFont{8}{16.8}{\rmdefault}{\mddefault}{\updefault}$m+\bell$}}}
\put(0,5775){\makebox(0,0)[lb]{{\SetFigFont{8}{16.8}{\rmdefault}{\mddefault}{\updefault}$m+\ell$}}}
\put(0,4875){\makebox(0,0)[lb]{{\SetFigFont{8}{16.8}{\rmdefault}{\mddefault}{\updefault}$m+\bell$}}}
\put(0,3075){\makebox(0,0)[lb]{{\SetFigFont{8}{16.8}{\rmdefault}{\mddefault}{\updefault}$m+\bell$}}}
\put(0,2175){\makebox(0,0)[lb]{{\SetFigFont{8}{16.8}{\rmdefault}{\mddefault}{\updefault}$m+\ell$}}}
\put(0,1275){\makebox(0,0)[lb]{{\SetFigFont{8}{16.8}{\rmdefault}{\mddefault}{\updefault}$m+\bell$}}}
\put(600,0){\makebox(0,0)[lb]{{\SetFigFont{8}{16.8}{\rmdefault}{\mddefault}{\updefault}$m+\ell$}}}
\put(1575,0){\makebox(0,0)[lb]{{\SetFigFont{8}{16.8}{\rmdefault}{\mddefault}{\updefault}$m+\bell$}}}
\put(2475,0){\makebox(0,0)[lb]{{\SetFigFont{8}{16.8}{\rmdefault}{\mddefault}{\updefault}$m+\ell$}}}
\put(3375,0){\makebox(0,0)[lb]{{\SetFigFont{8}{16.8}{\rmdefault}{\mddefault}{\updefault}$m+\bell$}}}
\put(4275,0){\makebox(0,0)[lb]{{\SetFigFont{8}{16.8}{\rmdefault}{\mddefault}{\updefault}$m+\ell$}}}
\put(5175,0){\makebox(0,0)[lb]{{\SetFigFont{8}{16.8}{\rmdefault}{\mddefault}{\updefault}$m+\bell$}}}
\put(6000,0){\makebox(0,0)[lb]{{\SetFigFont{8}{16.8}{\rmdefault}{\mddefault}{\updefault}$m+\ell$}}}
\put(6975,0){\makebox(0,0)[lb]{{\SetFigFont{8}{16.8}{\rmdefault}{\mddefault}{\updefault}$m+\bell$}}}
\put(7875,0){\makebox(0,0)[lb]{{\SetFigFont{8}{16.8}{\rmdefault}{\mddefault}{\updefault}$m+\ell$}}}
\put(8775,0){\makebox(0,0)[lb]{{\SetFigFont{8}{16.8}{\rmdefault}{\mddefault}{\updefault}$m+\bell$}}}
\put(9675,0){\makebox(0,0)[lb]{{\SetFigFont{8}{16.8}{\rmdefault}{\mddefault}{\updefault}$m+\ell$}}}
\put(4575,3900){\makebox(0,0)[lb]{{\SetFigFont{8}{16.8}{\rmdefault}{\mddefault}{\updefault}$\vep_1$}}}
\put(5550,3900){\makebox(0,0)[lb]{{\SetFigFont{8}{16.8}{\rmdefault}{\mddefault}{\updefault}$\vep_2$}}}
\put(6375,3900){\makebox(0,0)[lb]{{\SetFigFont{8}{16.8}{\rmdefault}{\mddefault}{\updefault}$\vep_3$}}}
\put(4200,7875){\makebox(0,0)[lb]{{\SetFigFont{8}{16.8}{\rmdefault}{\mddefault}{\updefault}   $m+\ell$}}}
\put(0,3750){\makebox(0,0)[lb]{{\SetFigFont{8}{16.8}{\rmdefault}{\mddefault}{\updefault}$m+\ell$}}}
\put(11100,4050){\makebox(0,0)[lb]{{\SetFigFont{8}{16.8}{\rmdefault}{\mddefault}{\updefault}$m+\bell$}}}
%%% horiz ar

%\input{FIG36.eepic}rows
\path(4500,7095)(4400,7045)(4500,6995)
\path(4500,6275)(4400,6225)(4500,6175)
\path(4500,5375)(4400,5325)(4500,5275)
\path(4500,4475)(4400,4425)(4500,4375)

%\path(10700,1025)(10600,975)(10700,925)
\path(10700,1775)(10600,1725)(10700,1675)
\path(10700,2675)(10600,2625)(10700,2575)
\path(10700,3500)(10600,3450)(10700,3400)

%%% vert arrows
\path(1315,7525)(1365,7425)(1415,7525)
\path(2215,7525)(2265,7425)(2315,7525)
\path(3115,7525)(3165,7425)(3215,7525)
\path(4015,7525)(4065,7425)(4115,7525)

\path(4900,3900)(4950,3800)(5000,3900)
\path(5800,3900)(5850,3800)(5900,3900)
\path(6700,3900)(6750,3800)(6800,3900)
\path(7600,3950)(7650,3850)(7700,3950)
\path(8500,3950)(8550,3850)(8600,3950)
\path(9400,3950)(9450,3850)(9500,3950)
\path(10300,3950)(10350,3850)(10400,3950)

% additional lines
\dashline{50}(900,7350)(900,7400)(900,7450)(900,7500)(900,7575)    % dashed
\dashline{50}(10575,375)(10625,375)(10675,375)(10725,375)(10800,375)  % dashed

\path(10575,825)(10800,825)
\path(10800,375)(10800,1275)
\path(1365,7350)(1365,7575)

\put(10800,375){\shade\ellipse{100}{100}}
\put(10800,375){\ellipse{100}{100}}
\put(900,7575){\shade\ellipse{100}{100}}
\put(900,7575){\ellipse{100}{100}}

\path(10700,875)(10600,825)(10700,775)

\end{picture}
}

 \vspace*{4mm}\begin{center}{\ft Figure 4.2: The correlation function after converting to SOS weights}
 \end{center}

The next step is to use the relation of Figure 2.7 $(a)$ to remove the 
dislocation step by step starting from the NW. We can do this until
we reach the leftmost central fixed edge variable $\vep_1$.
Hence the sum in Figure 4.2 is equal to that of
Figure 4.3 after taking the infinite volume limit and dividing into the parts 
corresponding to SOS model CTMs and HTMs.

\setlength{\unitlength}{0.00035in}
\begingroup\makeatletter\ifx\SetFigFont\undefined%
\gdef\SetFigFont#1#2#3#4#5{%
  \reset@font\fontsize{#1}{#2pt}%
  \fontfamily{#3}\fontseries{#4}\fontshape{#5}%
  \selectfont}%
\fi\endgroup%
{\renewcommand{\dashlinestretch}{30}
\begin{picture}(12573,10530)(-1500,1600)
\thinlines
\path(4725,9525)(4725,5925)
\path(1125,5925)(4650,5925)
\path(2925,8625)(4725,8625)
\path(2025,7725)(4725,7725)
\path(1125,6825)(1125,5925)
\path(2925,8625)(2925,5925)
\path(2025,7725)(2025,5925)
\path(1125,4725)(4650,4725)
\path(4725,4725)(4725,1125)
\path(1125,4725)(1125,3825)
\path(3825,1125)(4725,1125)
\path(2025,2925)(4725,2925)
\path(2925,2025)(4725,2025)
\path(2925,4725)(2925,2025)
\path(2025,4725)(2025,2925)
\path(8625,9525)(8625,5925)
\path(8625,5925)(12225,5925)
\path(8625,9525)(9525,9525)
\path(12225,6825)(12225,5925)
\path(8625,8625)(10425,8625)
\path(8625,7725)(11325,7725)
\path(3825,9525)(4725,9525)
\path(2250,4725)(2025,4500)
\path(3150,4725)(2925,4500)
\path(4050,4725)(3825,4500)
\path(2250,3825)(2025,3600)
\path(3150,3825)(2925,3600)
\path(4050,3825)(3825,3600)
\path(3150,2925)(2925,2700)
\path(4050,2925)(3825,2700)
\path(4050,2025)(3825,1800)
\path(4050,9525)(3825,9300)
\path(3150,8625)(2925,8400)
\path(4050,8625)(3825,8400)
\path(2250,7725)(2025,7500)
\path(3150,7725)(2925,7500)
\path(4050,7725)(3825,7500)
\path(1350,6825)(1125,6600)
\path(2250,6825)(2025,6600)
\path(3150,6825)(2925,6600)
\path(4050,6825)(3825,6600)
\dottedline{68}(5325,5925)(4725,5925)(4725,4725)(5325,4725)
\dottedline{68}(8025,5925)(8625,5925)(8625,5625)
\dottedline{68}(8625,5025)(8625,4725)(8025,4725)
\path(10425,8625)(10425,5925)
\put(12450,4275){\makebox(0,0)[lb]{\smash{{{\SetFigFont{8}{16.8}{\rmdefault}{\mddefault}{\updefault}$\hdots$}}}}}
\put(4425,5250){\makebox(0,0)[lb]{\smash{{{\SetFigFont{8}{16.8}{\rmdefault}{\mddefault}{\updefault}$a$}}}}}
\put(6075,5625){\makebox(0,0)[lb]{\smash{{{\SetFigFont{8}{16.8}{\rmdefault}{\mddefault}{\updefault}$a_1$}}}}}
\put(6675,5625){\makebox(0,0)[lb]{\smash{{{\SetFigFont{8}{16.8}{\rmdefault}{\mddefault}{\updefault}$a_{N-1}$}}}}}
\put(7875,5625){\makebox(0,0)[lb]{\smash{{{\SetFigFont{8}{16.8}{\rmdefault}{\mddefault}{\updefault}$a_N$}}}}}
\put(6075,4875){\makebox(0,0)[lb]{\smash{{{\SetFigFont{8}{16.8}{\rmdefault}{\mddefault}{\updefault}$a'_1$}}}}}
\put(6675,4875){\makebox(0,0)[lb]{\smash{{{\SetFigFont{8}{16.8}{\rmdefault}{\mddefault}{\updefault}$a'_{N-1}$}}}}}
\put(7875,4875){\makebox(0,0)[lb]{\smash{{{\SetFigFont{8}{16.8}{\rmdefault}{\mddefault}{\updefault}$a'_N$}}}}}
\put(5375,5250){\makebox(0,0)[lb]{\smash{{{\SetFigFont{8}{16.8}{\rmdefault}{\mddefault}{\updefault}$\vep_1$}}}}}
%\put(6300,5250){\makebox(0,0)[lb]{\smash{{{\SetFigFont{8}{16.8}{\rmdefault}{\mddefault}{\updefault}$\vep_2$}}}}}
\put(7075,5250){\makebox(0,0)[lb]{\smash{{{\SetFigFont{8}{16.8}{\rmdefault}{\mddefault}{\updefault}$\vep_N$}}}}}
\put(-800,5400){\makebox(0,0)[lb]{\smash{{{\SetFigFont{16}{16.8}{\rmdefault}{\mddefault}{\updefault}$\sum\limits_{a, a_j,a'_j}$}}}}}
\put(13400,5200){\makebox(0,0)[lb]{\smash{{{\SetFigFont{12}{16.8}{\rmdefault}{\mddefault}{\updefault}$\tgL_{a_N,a'_N}(u_0)$}}}}}
\put(1050,8700){\makebox(0,0)[lb]{\smash{{{\SetFigFont{12}{16.8}{\rmdefault}{\mddefault}{\updefault}$A_a(u)$}}}}}
\put(11250,8700){\makebox(0,0)[lb]{\smash{{{\SetFigFont{12}{16.8}{\rmdefault}{\mddefault}{\updefault}$A_{NE;\,a_N}(u)$}}}}}
\put(1050,1900){\makebox(0,0)[lb]{\smash{{{\SetFigFont{12}{16.8}{\rmdefault}{\mddefault}{\updefault}$A_{SW;\,a}(u)$}}}}}
\put(11550,1900){\makebox(0,0)[lb]{\smash{{{\SetFigFont{12}{16.8}{\rmdefault}{\mddefault}{\updefault}$A_{SE;\,a'_N}(u)$}}}}}
\put(5025,10850){\makebox(0,0)[lb]{\smash{{{\SetFigFont{8}{16.8}{\rmdefault}{\mddefault}{\updefault}$\Phi_{a,\,a_1}(u)$}}}}}
\put(6425,10850){\makebox(0,0)[lb]{\smash{{{\SetFigFont{8}{16.8}{\rmdefault}{\mddefault}{\updefault}$\hdots$}}}}}
\put(7025,10850){\makebox(0,0)[lb]{\smash{{{\SetFigFont{8}{16.8}{\rmdefault}{\mddefault}{\updefault}$\Phi_{a_{N-1},\,a_N}(u)$}}}}}
\put(4825,-200){\makebox(0,0)[lb]{\smash{{{\SetFigFont{8}{16.8}{\rmdefault}{\mddefault}{\updefault}$\Phi_{S;a'_1,a}(u)$}}}}}
\put(6425,9000){\makebox(0,0)[lb]{\smash{{{\SetFigFont{8}{16.8}{\rmdefault}{\mddefault}{\updefault}$\hdots$}}}}}
\put(6425,8100){\makebox(0,0)[lb]{\smash{{{\SetFigFont{8}{16.8}{\rmdefault}{\mddefault}{\updefault}$\hdots$}}}}}
\put(6425,7200){\makebox(0,0)[lb]{\smash{{{\SetFigFont{8}{16.8}{\rmdefault}{\mddefault}{\updefault}$\hdots$}}}}}
\put(6425,6300){\makebox(0,0)[lb]{\smash{{{\SetFigFont{8}{16.8}{\rmdefault}{\mddefault}{\updefault}$\hdots$}}}}}
\put(6425,4300){\makebox(0,0)[lb]{\smash{{{\SetFigFont{8}{16.8}{\rmdefault}{\mddefault}{\updefault}$\hdots$}}}}}
\put(6425,3400){\makebox(0,0)[lb]{\smash{{{\SetFigFont{8}{16.8}{\rmdefault}{\mddefault}{\updefault}$\hdots$}}}}}
\put(6425,2500){\makebox(0,0)[lb]{\smash{{{\SetFigFont{8}{16.8}{\rmdefault}{\mddefault}{\updefault}$\hdots$}}}}}
\put(6425,1600){\makebox(0,0)[lb]{\smash{{{\SetFigFont{8}{16.8}{\rmdefault}{\mddefault}{\updefault}$\hdots$}}}}}
\put(6425,-200){\makebox(0,0)[lb]{\smash{{{\SetFigFont{8}{16.8}{\rmdefault}{\mddefault}{\updefault}$\hdots$}}}}}

\put(7025,-200){\makebox(0,0)[lb]{\smash{{{\SetFigFont{8}{16.8}{\rmdefault}{\mddefault}{\updefault}$\Phi_{S;\,a'_N,a'_{N-1}}(u)$}}}}}
\put(12450,5325){\makebox(0,0)[lb]{\smash{{{\SetFigFont{8}{16.8}{\rmdefault}{\mddefault}{\updefault}$\hdots$}}}}}
\path(11325,7725)(11325,5925)
\path(8700,4725)(12225,4725)
\path(8625,4725)(8625,1125)
\path(8625,1125)(9525,1125)
\path(12225,4725)(12225,3825)
\path(11325,4725)(11325,2925)
\path(8625,2025)(10425,2025)
\path(8625,2925)(11325,2925)
\path(10425,4725)(10425,2025)
\put(4125,9900){\makebox(0,0)[lb]{\smash{{{\SetFigFont{8}{16.8}{\rmdefault}{\mddefault}{\updefault}$\vdots$}}}}}
\put(5700,9900){\makebox(0,0)[lb]{\smash{{{\SetFigFont{8}{16.8}{\rmdefault}{\mddefault}{\updefault}$\vdots$}}}}}
%\put(6525,9900){\makebox(0,0)[lb]{\smash{{{\SetFigFont{8}{16.8}{\rmdefault}{\mddefault}{\updefault}$\vdots$}}}}}
\put(7500,9900){\makebox(0,0)[lb]{\smash{{{\SetFigFont{8}{16.8}{\rmdefault}{\mddefault}{\updefault}$\vdots$}}}}}
\put(8925,9900){\makebox(0,0)[lb]{\smash{{{\SetFigFont{8}{16.8}{\rmdefault}{\mddefault}{\updefault}$\vdots$}}}}}
\put(4125,650){\makebox(0,0)[lb]{\smash{{{\SetFigFont{8}{16.8}{\rmdefault}{\mddefault}{\updefault}$\vdots$}}}}}
\put(5700,650){\makebox(0,0)[lb]{\smash{{{\SetFigFont{8}{16.8}{\rmdefault}{\mddefault}{\updefault}$\vdots$}}}}}
%\put(6600,650){\makebox(0,0)[lb]{\smash{{{\SetFigFont{8}{16.8}{\rmdefault}{\mddefault}{\updefault}$\vdots$}}}}}
\put(7500,650){\makebox(0,0)[lb]{\smash{{{\SetFigFont{8}{16.8}{\rmdefault}{\mddefault}{\updefault}$\vdots$}}}}}
\put(8925,650){\makebox(0,0)[lb]{\smash{{{\SetFigFont{8}{16.8}{\rmdefault}{\mddefault}{\updefault}$\vdots$}}}}}
\put(525,6450){\makebox(0,0)[lb]{\smash{{{\SetFigFont{8}{16.8}{\rmdefault}{\mddefault}{\updefault}$\hdots$}}}}}
\put(525,4275){\makebox(0,0)[lb]{\smash{{{\SetFigFont{8}{16.8}{\rmdefault}{\mddefault}{\updefault}$\hdots$}}}}}
\put(12450,6375){\makebox(0,0)[lb]{\smash{{{\SetFigFont{8}{16.8}{\rmdefault}{\mddefault}{\updefault}$\hdots$}}}}}
\path(1350,4725)(1125,4500)
\path(7125,9525)(7125,5925)
\path(9525,9525)(9525,5925)
\path(1125,6825)(4725,6825)
%\path(5325,8625)(8025,8625)
%\path(5325,7725)(8025,7725)
%\path(5325,6825)(8025,6825)
\path(8625,6825)(12225,6825)
\path(1125,3825)(4725,3825)
%\path(5325,3825)(8025,3825)
%\path(5325,2925)(8025,2925)
%\path(5325,2025)(8025,2025)
%%%%% new %%%%
\path(5325,4725)(6225,4725)
\path(5325,3825)(6225,3825)
\path(5325,2925)(6225,2925)
\path(5325,2025)(6225,2025)
\path(5325,1125)(6225,1125)

\path(7125,4725)(8025,4725)
\path(7125,3825)(8025,3825)
\path(7125,2925)(8025,2925)
\path(7125,2025)(8025,2025)
\path(7125,1125)(8025,1125)

\path(5325,9525)(6225,9525)
\path(5325,8625)(6225,8625)
\path(5325,7725)(6225,7725)
\path(5325,6825)(6225,6825)
\path(5325,5925)(6225,5925)

\path(7125,9525)(8025,9525)
\path(7125,8625)(8025,8625)
\path(7125,7725)(8025,7725)
\path(7125,6825)(8025,6825)
\path(7125,5925)(8025,5925)

\path(5325,5925)(5325,9525)
\path(5325,1125)(5325,4725)

\path(8025,5925)(8025,9525)
\path(8025,1125)(8025,4725)

%%%%%%%%%%%%%%
\path(8625,3825)(12225,3825)
\path(3825,4725)(3825,1125)
\path(6225,4725)(6225,1125)
\path(7125,4725)(7125,1125)
\path(5775,5925)(5775,4725)
%\path(6675,5925)(6675,4725)
\path(7575,5925)(7575,4725)
\path(9075,5625)(9075,5025)
\path(9975,5625)(9975,5025)
\path(10875,5625)(10875,5025)
\path(6225,9525)(6225,5925)
\texture{55888888 88555555 5522a222 a2555555 55888888 88555555 552a2a2a 2a555555 
        55888888 88555555 55a222a2 22555555 55888888 88555555 552a2a2a 2a555555 
        55888888 88555555 5522a222 a2555555 55888888 88555555 552a2a2a 2a555555 
        55888888 88555555 55a222a2 22555555 55888888 88555555 552a2a2a 2a555555 }
\put(5775,5325){\shade\ellipse{100}{100}}
\put(5775,5325){\ellipse{100}{100}}
%\put(6675,5325){\shade\ellipse{100}{100}}
%\put(6675,5325){\ellipse{100}{100}}
\put(7575,5325){\shade\ellipse{100}{100}}
\put(7575,5325){\ellipse{100}{100}}
\put(5325,5925){\shade\ellipse{100}{100}}
\put(5325,5925){\ellipse{100}{100}}
\put(4725,5925){\shade\ellipse{100}{100}}
\put(4725,5925){\ellipse{100}{100}}
\put(4725,4725){\shade\ellipse{100}{100}}
\put(4725,4725){\ellipse{100}{100}}
\put(5325,4725){\shade\ellipse{100}{100}}
\put(5325,4725){\ellipse{100}{100}}
\put(8025,5925){\shade\ellipse{100}{100}}
\put(8025,5925){\ellipse{100}{100}}
\put(7125,5925){\shade\ellipse{100}{100}}
\put(7125,5925){\ellipse{100}{100}}
\put(6225,5925){\shade\ellipse{100}{100}}
\put(6225,5925){\ellipse{100}{100}}
\put(6225,4725){\shade\ellipse{100}{100}}
\put(6225,4725){\ellipse{100}{100}}
\put(7125,4725){\shade\ellipse{100}{100}}
\put(7125,4725){\ellipse{100}{100}}
\put(8025,4725){\shade\ellipse{100}{100}}
\put(8025,4725){\ellipse{100}{100}}
\put(8625,4725){\shade\ellipse{100}{100}}
\put(8625,4725){\ellipse{100}{100}}
\put(8625,5025){\shade\ellipse{100}{100}}
\put(8625,5025){\ellipse{100}{100}}
\put(8625,5625){\shade\ellipse{100}{100}}
\put(8625,5625){\ellipse{100}{100}}
\put(8625,5925){\shade\ellipse{100}{100}}
\put(8625,5925){\ellipse{100}{100}}
%\path(8025,9525)(5325,9525)(5325,5925)
%        (8025,5925)(8025,9525)
%\path(5325,4725)(8025,4725)(8025,1125)
%        (5325,1125)(5325,4725)
\path(8625,5625)(12225,5625)
\path(8625,5025)(12225,5025)
\path(3825,9525)(3825,5925)
\path(11775,5625)(11775,5025)
\path(8850,4725)(8625,4500)
\path(9750,4725)(9525,4500)
\path(10650,4725)(10425,4500)
\path(11550,4725)(11325,4500)
\path(8850,3825)(8625,3600)
\path(9750,3825)(9525,3600)
\path(10650,3825)(10425,3600)
\path(8850,2925)(8625,2700)
\path(9750,2925)(9525,2700)
\path(8850,2025)(8625,1800)
\path(5550,4725)(5325,4500)
%\path(6450,4725)(6225,4500)
\path(7350,4725)(7125,4500)
\path(5550,3825)(5325,3600)
%\path(6450,3825)(6225,3600)
\path(7350,3825)(7125,3600)
\path(5550,2925)(5325,2700)
%\path(6450,2925)(6225,2700)
\path(7350,2925)(7125,2700)
\path(5550,2025)(5325,1800)
%\path(6450,2025)(6225,1800)
\path(7350,2025)(7125,1800)
\path(11550,6825)(11325,6600)
\path(9525,4725)(9525,1125)
\path(9525,4725)(9525,1125)
\path(5550,9525)(5325,9300)
\path(5550,9525)(5325,9300)
%\path(6450,9525)(6225,9300)
%\path(6450,9525)(6225,9300)
\path(7350,9525)(7125,9300)
\path(7350,9525)(7125,9300)
\path(5550,8625)(5325,8400)
%\path(6450,8625)(6225,8400)
\path(7350,8625)(7125,8400)
\path(5550,7725)(5325,7500)
%\path(6450,7725)(6225,7500)
\path(7350,7725)(7125,7500)
\path(5550,6825)(5325,6600)
%\path(6450,6825)(6225,6600)
\path(7350,6825)(7125,6600)
\path(8850,9525)(8625,9300)
\path(8850,8625)(8625,8400)
\path(9750,8625)(9525,8400)
\path(8850,7725)(8625,7500)
\path(9750,7725)(9525,7500)
\path(10650,7725)(10425,7500)
\path(8850,6825)(8625,6600)
\path(9750,6825)(9525,6600)
\path(10650,6825)(10425,6600)

% vert arrows
\path(5725,5025)(5775,4875)(5825,5025)
%\path(6625,5025)(6675,4875)(6725,5025)
\path(7525,5025)(7575,4875)(7625,5025)

\path(9025,5425)(9075,5275)(9125,5425)
\path(9925,5425)(9975,5275)(10025,5425)
\path(10825,5425)(10875,5275)(10925,5425)
\path(11725,5425)(11775,5275)(11825,5425)
\end{picture}
}

\vspace*{14mm}\begin{center}{\ft Figure 4.3: The correlation function $F_m^{(\ell)}(\vep_1,\vep_2,\cdots,\vep_N)$}
\end{center}

\noi
This picture is similar to Figure 3.7 expressing the correlation function 
of fusion SOS model. But we here have two new ingredients. One is  
the L-matrix defined 
 in \eqref{Lmate}. 
The other is  the `tail operator'
$\tgL_{b,a}(u_0):\cH^{(\ell)}_{m,a}\to \cH^{(\ell)}_{m,b}$ graphically defined 
by

\setlength{\unitlength}{0.0004in}
\begingroup\makeatletter\ifx\SetFigFont\undefined%
\gdef\SetFigFont#1#2#3#4#5{%
  \reset@font\fontsize{#1}{#2pt}%
  \fontfamily{#3}\fontseries{#4}\fontshape{#5}%
  \selectfont}%
\fi\endgroup%
{\renewcommand{\dashlinestretch}{30}
\begin{picture}(4998,870)(-3000,600)
\thinlines
\texture{55888888 88555555 5522a222 a2555555 55888888 88555555 552a2a2a 2a555555 
        55888888 88555555 55a222a2 22555555 55888888 88555555 552a2a2a 2a555555 
        55888888 88555555 5522a222 a2555555 55888888 88555555 552a2a2a 2a555555 
        55888888 88555555 55a222a2 22555555 55888888 88555555 552a2a2a 2a555555 }
\put(675,690){\shade\ellipse{100}{100}}
\put(675,690){\ellipse{100}{100}}
\put(675,90){\shade\ellipse{100}{100}}
\put(675,90){\ellipse{100}{100}}
\path(675,690)(4275,690)
\path(675,690)(4275,690)
\path(675,90)(4275,90)
\path(2925,690)(2925,90)
\path(2025,690)(2025,90)
\path(1125,690)(1125,90)
\path(3825,690)(3825,90)
\put(4875,240){\makebox(0,0)[lb]{{\SetFigFont{12}{16.8}{\rmdefault}{\mddefault}{\updefault}$\hdots$}}}
\put(0,690){\makebox(0,0)[lb]{{\SetFigFont{12}{16.8}{\rmdefault}{\mddefault}{\updefault}$b$}}}
\put(0,15){\makebox(0,0)[lb]{{\SetFigFont{12}{16.8}{\rmdefault}{\mddefault}{\updefault}$a$}}}
\put(1275,350){\makebox(0,0)[lb]{{\SetFigFont{8}{16.8}{\rmdefault}{\mddefault}{\updefault}$u_0$}}}
\put(2125,350){\makebox(0,0)[lb]{{\SetFigFont{8}{16.8}{\rmdefault}{\mddefault}{\updefault}$u_0$}}}
\put(3050,350){\makebox(0,0)[lb]{{\SetFigFont{8}{16.8}{\rmdefault}{\mddefault}{\updefault}$u_0$}}}
\put(3925,350){\makebox(0,0)[lb]{{\SetFigFont{8}{16.8}{\rmdefault}{\mddefault}{\updefault}$u_0$}}}
% vertical arrows
\path(1075,450)(1125,300)(1175,450)
\path(1975,450)(2025,300)(2075,450)
\path(2875,450)(2925,300)(2975,450)
\path(3775,450)(3825,300)(3875,450)
\end{picture}
}

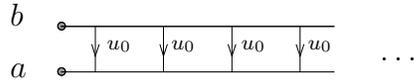
\vspace*{4mm}\begin{center}{\ft Figure 4.4: The graphical definition of the tail operator $\tgL_{b,a}(u_0)$}
\end{center}

\nin where the $\vep$'s on the vertical lines are summed over.

\ssect{The Tail Operator}
The tail operator is characterised by the following commutation relation 
\ben  
\sli_{d} \lmatsk{c}{d}{a}{b}{v}
\Phi_{E;c,d}(u)\tgL_{d,b}(u_0)=
\tgL_{c,a}(u_0) \Phi_{E;a,b}(u),
\een
where $u=u_0-v$ and the $L^{(k)}$-matrix is given in \eqref{Lmat}.
% \ben 
% \lmatsk{a}{b}{c}{d}{u}=\sli_{\vep} \lmatk{a}{b}{c}{d}{u}{\vep}.
% \een
This is due to the simple graphical argument 
given by Figure 4.5. The two steps
make successive use of the fundamental intertwining property shown in
Figure 2.6 (b).

\setlength{\unitlength}{0.00035in}
\begingroup\makeatletter\ifx\SetFigFont\undefined%
\gdef\SetFigFont#1#2#3#4#5{%
  \reset@font\fontsize{#1}{#2pt}%
  \fontfamily{#3}\fontseries{#4}\fontshape{#5}%
  \selectfont}%
\fi\endgroup%
{\renewcommand{\dashlinestretch}{30}
\begin{picture}(5973,7275)(-2000,600)
\thinlines
\path(1950,4290)(5550,4290)
\path(1350,3690)(1350,2790)(5550,2790)
\path(2400,4290)(2400,2790)
\path(3300,4290)(3300,2790)
\path(4200,4290)(4200,2790)
\path(5100,4290)(5100,2790)
\path(1350,3390)(5550,3390)
\texture{55888888 88555555 5522a222 a2555555 55888888 88555555 552a2a2a 2a555555 
        55888888 88555555 55a222a2 22555555 55888888 88555555 552a2a2a 2a555555 
        55888888 88555555 5522a222 a2555555 55888888 88555555 552a2a2a 2a555555 
        55888888 88555555 55a222a2 22555555 55888888 88555555 552a2a2a 2a555555 }
\put(1950,6090){\shade\ellipse{100}{100}}
\put(1950,6090){\ellipse{100}{100}}
\put(1950,6990){\shade\ellipse{100}{100}}
\put(1950,6990){\ellipse{100}{100}}
\put(1350,6690){\shade\ellipse{100}{100}}
\put(1350,6690){\ellipse{100}{100}}
\put(1350,5490){\shade\ellipse{100}{100}}
\put(1350,5490){\ellipse{100}{100}}
\put(1950,4290){\shade\ellipse{100}{100}}
\put(1950,4290){\ellipse{100}{100}}
\put(1350,3690){\shade\ellipse{100}{100}}
\put(1350,3690){\ellipse{100}{100}}
\put(1350,2790){\shade\ellipse{100}{100}}
\put(1350,2790){\ellipse{100}{100}}
\put(1950,1590){\shade\ellipse{100}{100}}
\put(1950,1590){\ellipse{100}{100}}
\put(1950,990){\shade\ellipse{100}{100}}
\put(1950,990){\ellipse{100}{100}}
\put(1950,90){\shade\ellipse{100}{100}}
\put(1950,90){\ellipse{100}{100}}
\path(1950,6990)(5550,6990)(5550,6090)
        (1950,6090)(1950,6990)
\path(5550,5490)(1350,5490)(1350,6390)
\path(1350,6690)(1350,6390)
\path(1350,6390)(1950,6390)
\path(2850,6990)(2850,6090)(2925,6090)
\path(3750,6990)(3750,6090)
\path(4650,6990)(4650,6090)
\path(2400,6090)(2400,5490)
\path(3300,6090)(3300,5490)
\path(4200,6090)(4200,5490)
\path(5100,6090)(5100,5490)
\path(2175,6990)(1950,6765)
\path(3075,6990)(2850,6765)
\path(3975,6990)(3750,6765)
\path(4875,6990)(4650,6765)
\path(1950,1590)(5550,1590)
\path(2400,1590)(2400,990)
\path(3300,1590)(3300,990)
\path(4200,1590)(4200,990)
\path(5100,1590)(5100,990)
\path(1950,990)(5550,990)(5550,90)
        (1950,90)(1950,990)
\path(2850,990)(2850,90)
\path(3750,990)(3750,90)
\path(4650,990)(4650,90)
\path(2175,990)(1950,765)
\path(3075,990)(2850,765)
\path(3975,990)(3750,765)
\path(4875,990)(4650,765)
\put(2075,6397){\makebox(0,0)[lb]{{\SetFigFont{8}{16.8}{\rmdefault}{\mddefault}{\updefault}$u$}}}
\put(2975,6397){\makebox(0,0)[lb]{{\SetFigFont{8}{16.8}{\rmdefault}{\mddefault}{\updefault}$u$}}}
\put(3875,6397){\makebox(0,0)[lb]{{\SetFigFont{8}{16.8}{\rmdefault}{\mddefault}{\updefault}$u$}}}
\put(4775,6397){\makebox(0,0)[lb]{{\SetFigFont{8}{16.8}{\rmdefault}{\mddefault}{\updefault}$u$}}}
\put(1575,6472){\makebox(0,0)[lb]{{\SetFigFont{8}{16.8}{\rmdefault}{\mddefault}{\updefault}$v$}}}
\put(2525,5647){\makebox(0,0)[lb]{{\SetFigFont{8}{16.8}{\rmdefault}{\mddefault}{\updefault}$u_0$}}}
\put(3450,5647){\makebox(0,0)[lb]{{\SetFigFont{8}{16.8}{\rmdefault}{\mddefault}{\updefault}$u_0$}}}
\put(4325,5647){\makebox(0,0)[lb]{{\SetFigFont{8}{16.8}{\rmdefault}{\mddefault}{\updefault}$u_0$}}}
\put(5250,5647){\makebox(0,0)[lb]{{\SetFigFont{8}{16.8}{\rmdefault}{\mddefault}{\updefault}$u_0$}}}
\put(5850,6472){\makebox(0,0)[lb]{{\SetFigFont{12}{16.8}{\rmdefault}{\mddefault}{\updefault}$\hdots$}}}
\put(5850,5647){\makebox(0,0)[lb]{{\SetFigFont{12}{16.8}{\rmdefault}{\mddefault}{\updefault}$\hdots$}}}
\put(900,6622){\makebox(0,0)[lb]{{\SetFigFont{12}{16.8}{\rmdefault}{\mddefault}{\updefault}$a$}}}
\put(900,5347){\makebox(0,0)[lb]{{\SetFigFont{12}{16.8}{\rmdefault}{\mddefault}{\updefault}$b$}}}
\put(1625,5947){\makebox(0,0)[lb]{{\SetFigFont{12}{16.8}{\rmdefault}{\mddefault}{\updefault}$d$}}}
\put(1650,7072){\makebox(0,0)[lb]{{\SetFigFont{12}{16.8}{\rmdefault}{\mddefault}{\updefault}$ c$}}}
\put(900,3622){\makebox(0,0)[lb]{{\SetFigFont{12}{16.8}{\rmdefault}{\mddefault}{\updefault}$a$}}}
\put(900,2647){\makebox(0,0)[lb]{{\SetFigFont{12}{16.8}{\rmdefault}{\mddefault}{\updefault}$b$}}}
\put(1600,4222){\makebox(0,0)[lb]{{\SetFigFont{12}{16.8}{\rmdefault}{\mddefault}{\updefault}$c$}}}
\put(1575,1522){\makebox(0,0)[lb]{{\SetFigFont{12}{16.8}{\rmdefault}{\mddefault}{\updefault}$c$}}}
\put(1575,847){\makebox(0,0)[lb]{{\SetFigFont{12}{16.8}{\rmdefault}{\mddefault}{\updefault}$a$}}}
\put(1575,22){\makebox(0,0)[lb]{{\SetFigFont{12}{16.8}{\rmdefault}{\mddefault}{\updefault}$b$}}}
\put(2475,3697){\makebox(0,0)[lb]{{\SetFigFont{8}{16.8}{\rmdefault}{\mddefault}{\updefault}$u_0$}}}
\put(3375,3697){\makebox(0,0)[lb]{{\SetFigFont{8}{16.8}{\rmdefault}{\mddefault}{\updefault}$u_0$}}}
\put(4275,3697){\makebox(0,0)[lb]{{\SetFigFont{8}{16.8}{\rmdefault}{\mddefault}{\updefault}$u_0$}}}
\put(5175,3697){\makebox(0,0)[lb]{{\SetFigFont{8}{16.8}{\rmdefault}{\mddefault}{\updefault}$u_0$}}}
\put(1650,3450){\makebox(0,0)[lb]{{\SetFigFont{8}{16.8}{\rmdefault}{\mddefault}{\updefault}$v$}}}
\put(2775,3450){\makebox(0,0)[lb]{{\SetFigFont{8}{16.8}{\rmdefault}{\mddefault}{\updefault}$v$}}}
\put(3675,3450){\makebox(0,0)[lb]{{\SetFigFont{8}{16.8}{\rmdefault}{\mddefault}{\updefault}$v$}}}
\put(4575,3450){\makebox(0,0)[lb]{{\SetFigFont{8}{16.8}{\rmdefault}{\mddefault}{\updefault}$v$}}}
\put(5850,3772){\makebox(0,0)[lb]{{\SetFigFont{10}{16.8}{\rmdefault}{\mddefault}{\updefault}$\hdots$}}}
\put(5850,3022){\makebox(0,0)[lb]{{\SetFigFont{10}{16.8}{\rmdefault}{\mddefault}{\updefault}$\hdots$}}}
\put(5775,1222){\makebox(0,0)[lb]{{\SetFigFont{10}{16.8}{\rmdefault}{\mddefault}{\updefault}$\hdots$}}}
\put(5775,397){\makebox(0,0)[lb]{{\SetFigFont{10}{16.8}{\rmdefault}{\mddefault}{\updefault}$\hdots$}}}
\put(2525,1147){\makebox(0,0)[lb]{{\SetFigFont{8}{16.8}{\rmdefault}{\mddefault}{\updefault}$u_0$}}}
\put(3425,1147){\makebox(0,0)[lb]{{\SetFigFont{8}{16.8}{\rmdefault}{\mddefault}{\updefault}$u_0$}}}
\put(4350,1147){\makebox(0,0)[lb]{{\SetFigFont{8}{16.8}{\rmdefault}{\mddefault}{\updefault}$u_0$}}}
\put(5225,1147){\makebox(0,0)[lb]{{\SetFigFont{8}{16.8}{\rmdefault}{\mddefault}{\updefault}$u_0$}}}
\put(3300,4747){\makebox(0,0)[lb]{{\SetFigFont{10}{16.8}{\rmdefault}{\mddefault}{\updefault}$=$}}}
\put(3375,2047){\makebox(0,0)[lb]{{\SetFigFont{10}{16.8}{\rmdefault}{\mddefault}{\updefault}$=$}}}
\put(-300,5400){\makebox(0,0)[lb]{{\SetFigFont{10}{16.8}{\rmdefault}{\mddefault}{\updefault}$\sum\limits_{d}$}}}
\put(2075,397){\makebox(0,0)[lb]{{\SetFigFont{8}{16.8}{\rmdefault}{\mddefault}{\updefault}$u$}}}
\put(2975,397){\makebox(0,0)[lb]{{\SetFigFont{8}{16.8}{\rmdefault}{\mddefault}{\updefault}$u$}}}
\put(3875,397){\makebox(0,0)[lb]{{\SetFigFont{8}{16.8}{\rmdefault}{\mddefault}{\updefault}$u$}}}
\put(4850,397){\makebox(0,0)[lb]{{\SetFigFont{8}{16.8}{\rmdefault}{\mddefault}{\updefault}$u$}}}
%vertical arrows
\path(2350,5800)(2400,5700)(2450,5800)
\path(3250,5800)(3300,5700)(3350,5800)
\path(4150,5800)(4200,5700)(4250,5800)
\path(5050,5800)(5100,5700)(5150,5800)

\path(2350,3100)(2400,3000)(2450,3100)
\path(3250,3100)(3300,3000)(3350,3100)
\path(4150,3100)(4200,3000)(4250,3100)
\path(5050,3100)(5100,3000)(5150,3100)

\path(2350,1300)(2400,1200)(2450,1300)
\path(3250,1300)(3300,1200)(3350,1300)
\path(4150,1300)(4200,1200)(4250,1300)
\path(5050,1300)(5100,1200)(5150,1300)

% horiz lines
\path(1850,6445)(1750,6395)(1850,6345)
\path(1950,3445)(1850,3395)(1950,3345)

\end{picture}
}

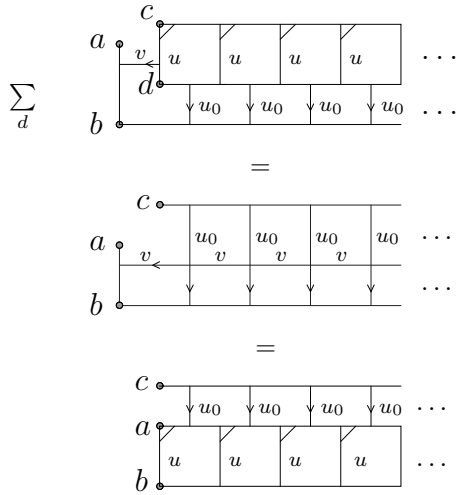
\vspace*{4mm}\begin{center}{\ft Figure 4.5: Commutation relations of the tail operator
and East half-transfer matrix}
\end{center}

\nin Expressed in terms of the North half-transfer matrix, this becomes 
after replacing $u_0 \to -u_0-1,\ v \to -v$ 
\bea  \sli_{d}
%\frac{g_d}{g_c} \, 
\lmatsk{c}{d}{a}{b}{u-u_0}
\Phi_{c,d}(u)\, \gL_{d,b}(u_0)=
%\frac{g_b}{g_a}\, 
\gL_{c,a}(u_0)\, \Phi_{a,b}(u).\lb{TailPhi}
\ee
where we define $\gL_{a,b}(u_0)\equiv g_a\Gamma^{-1} \tgL_{a,b}(-u_0-1)\Gamma g_b^{-1}$.

The other key property of the tail operator is
\bea
\Lambda_{a,a}(u_0)=\widetilde{\Lambda}_{a,a}(-u_0-1)=\id.\lb{Tailaa}
\ena 
This follows as a consequence of relation \mref{inversiona}.

\ssect{General Formula for Correlation Functions}
Now we return to the problem of computing the vertex model correlation function\\
$P^{(\ell)}(\vep_1,\vep_2,\cdots,\vep_N)$ introduced in Section 3.2.
In the infinite volume limit, we make the identification \mref{ci},
and hence compute $P^{(\ell)}(\vep_1,\vep_2,\cdots,\vep_N)$ via the equation
\bea P^{(\ell)} (\vep_1,\vep_2,\cdots,\vep_N)
= \frac{1}{ Z^{(\ell)}_m }
    F^{(\ell)}_m(\vep_1,\vep_2,\cdots,\vep_N),\lb{crf1}\ee
where $F^{(\ell)}_m(\vep_1,\vep_2,\cdots,\vep_N)$ is 
expressed as a trace over the composition of SOS model 
corner and half transfer matrices
and the tail operator shown in Figure 4.3. 
We hence find
\ben 
&&\hspace*{-8mm}F^{(\ell)}_m(\vep_1,\vep_2,\cdots,\vep_N)=
\sli_{a, a_j,a'_j} 
 \lmatk{a}{a_1}{a}{a'_1}{u_0}{\vep_1}  \lmatk{a_1}{a_2}{a'_1}{a'_2}{u_0}{\vep_1}\cdots \lmatk{a_{N-1}}{a_N}{a'_{N-1}}{a'_N}{u_0}{\vep_N}
\\
&&\hspace*{-8mm}\Tr_{\cH^{(\ell)}_{m,a}}\big(
\Phi_{a,a_1}(u)\cdots \Phi_{a_{N-1},a_N}(u) A_{NE;a_N}(u)\\ && \times \tgL_{a_N,a'_N}(u_0)A_{SE;a'_N}(u)\Phi_{S;a'_N,a'_{N-1}}(u) \cdots \Phi_{S;a'_1,a}(u) A_{SW;a}(u)A_{a}(u) \big)
\een
Here, the sum extends over
$a\in 2\Z+m+\ell$ and  $a_j,a'_j\in 2\Z+m+\sigma^j(\ell)$.
Rewriting this expression purely in terms of the North
half-transfer matrix $\Phi_{a,b}(u)$ and North-West corner-transfer
matrix $A_a(u)$, and making use of the cyclicity of the trace and 
relation \mref{Aphi} gives
\ben &&\hspace*{-8mm}F^{(\ell)}_{m,a}(\vep_1,\vep_2,\cdots,\vep_N)=
\sli_{a, a_j,a'_j} 
  \lmatk{a}{a_1}{a}{a'_1}{u_0}{\vep_1}  \lmatk{a_1}{a_2}{a'_1}{a'_2}{u_0}{\vep_1}\cdots \lmatk{a_{N-1}}{a_N}{a'_{N-1}}{a'_N}{u_0}{\vep_N}
\\
&&\hspace*{-8mm}[{a'_N}]\,
\Tr_{{\cal H}^{(\sigma^N(\ell))}_{m,a'_N}}\big(A_{a'_N}(u) A_{a'_N}(-1-u) A_{a'_N}(u) A_{a'_N}(-1-u) \\&&\hspace*{-8mm}\times
\Phi^*_{a'_N,a'_{N-1}}(-1)\cdots \Phi^*_{a'_1,a}(-1)\Phi_{a,a_1}(-1) \cdots
\Phi_{a_{N-1},a_N}(-1)
\gL_{a_N,a'_N}(-u_0-1) \big)
\een

The factor $Z^{(\ell)}_m $ in the denominator of \mref{crf1} is the partition function 
expressed by Figure 4.3 with $N=0$. This is the same as the SOS partition function 
given in \eqref{charpart} due to \eqref{Tailaa}. 
Using the infinite volume limit $A_a(u)\sim x^{-2uH^{(\ell)}_{m,a}}$ of
the SOS corner transfer matrix, we finally have 
\bea 
&&\hspace*{-10mm}P^{(\ell)}(\vep_1,\vep_2,\cdots,\vep_N)\nn\\
&&\hspace*{-10mm}=\frac{1}{[m]^*\chi^{(k)}_\ell(\bar{\tau})}
%\nn \\&&\times
\sli_{a, a_j,a'_j} 
  \lmatk{a}{a_1}{a}{a'_1}{u_0}{\vep_1}\lmatk{a_1}{a_2}{a'_1}{a'_2}{u_0}{\vep_1}\cdots \lmatk{a_{N-1}}{a_N}{a'_{N-1}}{a'_N}{u_0}{\vep_N}\nn \\
%[3mm]
&&
\hspace*{-10mm}
\times[a'_N]\,
Tr_{\cH^{(\sigma^N(\ell))}_{m,a'_N}}\big(x^{4H^{(\ell)}_{m,a}}
\Phi^*_{a'_N,a'_{N-1}}(-1)\cdots \Phi^*_{a'_1,a}(-1)\Phi_{a,a_1}(-1) \cdots
\Phi_{a_{N-1},a_N}(-1)
\gL_{a_N,a'_N}(-u_0-1) \big).\nn\\
\lb{gencorr}
\ena
The vertical spectral parameter $u_0$ introduced in the vertex model still
appears in this expression. However, we expect that after taking the
trace and all the summations that there will in fact be 
no $u_0$ dependence \cite{LaP98}. A more detailed
discussion of this point
will be given in \cite{KKW04b}.

\ssect{The Vertex-Face Correspondence for Lattice Operators}
From Figure 4.3, we can also extract the following vertex-face correspondence 
for lattice operators such as the half transfer matrices, the corner transfer matrix 
and the space of states.
\subsubsection{Half transfer matrices}
Let us define $\Phi^{(\bell,\ell)}_\vep(u;u_0)$ and 
$\Phi^{*(\bell,\ell)}_\vep(u;u_0)$ by
\bea
\Phi^{(\bell,\ell)}_\vep(u;u_0)&=&f^{(k)}(u-u_0)\bigoplus_{a\in m+\ell+2\Z \atop a'\in a+\{-k,-k+2,..,k\}} 
\psi^{(k)}_\vep(u-u_0)^{a'}_a\Phi_{a',a}(u),\lb{vfvertex}\\
\Phi^{*(\bell,\ell)}_\vep(u;u_0)&=&\sum_{\vep'}Q^{\vep'}_{\vep}\Phi^{
(\bell,\ell)}_{\vep'}(u-1;u_0)\nn\\
%&=&f^{(k)}(u-u_0-1)
%\bigoplus_{a'\in m+k-l+2\Z ? \atop a\in m+l+2\Z} \sum_{\vep'}Q^{\vep'}_{\vep}
%\psi^{(k)}_{\vep'}(u-u_0-1)^{a'}_a\Phi_{a',a}(u-1)\nn\\
&=&f^{(k)}(u-u_0)^{-1}\bigoplus_{a\in m+\ell+2\Z \atop a'\in a+\{-k,-k+2,..,k\}}
\psi^{*(k)}_\vep(u-u_0)^{a'}_a\Phi^*_{a',a}(u).\lb{vfdualvertex}
\ena
Here the function $f^{(k)}(u)$ is chosen such that 
\bea
&&f^{(k)}(u)f^{(k)}(u-1)= C^{(k)}(u),\lb{condf1}\\
&&f^{(k)}(u-2)=\frac{[u+k-1]}{[u-1]}f^{(k)}(u),\lb{condf2}
\ena
where $C^{(k)}(u)$ is a function appearing in \eqref{crosspsi}.
If we consider $u_0$ as a fixed constant and suppress it from the
notation from these lattice operators, then we have the following
theorem.

\begin{thm}
\begin{itemize}
\item[(i)]The lattice operator $\Phi^{(\bell,\ell)}_\vep(u)$
 satisfies the commutation relation \eqref{Rphiphi}. 
\item[(ii)] $\Phi^{(\bell,\ell)}_\vep(u)$ and $\Phi^{*(\ell,\bell)}_\vep(u)$ 
satisfy the inversion relation \eqref{8vinversion}.
\end{itemize}
\end{thm}
\noi{\it Proof)}\quad    
The commutation relation \eqref{Rphiphi} follows from \eqref{fusionvertexface}and \eqref{commrel}, whereas the inversion relation \eqref{8vinversion} 
follows from \eqref{fusioninv1} and \eqref{Eqn:Inversion}.  

%\vspace{3mm}
\subsubsection{ The corner transfer matrix}
Let us define
\be
&&\rho^{(\ell)}=\bigoplus_{a' \atop a\in m+l+2\Z} 
\Lambda_{a',a}(u_0)\ [a] x^{4H^{(\ell)}_{m,a}}.
\en
%Then 
\begin{thm}
$\rho^{(\ell)}$ and $\Phi^{(\bell,\ell)}_\vep(u;u_0)$ satisfy 
\bea
&&\rho^{(\bell)}\Phi^{(\bell,\ell)}_\vep(u-2;u_0)=\Phi^{(\bell,\ell)}_\vep(u;u_0)
\rho^{(\ell)}.\lb{rhophiSOS}
\ena
\end{thm}
\noi This should be compared with \eqref{rhophi}.

\noi{\it Proof)}\quad   
Multiply \eqref{TailPhi} by $f^{(k)}(u-u_0-2)\psi^{'(k)}_\vep(u-u_0)^a_b$ 
and take 
the summation in $a$. Then from \eqref{compopsip}, \eqref{fusioninv2p} and \eqref{condf2}, we obtain 
\be
&&\sum_a \Lambda_{c,a}(u_0)[a] f^{(k)}(u-u_0-2)\psi^{(k)}_{\vep}(u-u_0-2)^a_b
\Phi_{a,b}(u)\nn\\
&&=\sum_a  f^{(k)}(u-u_0)\psi^{(k)}_{\vep}(u-u_0)^c_a\Phi_{c,a}(u)
\Lambda_{a,b}(u_0)[b].
\en
Now act with $x^{4H^{(\ell)}_{m,a}}$ on right of this expression and use \eqref{cornerPhi}. Then taking direct sum 
$\displaystyle{\bigoplus_{b,c}}$ of the result, we have 
\be
&&\sum_a \bigoplus_{c}\Lambda_{c,a}(u_0)[a] x^{4H^{(\bell)}_{m,a}}\cdot \bigoplus_{b}
f^{(k)}(u-u_0-2)\psi^{(k)}_{\vep}(u-u_0-2)^a_b\Phi_{a,b}(u-2)\nn\\
&&=\sum_a  \bigoplus_{c}f^{(k)}(u-u_0)\psi^{(k)}_{\vep}(u-u_0)^c_a\Phi_{c,a}(u)
\cdot \bigoplus_{b}\Lambda_{a,b}(u_0)[b]x^{4H^{(\ell)}_{m,a}}.
\en 
We hence obtain \eqref{rhophiSOS}.

%\vspace{3mm}
\subsubsection{ The space of states} 
Figure 4.3 with $N=0$ gives the partition function 
for the dressed vertex model. 
%\be
%&&\sum_{a}[a]\Tr_{\cH_{m,a}^{(\ell)}}\Lambda_{a,a}(-u_0-1)\ x^{4H_{m,a}^{(\ell)}} 
%\en
From \eqref{Tailaa}, this coincides with the SOS partition function $Z^{(\ell)}_m$ 
as mentioned above. 
Then from Theorem \ref{charpartSOS}, we have 
\be
&&Z^{(\ell)}_m=[m]^* Z^{(\ell)}
\en
The RHS is the partition 
function of the $k\times k$ fusion eight-vertex model multiplied by $[m]^*$. 
The multiplicity by the factor $[m]^*$ can be regarded as the effect of 
dressing the boundary of the vertex model.
Understanding this multiplicity and considering \eqref{rhophiSOS}, we can roughly 
make the following vertex-face correspondences for the spaces of states
and 
corner transfer matrices.
\be
&&\H^{(\ell)} \ \longleftrightarrow \ \bigoplus_{a\in \ell+m+2\Z}\  \H^{(\ell)}_{m,a}, 
%\\&&
\qquad x^{4H^{(\ell)}} \ \longleftrightarrow \ \rho^{(\ell)}
\en
%%%%%%%%%%%%%%%%%%%%%%%%%%%%%%%%%%%%%%%%%%%%
\section{The Realization of Fusion SOS Models in Terms of 
The Elliptic Algebra $U_{x,p}(\widehat{\goth{sl}}_2)$}
%%%%%%%%%%%%%%%%%%%%%%%%%%%%%%%%%%%%%%%%%%%%
In this section, 
we show how the space of the states ${\cal H}_{m,a}^{(\ell)}$, the corner Hamiltonian 
$H_{m,a}^{(\ell)}$, the half transfer matrix $\Phi_{b,a}(u)$ and the tail operator 
${\Lambda}_{b,a}(u)$ appearing in \eqref{gencorr} can all be constructed in terms of the representation theory of the elliptic algebra $U_{x,p}(\slth)$. One of the key results is the remarkably simple algebraic form of the tail operator given by Conjecture \ref{conj2}. 

\subsection{The Elliptic Algebra $U_{x,p}(\widehat{\goth{sl}}_2)$}
We first present a brief review of the elliptic algebra
$U_{x,p}(\widehat{\goth{sl}}_2)$ and its associated vertex 
operators 
\cite{Ko98,JKOS}\footnote{$U_{x,p}(\widehat{\goth{sl}}_2)$ is conventionally
referred to as $U_{q,p}(\widehat{\goth{sl}}_2)$ - the difference
is merely a change of notation.} . The elliptic algebra
$U_{x,p}(\widehat{\goth{sl}}_2)$ provides the Drinfeld realization of
the face-type elliptic quantum group 
${\cal B}_{x,\lambda}(\widehat{\goth{sl}}_2)$ 
\cite{JKOS1} tensored
by a Heisenberg algebra. 

\subsubsection{The Definition and a realization of $U_{x,p}(\slth)$}

\begin{df}\cite{Ko98,JKOS}~{\bf
(The Elliptic Algebra $U_{x,p}(\widehat{\goth{sl}}_2)$)}~~
The elliptic algebra $U_{x,p}(\widehat{\goth{sl}}_2)\ (p=x^{2r},\ r\in \C)$
is the associative algebra of the currents
$E(u), F(u)$, $K(u)$ and the grading operator $\widehat{d}$  
satisfying the following relations:
\ben
K(u_1)K(u_2)&=&\rho(u_1-u_2)K(u_2)K(u_1),\\
K(u_1)E(u_2)&=&\frac{[u_1-u_2+\frac{1-r^*}{2}]^*}
{[u_1-u_2-\frac{1+r^*}{2}]^*}E(u_2)K(u_1),\\
K(u_1)F(u_2)&=&\frac{[u_1-u_2-\frac{1+r}{2}]}{
[u_1-u_2+\frac{1-r}{2}]}F(u_2)K(u_1),\\
E(u_1)E(u_2)&=&\frac{[u_1-u_2+1]^*}{
[u_1-u_2-1]^*}E(u_2)E(u_1),\\
F(u_1)F(u_2)&=&\frac{[u_1-u_2-1]}{
[u_1-u_2+1]}F(u_2)F(u_1),\\\een
\ben
&&[\widehat{d},E(u)]=\left(-z\frac{\partial}{\partial z}+\frac{1}{r^*}\right)E(u),
\quad  [\widehat{d},F(u)]=\left(-z\frac{\partial}{\partial z}+\frac{1}{r^*}\right)F(u),
\\
&&[E(u_1),F(u_2)]=\frac{1}{x-x^{-1}}
\left(\delta(x^{-k}z_1/z_2)
H^+\left(u_2+\frac{k}{4}\right)-
\delta(x^{k}z_1/z_2)
H^-\left(u_2-\frac{k}{4}\right)
\right),\nonumber\\
\\
&&H^\pm(u)=\kappa
K\left(u\pm
\left(\frac{r}{2}-\frac{k}{4}\right)+\frac{1}{2}
\right)K\left(
u\pm
\left(\frac{r}{2}-\frac{k}{4}\right)-\frac{1}{2}
\right).
\een
Here $r^*=r-k$, $z=x^{2u}, z_i=x^{2u_i}\ (i=1,2)$ and $\delta(z)=\sum_{n \in {\mathbb{Z}}}z^n$. The constant $\kappa$ is given by
\ben
\kappa&=&\frac{\xi(x^{-2};p^*,x)}{
\xi(x^{-2};p,x)}\quad \hb{with}\quad 
\xi(z;p,x)=\frac{(x^2z;p,x^4)_\infty
(px^2z;p,x^4)_\infty}
{(x^4z;p,x^4)_\infty 
(pz;p,x^4)_\infty},
\een
and the scalar function
$\rho(v)$ is given by
\bea
\rho(v)&=&\frac{\rho^{+*}(v)}{
\rho^+(v)},\quad\hb{with}\quad
\rho^+(v)=z^{\frac{1}{2r}}x^{\frac{1}{2}}
\frac{(px^2z;p,x^4)_\infty^2
(z^{-1};p,x^4)_\infty
(x^4z^{-1};p,x^4)_\infty
}{(pz;p,x^4)_\infty
(px^4z;p,x^4)_\infty
(x^2z^{-1};p,x^4)_\infty^2}.
\lb{rho}\ee
The $*$ symbols  
always
indicates the replacement $r\ra r^*$. For example, $p^*=x^{2r^*}$, $[u]^*=x^{\frac{u^2}{r^*}-u}\Theta_{x^{2r^*}}(x^{2u})$.
\end{df}

The algebra 
$U_{x,p}(\widehat{\goth{sl}}_2)$
is realized by tensoring
 $U_{x}(\widehat{\goth{sl}}_2)$ and a
Heisenberg algebra \cite{JKOS}. 
For this realization, it is convenient to introduce the Drinfeld realization of $U_x(\slth)$.
\begin{df}~~
{\bf (The Drinfeld Realization of $U_x(\slth)$)}
The quantum affine algebra $U_x(\slth)$ is the associative algebra generated by 
 $h, a_m, x_n^\pm (m \in {\mathbb{Z}}_{\neq 0}, n 
\in {\mathbb{Z}}), d$ and the central element $k$ satisfying the relations 
\ben
&&~[h,d]=0,~[d,a_n]=na_n,~[d,x_n^\pm]=nx_n^\pm,\\
&&~[h,a_n]=0,~[h,x^\pm(z)]=\pm 2 x^\pm(z),\\
&&~[a_n,a_m]=\frac{[\![2n]\!]_x[\![kn]\!]_x}{n}x^{-k|n|}\delta_{n+m,0},\\
&&~[a_n,x^+(z)]=\frac{[\![2n]\!]_x}{n}x^{-k|n|}z^n x^+(z),\\
&&~[a_n,x^-(z)]=-\frac{[\![2n]\!]_x}{n}z^n x^-(z),\\
&&~(z-x^{\pm 2}w)x^\pm(z)x^\pm(w)=
(x^{\pm 2}z-w)x^\pm(w)x^\pm(w),\\
&&~[x^+(z),x^-(w)]=\frac{1}{x-x^{-1}}
\left(\delta(x^{-k}z/w)\psi(x^{k/2}w)-
\delta(x^{k}z/w)\varphi(x^{-k/2}w)
\right),
\een
where  $x^\pm(z), \psi(z)$ and $\varphi(z)$ denote the Drinfeld currents
 defined by 
\ben
x^\pm(z)&=&\sum_{n\in {\mathbb Z}}x^\pm_n z^{-n},\\
\psi(x^{k/2}z)&=&x^h
\exp\left(
(x-x^{-1})\sum_{n>0}a_n z^{-n}\right),\\
\varphi(x^{-k/2}z)&=&x^{-h}
\exp\left(-(x-x^{-1})\sum_{n>0}a_{-n}z^n\right).
\een
\end{df}
Let us denote by ${\mathbb{C}}\{\hat{\cal H}\}$
 the Heisenberg algebra generated by the pair 
$P, Q$ with $[Q,P]=1$. Then we have the following realization of $U_{x,p}(\slth)$.
\begin{thm}\lb{UxHeisen}\cite{JKOS}~~
The elliptic algebra $U_{x,p}(
\widehat{\goth{sl}}_2)$ is realized by tensoring
$U_{x}(\widehat{\goth{sl}}_2)$ 
and the Heisenberg algebra
${\mathbb{C}}\{\hat{\cal H}\}$.
The generators $E(u), F(u)$, $K(u)$ and $\widehat{d}$ are given by
\begin{eqnarray*}
K(u)&=&k(z)e^Q  z^{(\frac{1}{r}-\frac{1}{r^*})\frac{P}{2}+
\frac{h}{2r}+\frac{1}{4}(\frac{1}{r^*}-\frac{1}{r})},\\
E(u)&=&u^+(z,p)x^+(z) e^{2Q} z^{\frac{1}{r^*}(-P+1)},\\
F(u)&=&x^-(z)u^-(z,p) z^{\frac{1}{r}(P+h-1)},\\
\widehat{d}&=&d-\frac{1}{4r^*}(P-1)(P+1)+\frac{1}{4r}(P+h-1)(P+h+1),
\end{eqnarray*}
where $k(z)$ and $u^\pm(z,p)$ are given by
\begin{eqnarray*}
&&k(z)=\exp\left(\sum_{n>0}\frac{[\![n]\!]_x}{[\![2n]\!]_x [\![r^*n]\!]_x}a_{-n}
(x^kz)^n\right)
\exp\left(-\sum_{n>0}\frac{[\![n]\!]_x}{[\![2n]\!]_x [\![rn]\!]_x}a_nz^{-n}\right),\\
&&u^+(z,p)=\exp\left(\sum_{n>0}\frac{1}{[\![r^*n]\!]_x}a_{-n}(x^rz)^n
\right),\quad
u^-(z,p)=\exp\left(
-\sum_{n>0}\frac{1}{[\![rn]\!]_x}a_{n}(x^{-r}z)^{-n}\right).
\end{eqnarray*}
\end{thm}

\noi
The following commutation relations are important.
\begin{prop}
\begin{eqnarray}
&&~[K(u),P]=K(u),~[E(u),P]=2E(u),~[F(u),P]=0,\\
&&~[K(u),P+h]=K(u),~[E(u),P+h]=0,~[F(u),P+h]=2F(u).
\end{eqnarray}
\end{prop}

Next we define the `half currents' and the $L$-operator.
\begin{df}~\cite{JKOS}~{\bf (Half Currents)}
~~We define the half currents
$E^+(u), F^+(u)$ and $K^+(u)$ by
\begin{eqnarray*}
K^+(u)&=&K\left(u+\frac{r+1}{2}\right),\quad  
E^+(u)= \widehat{E}^+(u)
\frac{[1]^*}{[P-1]^*},\quad  
F^+(u)= \widehat{F}^+(u)\frac{[1]}{[P+h-1]},
\end{eqnarray*}
where 
\be
&&\widehat{E}^+(u)=a^* \oint_{C^*}\frac{dw}{2\pi i w} E(v)
\frac{[u-v+\frac{k}{2}-P+1]^*}{
[u-v+\frac{k}{2}]^*},\\
&&\widehat{F}^+(u)=a \oint_{C}\frac{dw}{2\pi i w}
F(v)\frac{[u-v+P+h-1]}{[u-v]}.
\en
The contours $C^*$ and $C$ are defined as follows 
\begin{eqnarray*}
C^*:|p^*x^k z|<|w|<|x^kz|,~~~C:|pz|<|w|<|z|,
\end{eqnarray*}
and the constant $a,a^*$ are chosen to satisfy
\begin{eqnarray*}
\frac{a^* a [1]^* \kappa}{x-x^{-1}}=1.
\end{eqnarray*}
\end{df}

\begin{df}~\cite{JKOS}~{\bf (L-operator)}
~~We define the $L$-operator
$\widehat{L}^+(u) \in {\rm End}({\mathbb{C}}^2)\otimes
U_{x,p}(\widehat{\goth{sl}}_2)$ by
\begin{eqnarray*}
\widehat{L}^+(u)=
\left(\begin{array}{cc}
1&F^+(u)\\
0&1
\end{array}
\right)
\left(\begin{array}{cc}
K^+(u-1)&0\\
0&K^+(u)^{-1}
\end{array}
\right)
\left(\begin{array}{cc}
1&0\\
E^+(u)&1
\end{array}
\right).\lb{Lop}
\end{eqnarray*}
\end{df}

\subsubsection{The vertex operators of ${U}_{x,p}(\widehat{\goth{sl}}_2)$}

Let $V(\la_\ell)$ be the level $k$ irreducible $U_x(\slth)$-module 
 with highest weight $\la_\ell$. 
Let us denote by 
$(\pi_{n,z}, V_{n,z})\ (n=0,1,..,k)$ the $n+1$-dimensional evaluation representation of
$U_x(\widehat{\goth{sl}}_2)$: $V_{n,z}=V^{(n)}\otimes
{\mathbb{C}}[z,z^{-1}],~V^{(n)}=\oplus_{m=0}^n  {\mathbb{C}}\,v^{(n)}_m$. 
In physical applications, we only have to consider the case $n=k$. 
The co-algebra structure of 
${\cal B}_{x,\lambda}(\widehat{\goth{sl}}_2)$
allows us to define the type I intertwining operator 
$\Phi(u,P) : V(\la_\ell)\to V(\la_{k-\ell})\otimes V_{k,z}$ of 
${\cal B}_{x,\lambda}(\widehat{\goth{sl}}_2)$-modules. Note that 
${\cal B}_{x,\lambda}(\widehat{\goth{sl}}_2)\cong U_x(\slth)$.

Now let us consider the $U_{x,p}(\slth)$-modules.
According to the realization of $U_{x,p}(\slth)$ given by 
Theorem \ref{UxHeisen}, we define the the level $k$  
$U_{x,p}(\slth)$-module  $\widehat{ V}(\la_\ell)$ 
by 
\be
&&\widehat{ V}(\la_\ell)=\bigoplus_{m \in {\mathbb Z}}
V(\lambda_\ell) \otimes e^{-m Q}.
\en 
The type I vertex operator $\widehat{\Phi}(u)$  
of  ${U}_{x,p}(\widehat{\goth{sl}}_2)$ is simply 
defined to be the same as $\Phi(u,P)$: 
\be
&&\widehat{\Phi}(u)=\Phi(u,P) : \widehat{V}(\la_\ell)\to \widehat{V}(\la_{k-\ell})
\otimes V_{k,z}.
\en
Throughout this paper we consider only the type I vertex operator. 

From the intertwining relation for $\Phi(u,P)$, we obtain
 the following relation which characterises the vertex operator uniquely up to a 
normalisation factor:
\begin{eqnarray}
\widehat{\Phi}(u_2)
\widehat{L}^+(u_1)=
R_{1,k}^{+(13)}(u_1-u_2,P+h)\widehat{L}^+(u_1)
\widehat{\Phi}(u_2).
\lb{int1}
\end{eqnarray}
Here $R_{1k}^{+}(u,s)$ is an image of the $L$-operator:
$R_{1k}^+(u-v,s)=(\id \otimes \pi_{k,w})L^+(u,s)$ with $z=x^{2u}, w=x^{2v}$.
The finite dimensional representations of the elliptic currents as well as 
 the expression for the matrix $R_{1k}^{+}(u,s)$ can be found in 
Appendix C of \cite{JKOS}.
Inputting the realization of $\widehat{L}^+(u)$ given by
Definition \ref{Lop} 
into \eqref{int1}, we can solve \eqref{int1} for $\widehat{\Phi}(u)$. 
Let us define the components of the vertex operator 
$\widehat{\Phi}(u)$ as follows:
\ben
\widehat{\Phi}\left(u-\frac{1}{2}\right)=\sum_{m=0}^k
{\Phi}_{k,m}(u)\otimes v_m.
\een
Then we obtain the following realization of the vertex operators\cite{JKOS}.
\begin{thm}\lb{realvertex}
The highest component ${\Phi}_{k,k}(u)$ is given by  
\ben
&&{\Phi}_{k,k}(u)=:\exp\left\{
%\sum_{n\not=0}\frac{\al_m}{[\![2m]\!]_x}z^{-n}
\sum_{n\not=0} \frac{ [\![r^*n]\!]_x \alpha_n}{[\![2n]\!]_x[\![rn]\!]_x} z^{-n}
\right\}: e^{\frac{k}{2}\al}(-z)^{\frac{h}{2}}z^{-\frac{k}{2r}(P+h)},
\en
where\\[-12mm]
\ben 
&& \alpha_n=\left\{\mmatrix{a_n&{\rm for}\ws n>0\cr
  \frac{[\![rn]\!]_x}{[\![r^*n]\!]_x}x^{k|n|}a_n&{\rm for}\ws n<0}\right. ,\\[1mm]
&& \hb{such that}\ws
[\alpha_m,\alpha_n]=\delta_{m+n,0}\frac{[\![2m]\!]_x[\![km]\!]_x}{m}
\frac{[\![rm]\!]_x}{[\![r^*m]\!]_x}.\een
For the remaining components $m=0,1,\cdots,k$, we have the formula
\begin{eqnarray}
{\Phi}_{k,m}(u)&=&{\Phi}_{k,k}(u)
\widehat{F}^+
\left(u+\frac{k}{2}-r
\right)^{k-m}\,\,
\prod_{j=1}^{k-m}
\frac{[1][P+h+2k-m-1+2j]}{
[P+h+k-1+2j][P+h+k+2j]}\nonumber\\
&=&\oint_{C_1} \frac{dw_1}{2\pi i w_1}\cdots
\oint_{C_{k-m}} \frac{dw_{k-m}}{2\pi i w_{k-m}}\,
{\Phi}_{k,k}(u)F(v_1)\cdots F(v_{k-m})\nonumber\\
&\times&\prod_{j=1}^{k-m}
\frac{[u-v_j+P+h+\frac{k}{2}-1-2(k-m-j)][1][P+h+2k-m-1+2j]}
{[u-v_j+\frac{k}{2}][P+h+k-1+2j][P+h+k+2j]},\nn\\ \lb{vertexcomp}
\end{eqnarray}
where $z=x^{2u}, w_j=x^{2v_j}$ and  the integral contours $C_j~~(1\leq j \leq k-m)$ 
are given by
\begin{eqnarray*}
&&C_{1}:~|x^kz|<|w_{1}|<|p^{-1}x^kz|, |x^{-k}z|,\\
&&C_j:~|x^kz|<|w_{j}|<|p^{-1}x^kz|, |x^{-k}z|, |x^2w_{j-1}|,~~
(2 \leq j \leq k-m).
\end{eqnarray*}
\end{thm}
\noi The expression of $\Phi_{k,k}(u)$ is equivalent but slightly different from the 
one given in \cite{JKOS} in the zero-mode part. 
Note the following properties.
\begin{prop}
\begin{eqnarray}
&&[P,{\Phi}_{k,k}(u)]=0,
~[h,{\Phi}_{k,k}(u)]=k {\Phi}_{k,k}(u),\nn\\
&&F(v){\Phi}_{k,k}(u)=
\frac{[u-v-\frac{k}{2}]}{[u-v+\frac{k}{2}]}{\Phi}_{k,k}(u)F(v),\lb{PhiF}\\
&&E(v){\Phi}_{k,k}(u)={\Phi}_{k,k}(u)E(v).\nn
\ena
\end{prop}

Let us next consider the commutation relations of the 
vertex operators of $U_{x,p}(\slth)$. In \cite{JKOS1}, 
they are expected to be commutation 
relations with exchange coefficients being exactly the fused face weights 
$W^{(k,k)}$ \eqref{kkfusionW}. In order to derive such relations, we
 consider the following gauge transformation of the vertex operators 
${\Phi}_{k,m}(u) \ra \widehat{\Phi}_{\vep}(u)$ with $\vep=2m-k$.
\be
\widehat{\Phi}_{\vep}(u)&=&\Phi_{k,m}(u)
\prod_{j=1}^{k-m}
\frac{[P+h+k-1+2j][P+h+k+2j]}{[1][P+h+2k-m-1+2j]}\quad (\vep=-k,-k+2,..,k).
\en
From \eqref{vertexcomp}, we find 
\bea
\widehat{\Phi}_{\vep}(u)={\Phi}_{k,k}(u)
\widehat{F}^+
\left(u+\frac{k}{2}-r
\right)^{\frac{k-\vep}{2}}.\lb{fvo}
\ena
By using the realisation obtained in Theorem 
\ref{realvertex} and the commutation relation \eqref{PhiF}, 
we have checked the following commutation relations 
for levels $k=1,2,3$. 

\begin{conj}~~{\bf (Commutation relations)}~~The vertex
operators $\widehat{\Phi}_\vep(u)$, \\$(\vep=-k,-k+2,\cdots, k)$ 
satisfy the commutation
relations
\begin{eqnarray}
&&\hspace*{-20mm}\widehat{\Phi}_{\vep_2}(u_2)
\widehat{\Phi}_{\vep_1}(u_1)\nn\\
&&\hspace*{-20mm}=
\sum_{\nu_1\in\{-k,-k+2,\cdots, k\}\atop \nu_1+\nu_2=\vep_1+\vep_2}
W^{(k,k)}\left(\left.
\begin{array}{cc}
P+h-\vep_1-\vep_2
&P+h-\vep_2\\
P+h-\nu_1&P+h
\end{array}
\right|u_1-u_2\right)
\widehat{\Phi}_{\nu_1}(u_1)
\widehat{\Phi}_{\nu_2}(u_2).
%\nonumber\\[-2mm]
\label{eqn:VO1}
\end{eqnarray}
where the 
coefficients 
$W^{(k,k)}$ are given by \mref{kkfusionW}.
\end{conj}

\subsection{The Realisation of Fusion SOS Models and the Tail Operator}
Now we formulate the fusion SOS model in terms of the 
representation theory of $U_{x,p}(\slth)$ and give a realisation of 
the space of states $\H^{(\ell)}_{m,a}$, the corner Hamiltonian $H^{(\ell)}_{m,a}$, 
the half transfer matrix $\Phi_{b,a}(u)$ and the tail operator 
${\Lambda}_{b,a}(u)$.

\subsubsection{The space of states and the CTM Hamiltonian}
We first show that the level $k$ $U_{x,p}(\widehat{\goth{sl}}_2)$-modules have a 
natural decomposition 
into the  Virasoro highest-weight modules associated with the coset 
$(\slth)_k \oplus (\slth)_{r-k-2} /(\slth)_{r-k}$. 
Irreducible Virasoro modules are identified with the 
spaces of states $\H^{(\ell)}_{m,a}$ of the $k\times k$ fusion SOS model. 

In order to see such a decomposition, 
% of $U_{x,p}$-modules into the Virasoro modules, 
it is convenient to realize the level $k$ 
%Drinfeld currents of 
$U_x(\slth)$-module $V(\la_\ell)$ 
in terms of a $q$-deformed $\Z_k$-parafermion module $\H^{PF}_{\ell,M}$ and 
the Fock module $\F^a$ of the Drinfeld bosons $a_n$\cite{Mat94,Ko98}
(see also \cite{GeQiu} for the CFT case). 

The $q$-deformed $\Z_k$-parafermion algebra is conveniently introduced through the  
$q$-deformed  $\cZ$-algebra associated with the level 
$k$ Drinfeld currents of $U_x(\slth)$. The algebraic structure of the $q$-deformed $\cZ$-algebra 
is quite parallel to the classical case\cite{LePr85}.
The 
$q$-deformed case was considered in \cite{Jing}. 
The deformed 
$\cZ$-algebra is generated by ${\cal Z}_{\pm,n} \ (n\in \Z)$ whose generating 
functions $\cZ_\pm(z)=\sum_{n\in\Z}\cZ_{\pm,n} z^{-n}$ are defined by 
\be
&&\cZ_+(z)=\EXP{-\sum_{n>0}\frac{1}{[\![kn]\!]_x}a_{-n} z^{n}}\, x^+(z)\, \EXP{\sum_{n>0}\frac{1}{[\![kn]\!]_x}a_{n} z^{-n}},\nn\\
&&\cZ_-(z)=\EXP{\sum_{n>0}\frac{x^{kn}}{[\![kn]\!]_x}a_{-n} z^{n}}\, x^-(z)\, 
\EXP{-\sum_{n>0}\frac{x^{kn}}{[\![kn]\!]_x}a_{n} z^{-n}}. \lb{zalg}
\en
The $\cZ$-algebra commutes with the Drinfeld bosons $a_n, \ n\not=0$. 
Then the  level $k$ highest-weight $U_q(\slth)$-module $V(\la_\ell)$ with 
highest weight $\la_\ell$ has the structure 
\bea
&&V(\la_\ell)=\F^a\otimes \Omega_\ell,\lb{VlatoOmega}
\ena
where $\F^a=\C[a_{-n}\ (n>0)]$.  
The space $\Omega_\ell$ is called the vacuum space defined by 
\be
&&\Omega_\ell=\left\{v\in V(\la_\ell)\ \bigl|\ a_n v=0\ (n>0)  \right\}.
\en
The space $\Omega_\ell$ is spanned by the vectors $v_\ell(\vep_1,..,\vep_s;n_1,..,n_s)\ (s\geq 0, \vep_j\in \{\pm\}, n_s\leq0,n_{s-1}+n_s\leq 0,..,n_1+\cdots+n_s\leq0 )$ given by 
\be
&&\prod_{1\leq i<j\leq s} \left[\frac{(x^{-2}x^{k+\frac{\vep_i+\vep_j}{2}k};x^{2k})_\infty}{(x^2x^{k+\frac{\vep_i+\vep_j}{2}k};x^{2k})_\infty}\right]^{\vep_i\vep_j}\cZ_{\ep_1}(z_1)\cdots \cZ_{\ep_s}(z_s)\cdot 1\otimes e^{\frac{\ell}{2}\al}\nn\\
&&\qquad\quad =\sum_{n_1,..,n_s\in\Z}
v_\ell(\ep_1,..,\ep_s;n_1,..,n_s)z_1^{-n_1}\cdots z_s^{-n_s}.
\en
Here the action of $\cZ_{\pm,n}$ is defined as follows.
\be
&&\cZ_{\pm,n}\cdot (f\otimes e^{{\ell}\frac{\alpha}{2}})
=\left\{\mmatrix{\cZ_{\pm,n}f\otimes e^{{\ell}\frac{\alpha}{2}}&n\leq 0
\cr
                 [\cZ_{\pm,n} ,f]\otimes e^{{\ell}\frac{\alpha}{2}}&n\geq 1}  \right.  
\en
for $f\in\C[\cZ_{+,n}, \cZ_{-,n}\ (n\leq 0)]$.
The weight of $v_\ell(\vep_1,..,\vep_s;n_1,..,n_s)$ is $\la_\ell+\sum_{j=1}^s\vep_j \al$ 
and its degree is $-\frac{\ell(\ell+2)}{4(k+2)}+n_1+\cdots+n_s$. 

Now let us consider the $q$-deformed $\Z_k$-parafermion. 
Define the basic $\Z_k$-parafermion currents $\Psi(z)$ and $\Psi^{\dagger}(z)$ through 
the following relations.
\be
&&\cZ_{+}(z)=\Psi(z)\otimes e^{\al}z^{\frac{h}{k}},\\
&&\cZ_{-}(z)=\Psi^\dagger(z)\otimes e^{-\al}z^{-\frac{h}{k}},\\
&&[\Psi(z),\al]=[\Psi(z),h]=[\Psi^\dagger(z),\al]=[\Psi^\dagger(z),h]=0.
\en
To make this expression well-defined, 
$\Psi(z)$ and $\Psi^{\dagger}(z)$ should have their  mode expansions 
depending on the weight of vectors on which they act. Namely, on the vector 
with weight $\la$ such that $(h,\la)=m$, we have 
\be
&&\Psi(z)\equiv\Psi^+(z)=\sum_{n\in \Z}\Psi_{+,\frac{m}{k}-n} z^{-\frac{m}{k}+n-1},\\
&&\Psi^\dagger(z)\equiv\Psi^-(z)=\sum_{n\in \Z}\Psi_{-,\frac{m}{k}-n} 
z^{\frac{m}{k}+n-1}.
\en
The  $q$-deformed $\Z_k$-parafermion algebra is generated by 
$\Psi_{+,\frac{m}{k}-n},\ \Psi_{-,\frac{m}{k}-n} \ 
( n\in \Z) $. Its relations can be expressed as follows.
\be
&&\left(\frac{z}{w}\right)^{\frac{2}{k}}\frac{(x^{-2}w/z;x^{2k})_{\infty}}
{(x^{2+2k}w/z;x^{2k})_{\infty}}\Psi^{\pm}(z)\Psi^{\pm}(w)=
\left(\frac{w}{z}\right)^{\frac{2}{k}}\frac{(x^{-2}z/w;x^{2k})_{\infty}}
{(x^{2+2k}z/w;x^{2k})_{\infty}}\Psi^{\pm}(w)\Psi^{\pm}(z),\\
&&\left(\frac{z}{w}\right)^{-\frac{2}{k}}\frac{(x^{2+k}w/z;x^{2k})_{\infty}}
{(x^{-2+k}w/z;x^{2k})_{\infty}}\Psi^{\pm}(z)\Psi^{\mp}(w)-
\left(\frac{w}{z}\right)^{-\frac{2}{k}}\frac{(x^{2+k}z/w;x^{2k})_{\infty}}
{(x^{-2+k}z/w;x^{2k})_{\infty}}\Psi^{\mp}(w)\Psi^{\pm}(z)\nn\\
&&\qquad\qquad\qquad =\frac{1}{x-x^{-1}}\left(\delta\left(x^k\frac{w}{z}\right)-\delta\left(x^{-k}\frac{w}{z}\right)\right).
\en
By construction, the following statement is obvious.
\begin{thm}
The following currents $x^{\pm}(z)$ and operator $d$ 
with $h$ give a level $k$ representation of $\uq$ .
\bea
&&x^+(z)= 
{\Psi}(z)\ :\EXP{-\sum_{n\neq 0}\frac{1}{[\![kn]\!]_x}a_n z^{-n}}:e^{\al}
z^{\frac{1}{k}h}
,\label{uqe2}\\
&&x^-(z)= 
{\Psi}^{\dagger}(z)\ 
:\EXP{\sum_{n\neq 0}\frac{x^{k|n|}}{[\![kn]\!]_x}a_n z^{-n}}:e^{-\al}
z^{-\frac{1}{k} h} 
\label{uqf2},\\
&& d=d^{PF}+d^{a},\label{lzerobar}
\ena
where 
\bea
&&d^{a}=
-\sum_{m>0}\frac{m^2x^{km}}{[\![2m]\!]_x [\![km]\!]_x}a_{-m}a_{m}-\frac{{h}^2}{4k}
\ena
and $d^{PF}$ is an operator such that 
\be
&&d^{PF}\cdot 1\otimes e^{\frac{\ell}{2}\al}=-\frac{\ell(k-\ell)}{2k(k+2)} 1\otimes e^{\frac{\ell}{2}\al},\\
&&[d^{PF},\Psi(z)]=-z\frac{\partial}{\partial z}\Psi(z), \quad 
[d^{PF},\Psi^\dagger(z)]=-z\frac{\partial}{\partial z}\Psi^\dagger(z).
\en
\end{thm}
We define the $\Z_{2k}$ charge of 
$\Psi_{\pm,\frac{m}{k}-n}$
 and $1\otimes e^{\frac{\ell}{2}\al}$ to be 
$\pm 2$ and  $\ell$ mod $2k$ respectively. 
For example, the  $\Z_{2k}$ charge of the vector 
\bea
&&\Psi_{\vep_1,\frac{\ell+2(\ep_2+\cdots+\ep_s)}{k}-n_1}\Psi_{\vep_2,\frac{\ell+2(\ep_3+\cdots+\ep_s)}{k}-n_2}\cdots
\Psi_{\vep_s,\frac{\ell}{k}-n_s}\otimes e^{\frac{\ell}{2}\al}
 \lb{pfvector}
\ena
 is $\ell+2\sum_{j=1}^s \vep_j$. 
Let us denote by $\H^{PF}_{\ell,M}$ the irreducible 
parafermion module of the $\Z_{2k}$ charge $M$ defined by the relation
\bea
\Omega_\ell&=&\bigoplus_{\tilde{M}\in \ell+2\Z}\H^{PF}_{\ell,\tilde{M}}\otimes 
e^{\frac{\tilde{M}}{2}\al}
%\nn\\&=&
=\bigoplus_{n\in\Z}\bigoplus_{M=0\atop {\rm mod} 2k}^{2k-1}
\H^{PF}_{\ell,M}\otimes e^{\frac{M+2kn}{2}\al}
.\lb{OmegatoPF}
\ena
Here $\H^{PF}_{\ell,M}=\{0\}$ for $M\not\equiv \ell\ \mod 2$ and 
$\H^{PF}_{\ell,M}=\H^{PF}_{\ell,M+2k}$. We also 
assume the symmetry\cite{GeQiu}
\be
&&\H^{PF}_{\ell,M}=\H^{PF}_{k-\ell,M+k}=\H^{PF}_{\ell,-M}.
\en
The basic parafermion currents act on the space $\H^{PF}_{\ell,M}$ as the following linear operators.
\be
&&\Psi(z)\ :\ \H^{PF}_{\ell,M} \to \H^{PF}_{\ell,M+2},\\
&&\Psi^\dagger(z)\ :\ \H^{PF}_{\ell,M} \to \H^{PF}_{\ell,M-2}.
\en

The character of the $q$-$\Z_k$-parafermion space $\H^{PF}_{\ell,M}$ is known to 
be\cite{GeQiu}  
\bea
&&(x^4)^{-\frac{c^{PF}}{24}}\trace_{\H^{PF}_{\ell,M}} x^{-4d^{PF}}=\eta(\bar{\tau})\ 
c^{\la_\ell}_{\la_M}(\bar{\tau})\lb{chPF}
\ena
where $c^{PF}=\frac{2(k-1)}{k+2}$. $c^{\la_\ell}_{\la_M}(\bar{\tau})$ and $\eta(\bar{\tau})$ are the string function and 
Dedekind's $\eta$-function, given by \mref{stfn} and 
\be
&&\eta(\bar{\tau})=(x^4)^{\frac{1}{24}}(x^4;x^4)_\infty.
\en

From \eqref{VlatoOmega} and \eqref{OmegatoPF}, the level $k$ irreducible highest-weight module $V(\la_\ell)$ of $U_x(\slth)$ with highest weight $\la_\ell$
is realized as follows: 
\bea
&&V(\la_\ell)=\F^a\otimes \bigoplus_{n\in\Z}\bigoplus_{M=0\atop {\rm mod} 2k}^{2k-1}\H^{PF}_{\ell,M}\otimes 
e^{{(M+2kn)}\frac{\alpha}{2}}.\lb{Fockslth}
\ena
In particular, the highest-weight vector is given by 
\bea
&&1\otimes 1\otimes e^{{\ell}\frac{\alpha}{2}}.
\ena
From \eqref{Fockslth}, 
the normalised character of $V(\la_\ell)$ is evaluated as follows:
\bea
\chi^{(k)}_{\ell}(x^{4},y)&=&(x^4)^{-\frac{c}{24}}\trace_{V(\la_\ell)}x^{-4d}y^{\frac{h}{2}}\nn\\
&=&\sum_{n\in\Z}\sum_{M=0\atop {\rm mod} 2k}^{2k-1}\ c^{\la_\ell}_{\la_{M}}(\bar{\tau}) x^{4k(n+\frac{M}{2k})^2}y^{k(n+\frac{M}{2k})}.
\ena
By setting $y=x^{-2}$, we reproduce the level $k$ principally specialised character
 \eqref{ch1}.

Now let us consider the $U_{x,p}(\slth)$-modules 
$\displaystyle{\widehat{V}(\la_\ell)=\bigoplus_{m\in \Z}\ V(\la_\ell)\otimes e^{-mQ}}.$
From \eqref{Fockslth}, we have
\bea
&&\widehat{V}(\la_\ell)=\bigoplus_{m\in \Z}\bigoplus_{n\in\Z}
\bigoplus_{M=0\atop {\rm mod}\ 2k}^{2k-1}\F_{M;m,\ell,n}
\ena
with
\bea
&&\F_{M;m,\ell,n}=\F^a\otimes \H^{PF}_{\ell,M}\otimes 
e^{{(M+2kn)}\frac{\alpha}{2}}\otimes e^{-mQ}.\lb{FockSOS}
\ena

Let $r$ be generic and  
%, $\F_{M;m,\ell,n}$ is irreducible.
note that
\bea
&&P|_{\F_{M;m,\ell,n}}=m,\qquad  
P+h|_{\F_{M;m,\ell,n}}=M+m+2kn.\lb{PPh}
\ena
From \eqref{FockSOS} and \eqref{chPF}, 
the character of the space $\F_{M;m,\ell,n}$ is evaluated as follows.
\be
&&(x^4)^{-\frac{c_{Vir}}{24}}\trace_{\F_{M;m,\ell,n}}x^{-4\widehat{d}}
= c^{\la_\ell}_{\la_{M}}(\bar{\tau})\  
x^{\frac{(mr-(M+m+2kn))r^*)^2}{krr^*}},
\en
where $c_{Vir}=\frac{3k}{k+2}\left(1-\frac{2(k+2)}{rr^*}\right)$.
This coincides with the one point function \eqref{foneptf}  
of the fusion SOS model for $a=M+m+2kn$. 
We hence make the following identification:
\be
{\rm the\ SOS\ space\ of\ states : }\quad \H_{m,a}^{(\ell)} &\longleftrightarrow& \F_{M;m,\ell,n}\qquad a=M+m+2kn,\ M\equiv \ell \mod 2,\\
{\rm the\ corner\ Hamiltonian :}\quad H^{(\ell)}_{m,a}&\longleftrightarrow& -\widehat{d}-\frac{1}{24}c_{Vir}.
\en

Furthermore let us set 
\be
&&\F_{m,\ell}(n)\equiv\bigoplus_{M=0\atop {\rm mod}\ 2k}^{2k-1}\F_{M;m,\ell,n}.
\en
When $r$ is generic, the character of $\F_{m,\ell}(n)$
%given by 
%\be
%&&\sum_{M=0\atop {\rm mod}\ 2k}^{2k-1}c^{\la_\ell}_{\la_{M}}(\bar{\tau})\  
%x^{\frac{(mr-(M+m+2kn))r^*)^2}{krr^*}}.
%\en 
%This 
coincides with the one of the irreducible 
Virasoro module $Vir_{m,a}\ (a\equiv \ell+m\ \mod 2) $ associated with the 
coset $(\slth)_k \oplus (\slth)_{r-k-2} /(\slth)_{r-2}$. 
In addition, in Appendix \ref{RSOS},
 we consider the case when $r$ is an integer $>k+2$. In this case, 
$\F_{m,\ell}(n)$ is reducible.  We observe that the BRST resolution 
of the complex formed by $\F_{m,\ell}(n)$ yields  the irreducible 
coset Virasoro {\it minimal} module
$Vir_{m,a},\ (a\equiv m+\ell \ \mod 2) $. 
These considerations leads us to the following conjecture:
\begin{conj}
The space $\F_{m,\ell}(n)$  is isomorphic to the irreducible 
coset Virasoro module $Vir_{m,a},$\ $a\equiv m+\ell \ \mod 2 $ 
with the central charge $c_{Vir}=\frac{3k}{k+2}\left(1-\frac{2(k+2)}{4rr^*}\right)$ 
and the highest weight $h_{m,a}=\frac{\ell(k-\ell)}{2k(k+2)}+\frac{(mr-ar^*)^2-k^2}{4krr^*}$.
\end{conj}
 
\subsubsection{The vertex operators}
The vertex operator $\widehat{\Phi}_\ep(u)\ (\ep=-k,-k+2,..,k)$  
of the elliptic algebra $U_{x,p}(\slth)$ in \eqref{fvo} acts 
on the space  $\F_{M;m,\ell,n}$ as
\be
&&\widehat{\Phi}_\ep(u)\ :\ \F_{M;m,\ell,n} \to \F_{M+\vep;m,k-\ell,n}
\en
and satisfies the commutation relation \eqref{eqn:VO1}.
The relation \eqref{eqn:VO1} is similar to 
that of the lattice vertex operators \eqref{commrel} but not precisely 
the same. Noting the symmetry \mref{gaugek}, it turns out that the following 
 gauge transformation $\widehat{\Phi}(u) \to \Phi_\ep(u)$  resolves
  this discrepancy.
\bea
\hspace*{-10mm}\Phi_\ep(u)=\widehat{\Phi}_\ep(u)\frac{1}{(P+h,P+h+\ep)_k}=\Phi_{k,k}(u)
\widehat{F}^+\left(u+\frac{k}{2}-r\right)^{\frac{k-\vep}{2}}\frac{1}{(P+h,P+h+\ep)_k}.
%\nn \\
\lb{fvoSOS}
\ee
In fact $\Phi_{\ep}(u)$ satisfies 
the commutation relation 
\be
&&\Phi_{\ep_2}(u_2)
\Phi_{\ep_1}(u_1)=\\
&&\sum_{\nu_1+\nu_2=\ep_1+\ep_2}
W^{(k,k)}\left(\left.
\begin{array}{cc}
P+h
&P+h-\nu_1\\
P+h-\ep_2&P+h-\ep_1-\ep_2
\end{array}
\right|u_1-u_2\right)
\Phi_{\nu_1}(u_1)
\Phi_{\nu_2}(u_2).\nonumber\\
\en
This is exactly the same commutation relation as \eqref{commrel} if 
we make the identification
\bea
&&\Phi_{a+\ep,a}(u)={\Phi}_\ep(u) \lb{relPhi}
\ena
on $\F_{M;m,\ell,n}= \H^{(\ell)}_{m,a}\ (a=M+m+2kn,\ M\equiv \ell \mod\ 2)$. 
This is the realization of the half transfer matrix in terms of the vertex operator 
of $U_{x,p}(\slth)$.

As discussed in \cite{LaP98}, there is a second realization of the vertex operators. 
This is due to 
the symmetries of the space of states $\H^{(\ell)}_{m,a}=\H^{(\ell)}_{-m,-a}$ 
and the 
Boltzmann weights
\begin{eqnarray*}
W^{(k,k)}\left(\left.
\begin{array}{cc}
a&b\\
c&d
\end{array}
\right|u\right)=
W^{(k,k)}\left(\left.
\begin{array}{cc}
-a&-b\\
-c&-d
\end{array}
\right|u\right).\lb{BWminus}
\end{eqnarray*}
In fact from \eqref{commrel}, we have 
\bea
&&\Phi_{c,a}(u_2)\Phi_{a,b}(u_1)=\sum_{d}W^{(k,k)}\BW{-c}{-d}{-a}{-b}{u_1-u_2}
\Phi_{c,d}(u_1)\Phi_{d,b}(u_2).\lb{phiphim}
\ena
On the other hand, we have an operator 
\bea
&&\Phi_{-a',-a}(u)\ :\ \H^{(\ell)}_{-m,-a} \to \H^{(k-\ell)}_{-m,-a'},
\ena 
which can be shown, using the same argument that leads to
\eqref{commrel}, 
to satisfy the commutation relation   
\bea
\Phi_{-c,-a}(u_2)\Phi_{-a,-b}(u_1)=\sum_{d}W^{(k,k)}\BW{-c}{-d}{-a}{-b}{u_1-u_2}
\Phi_{-c,-d}(u_1)\Phi_{-d,-b}(u_2).
\ena
Comparing this with \eqref{phiphim},  
we can simply make the following identification.
\bea
&&\Phi_{a,b}(u)=\Phi_{-a,-b}(u)
\ena
on $\H^{(\ell)}_{-m,-b}$. From \eqref{relPhi}, we have
\bea
\hspace*{-10mm}&&\Phi_{\ep}(u)
=
\Phi_{-\ep}(u)\biggl|_{\H^{(\ell)}_{-m,-a}}
%\nn\\&=&
=\Phi_{k,k}(u)\widehat{F}^+\left(u+\frac{k}{2}-r\right)^{\frac{k+\ep}{2}}\frac{1}{(P+h,P+h-\ep)_k}\biggl|_{\H^{(\ell)}_{-m,-a}}.
\lb{secondPhi}
\ena

\subsubsection{The tail operator }\lb{tailoprel}
We next consider the realisation of the tail operator introduced in Section 3.5. 
The tail operator ${\Lambda}_{a',a}(u) : \H_{m,a}^{(\ell)}\to  \H_{m,a'}^{(\ell)}$  
is characterised by the commutation relation \eqref{TailPhi}. 
In addition, it follows from formula \eqref{gencorr} that we only have to consider the case 
$a'-a\in 2\Z$. In a similar way to the case of vertex operators, 
we seek to realise the tail operator  
in the following form 
\bea
{\Lambda}_{a+\ep,a}(u)={\Lambda}_{\ep}(u)\ (\ep\in 2\Z)\lb{reltail}
\ena
on the space $\F_{M;m,\ell,n}= \H^{(\ell)}_{m,a}\ (a=M+m+2kn, \ M\equiv \ell \mod 2)$.  
Note that 
from \eqref{PPh}, the tail operator should satisfy 
\be
&&[P,{\Lambda}_{\ep}(u)]=0,\quad [P+h,{\Lambda}_{\ep}(u)]=\ep.
\en
Substituting \eqref{relPhi} and \eqref{reltail} into \eqref{TailPhi}, 
we obtain the following commutation relation.
\begin{eqnarray}
&&\hspace*{-20mm}{\Lambda}_{\ep_1}(u_1)
\Phi_{\ep_2}(u_2)
%\nn\\ \hspace*{-10mm}&&
=\hspace*{-8mm}
\sum_{\nu_2\in\{-k,-k+2,..,k\}\atop 
\nu_1+\nu_2=\ep_1+\ep_2}
L^{(k)}\left(\left.\hspace*{-2mm}
\begin{array}{cc}
P+h&P+h-\nu_2\\
P+h-\ep_1&P+h-\ep_1-\ep_2
\end{array}\right|u_2-u_1\right)
\Phi_{\nu_2}(u_2) {\Lambda}_{\nu_1}(u_1).\lb{TailPhimu}
\end{eqnarray}
Here $L^{(k)}
%\left(
%\left.
%\begin{array}{cc}
%a&b\\
%c&d
%\end{array}\right|u
%\right)
$ are
the $k$-fusion $L$-matrices, explicit formulae for which are given 
in Appendix \ref{Lfusion}.

Let us first consider the case $\ep<0$ in ${\Lambda}_{\ep}(u)$.
Setting $\ep_1=-2k, 
\ep_2=k, \nu_1=-2s, \nu_2=-k+2s\ (s\in {\mathbb{N}})$ in \eqref{TailPhimu}, we have  
\begin{eqnarray}
&&\hspace*{-20mm}{\Lambda}_{-2k}(u_1) \Phi_k(u_2)
%\nn\\&&\hspace*{-14mm}
=\sum_{s=0}^k
L^{(k)}\left(\left.
\begin{array}{cc}
P+h&P+h+k-2s\\
P+h+2k&P+h+k
\end{array}
\right|u_2-u_1\right)\Phi_{-k+2s}(u_2)
{\Lambda}_{-2s}(u_1),
\label{eqn:tailfull}
\end{eqnarray}
where the coefficient $L^{(k)}$ is given by \eqref{Lmatk2} with $m=P+h, n=P+h+2k$.
%\begin{eqnarray*}
%&&L^{(k)}\left(\left.\begin{array}{cc}P+h&
%P+h+k-2s\\
%P+h+2k&P+h+k
%\end{array}\right|u
%\right)\nn\\[2mm]
%&&=\frac{[P+h+k]_k[k]_{k-s}[-u+s]_{s}[-u+P+h+k-s]_{k-s}}
%{[P+h+2k]_k[-u]_k}.
%\end{eqnarray*}
This $L^{(k)}$ has simple poles at $u=0,-1,\cdots,-k+1$.
We take the residue at $u_1=u_2+k-1$.
Assuming that
${\Lambda}_{-2k}(u_1)\Phi_k(u_2)$ on the left hand side
of \mref{eqn:tailfull} does not have a pole at
$u_1=u_2+k-1$,
we have the necessary condition
\begin{eqnarray}
0&=&\sum_{s=0}^k
[P+h+k]_s[P+h+2k-s-1]_{k-s}[k+s-1]_s[k]_{k-s}\nonumber\\
&\times&\Phi_{-k+2s}(u_2)
\Lambda_{-2s}(u_2+k-1).\label{eqn:necessary}
\end{eqnarray}
Substituting in formula \eqref{fvoSOS}, we obtain the recursion relation for the 
tail operator   
\be
0&=&\sum_{s=0}^k
[P+h+2k]_s[P+h+3k-s-1]_{k-s}[k+s-1]_s[k]_{k-s}\nonumber\\
&\times&F^+\left(u+\frac{k}{2}-r\right)^{k-s}
\Lambda_{-2s}(u_2+k-1)\frac{1}{(P+h-2s,P+h-k)_k}.
\label{eqn:necessary2}
\en
To solve this for $\Lambda_{\ep}(u+k-1)$, we make the following ansatz. 
\begin{eqnarray}
{\Lambda}_{-2s}(u+k-1)&=&
\widehat{F}^+\left(u+\frac{k}{2}-r\right)^{s}
\lambda_s(P+h),
\lb{tail1}
\end{eqnarray}
where $\lambda_s(P+h)$ is a function to be determined.
Then the necessary condition
(\ref{eqn:necessary2}) reduces to the relation
\begin{eqnarray*}
&&\sum_{s=0}^k
[P+h]_s[P+h+k-s-1]_{k-s}
[k+s-1]_{s}[k]_{k-s}\nonumber
\frac{\lambda_s(P+h)}{
(P+h-2s,P+h-k)_k}=0.
\end{eqnarray*}
This equation is satisfied if we choose 
$\lambda_s(P+h)=g_{P+h}^{-1} \cdot g_{P+h-2s}$ (where $g_a$
is defined below \mref{crossing}).
We hence obtain the identification 
\be
&&{\Lambda}_{-2s}(u+k-1)=
g_{P+h}\widehat{F}^+\left(u+\frac{k}{2}-r\right)^{s}g_{P+h}^{-1}.
\lb{tail1fin}
\en
To check that this is also a sufficient condition, 
we substitute \mref{tail1} back into the 
full commutation relation \mref{eqn:tailfull}. We find that
( \ref{eqn:tailfull}) then reduces to the same theta function
identity that occurs in 
the commutation relation \mref{eqn:VO1} for vertex operators, which we have checked 
for the cases $k=1,2,3$. See Appendix \ref{commrelapp} for details.

Let us next study ${\Lambda}_{\ep}(u)$ for the case $\ep>0$.
For this purpose, we use the second realization of 
the vertex operators \eqref{secondPhi}. 
In addition we have the symmetry 
\be
&&\psi^{(k)}(u)^{-a}_{-b}=\psi^{(k)}(u)^{a}_{b},\quad 
\psi^{*(k)}(u)^{-a}_{-b}=\psi^{*(k)}(u)^{a}_{b}.
\en
Hence 
\be
&&L^{(k)}\BW{-a}{-b}{-c}{-d}{u} =L^{(k)}\BW{a}{b}{c}{d}{u}.
\en
Therefore from \eqref{TailPhi}, we have
\bea
&&{\Lambda}_{c,a}(u_1)
\Phi_{a,b}(u_2)
%\nn\\
=
\sum_{d}
L^{(k)}\left(\left.
\begin{array}{cc}
-c&-d\\
-a&-b
\end{array}\right|u_2-u_1\right)
\Phi_{c,d}(u_2) {\Lambda}_{d,b}(u_1).\lb{TailPhiplus}
\ee
On the other hand, consider the operator 
\be
&&{\Lambda}_{-a',-a}(u):\H^{(\ell)}_{-m,-a} \to \H^{(\ell)}_{-m,-a'}.
\en
From the same argument that leads to \eqref{TailPhi}, we have the commutation relation
\bea
&&{\Lambda}_{-c,-a}(u_1)
\Phi_{-a,-b}(u_2)
%\nn\\
=
\sum_{d}
L^{(k)}\left(\left.
\begin{array}{cc}
-c&-d\\
-a&-b
\end{array}\right|u_2-u_1\right)
\Phi_{-c,-d}(u_2) {\Lambda}_{-d,-b}(u_1).\lb{TailPhiplus2}
\ee
Comparing \eqref{TailPhiplus} and \eqref{TailPhiplus2} and using the second 
realisation of the vertex operator \eqref{secondPhi}, we identify 
\be
&&{\Lambda}_{-a,-b}(u)={\Lambda}_{a,b}(u)
\en
under the identification of $\H^{(\ell)}_{-m,-a}$ with $\H^{(\ell)}_{m,a}$. 
Therefore, we obtain the following realisation
\be 
{\Lambda}_{a+2s,a}(u+k-1)&=&{\Lambda}_{-a-2s,-a}(u+k-1)\nn\\
&=&g_{P+h}\widehat{F}^+\left(u+\frac{k}{2}-r\right)^{s}g_{P+h}^{-1}
\biggl|_{\H^{(\ell)}_{-m,-a}}.
\lb{tail1plusfin}
\en
We are thus lead to the following simple conjecture.

\begin{conj}
\label{conj2}~~{\bf (The Realisation of the Tail Operator)}~~
The tail operator ${\Lambda}_{a\pm 2s,a}(u)~(s \in {\mathbb{N}})$ 
is realized by the following power of the half-current:
\begin{eqnarray*}
{\Lambda}_{a\pm 2s,a} (u)&=&
g_{P+h} \widehat{F}^+\left(u-\frac{k}{2}-r+1\right)^{s}
g_{P+h}^{-1}\biggl|_{\H^{(\ell)}_{\mp m,\mp a}}
\\
&=&
\oint_{J_1}\frac{dw_1}{2\pi i w_1} \cdots 
\oint_{J_{s}}\frac{dw_{s}}{2\pi i w_{s}} 
F(v_1) \cdots F(v_{s})\nonumber\\
&\times&(-1)^s \sqrt{\frac{[P+h-2s]}{[P+h]}}\prod_{j=1}^s
\frac{[u-v_j+P+h-\frac{k}{2}-2s+2j]}{
[u-v_j-\frac{k}{2}+1]}
\biggl|_{\H^{(\ell)}_{\mp m,\mp a}}.\nonumber
\end{eqnarray*}
The integrations contours $J_j~(1\leqq j \leqq s)$
are given by
\begin{eqnarray*}
J_j:~|x^{-k}z^{-1}|, |x^{-2}w_{j+1}|<|w_j|<|p^{-1}x^{-k}z^{-1}|.
\end{eqnarray*}
\end{conj}

\section{Summary}
In this paper, we have first generalised the approach of \cite{LaP98}
in order to obtain the trace formula \mref{gencorr} for N-point
correlation functions of the level $k$ fusion analogue of the the
eight-vertex model. The objects that appear in this trace are the
space of states $\cH^{(\ell)}_{m,a}$, the corner Hamiltonian
$H^{(\ell)}_{m,a}$, the half transfer matrices $\Phi_{b,a}(u)$ and 
$\Phi^*_{b,a}(u)$, and the tail operator $\gL_{b,a}(u)$. We have 
constructed each of these objects in terms of the algebra
$U_{x,p}(\slth)$ in Section 5. A multiple
integral formula for \mref{gencorr} then follows rather simply.

In a following paper \cite{KKW04b}, we shall 
examine the $k=2$ case in detail.
We shall make use of the rather simpler 1-boson/1-fermion free field 
realisation that exists in this case in order to produce and 
analyse explicit expressions for certain correlation functions.

\subsection*{Acknowledgements}
The authors would like to thank M. Jimbo for stimulating discussions.
They also thank A. Kuniba, M. Lashkevich, T. Nakanishi, A. Nakayashiki, M. Okado,  
Y. Pugai, Y.H. Quano, M. Rossi, J. Shiraishi, and T. Takebe  for
useful conversations. 
TK and HK are grateful to colleagues at Heriot-Watt University for 
their kind hospitality during the period when this work was started.
HK and RW respectively thank P. Goddard,  
in  DAMTP, Cambridge Univ., and M. Jimbo at Tokyo University for
their warm hospitality.

TK is supported by the Grant-in-Aid for Young Scientists (B)
(14740107) from the JSPS. HK thanks the JSPS and Royal Society for an
exchange fellowship, and acknowledges support
from the Grant-in-Aid for Scientific Research (C) 15540033, JSPS.
RAW acknowledges partial support given by the EUCLID research
training network funded by the European 
Commission under contract HPRN-CT-2002-00325.        

\newpage
\appendix
\setcounter{equation}{0}
\begin{appendix}
\section{The Proof of Formula (3.19)}\lb{ChPartition}
It is convenient to use the parametrisation 
$a=\ell+m+2(kn+t)\ (n\in \Z,\ 0\leq t\leq k-1 \mod\ 2)$.
From \eqref{foneptf}, we have 
\be
\sum_{a\in m+\ell+2\Z}[a] \, \Tr_{\H^{(\ell)}_{m,a}} x^{4H^{(\ell)}_{m,a}}
&=&\sum_{n\in\Z}\sum_{t=0\atop {\rm mod} k}^{k-1}[\ell+m+2(kn+t)]
c^{\la_\ell}_{\la_{M}}(\bar{\tau})
x^{\frac{(mr-(\ell+m+2(kn+t))r^*)^2}{krr^*}}\nn\\
&=&\sum_{n\in\Z}\sum_{M=0\atop {\rm mod} 2k}^{2k-1}[M+m+2kn]c^{\la_\ell}_{\la_{M}}(\bar{\tau})
x^{\frac{(mr-(M+m+2kn))r^*)^2}{krr^*}}
\en
with $M\equiv a-m \mod 2$. Note that 
$c^{\la_\ell}_{\la_{M}}(\bar{\tau})=0$ for  $M\not\equiv \ell$ mod 2.
Using the formula $[u]=x^{\frac{u^2}{r}-u}\sum_{s\in\Z}(-)^s x^{rs^2}x^{s(2u-r)}$, 
this can be rewritten as follows:
\be
&&x^{\frac{m^2}{r^*}-m}\sum_{s\in\Z}(-)^s x^{r^*s(s+1)}x^{-2ms}\ I(s),\nn\\
\en
where we define
\be
&&I(s)=\sum_{n\in\Z}\sum_{M=0\atop {\rm mod} 2k}^{2k-1}c^{\la_\ell}_{\la_{M}}(\bar{\tau})
x^{\frac{1}{k}{(k(2n-s)+M+\frac{k}{2})^2}}x^{-\frac{k^2}{4}}.
\en
We show that $I(s)\ (s\in \Z)$ is independent of $s$. 
In fact, for
the case $s=2u+1\ (u\in \Z)$, we can eliminate $u$ by shifting $n \to n+u$. 
Then we obtain
\be
I(2u+1)
&=&\sum_{n\in\Z}\sum_{M=0\atop {\rm mod} 2k}^{2k-1}c^{\la_\ell}_{\la_{M}}(\bar{\tau})
x^{{4k(n+\frac{M}{2k})^2}}x^{-2k(n+\frac{M}{2k})}\nn\\
&=&\chi^{(k)}_\ell(\bar{\tau}).\lb{evenI}
\en
Let us next set $s=2u$. We have 
\be
I(2u)&=&\sum_{n\in\Z}\sum_{M=0\atop {\rm mod} 2k}^{2k-1}
c^{\la_\ell}_{\la_{M}}(\bar{\tau})
x^{{4k(n+\frac{M}{2k})^2}}x^{2k(n+\frac{M}{2k})}.\nn\\
\en
By changing $n\to -n,\ M\to -M$ and using the symmetry 
$c^{\la_\ell}_{\la_{M}}(\bar{\tau})=c^{\la_\ell}_{\la_{-M}}(\bar{\tau})$, we find that 
$I(2u)$ coincides with $I(2u+1)$. Therefore the LHS of \mref{charpart}
is 
\be
{\rm LHS}&=&x^{\frac{m^2}{r^*}-m}\sum_{s\in\Z}(-)^s x^{r^*s(s+1)}x^{-2ms}\ 
\chi^{(k)}_\ell(\bar{\tau})\nn\\
&=&[m]^*\chi^{(k)}_\ell(\bar{\tau}).
\en
This coincides with the RHS.

\section{The BRST Resolution of $\F_{m,\ell}(n)$}\lb{RSOS}
Let $r>k+2 \in\Z$ and fix $m,\ell \in \Z$ with $0\leq \ell \leq k$. 
Note that regarding the half currents $E^+(u)$ and $F^+(u)$ as 
the screening currents\cite{Ko98}, 
we can define the $q$-Virasoro algebra associated with 
the coset $(\slth)_k\oplus (\slth)_{r-k-2}/(\slth)_{r-2}$ as their
commutant \cite{FrRe}. As such a $q$-Virasoro module, $\F_{m,\ell}(n)$ is reducible. 
Consider the BRST operator given by
\be
&&Q^+_s=\widehat{E}(u)^s \ :\   \F_{m,\ell}(n)\ \to \F_{m-2s,\ell}(n).
\en
\begin{prop}\cite{Ko98}
The BRST operator $Q^+_s$ is independent of $u$ and 
is  nilpotent in the following sense.
\be
 &&Q^+_sQ^+_{r^*-s}=Q^+_{r^*-s}Q^+_s=0.
\en
\end{prop}
Setting $Q_{2j}=Q^+_m,\ Q_{2j+1}=Q^+_{r^*-m}\ (j\in \Z)$, we then have the following 
 complex ${\cal C}_{m,\ell}$ 
\be
&&\cdots {\xrightarrow{Q_{-2}}}\ \F_{-m+2r^*,\ell}(n)\ {\xrightarrow{Q_{-1}}}\ 
\F_{m,\ell}(n)\ 
{\xrightarrow{Q_{0}}}\ \F_{-m,\ell}(n)\ {\xrightarrow{Q_{1}}}\  \cdots .
\en
We conjecture the following statement about the cohomology group
$H^j({\cal C}_{m,\ell})$ of this complex \cite{Ko98}:
\be
&&H^j({\cal C}_{m,\ell})=0 \qquad ( j\not=0 ).
\en
Then by the Euler-Poincar${\rm \acute{e}}$ principle, we can evaluate 
the character of the 0-th cohomology group as follows:
\bea
&&
\Tr_{H^0({\cal C}_{m,\ell})}\ x^{-4\left(\widehat{d}+\frac{c_{Vir}}{24}\right)} =\sum_{j\in \Z}\left( 
\Tr_{\F_{m-2r^*j,\ell}(n)} \ x^{-4\left(\widehat{d}+\frac{c_{Vir}}{24}\right)}-
\Tr_{\F_{-m-2r^*j,\ell}(n)} \ x^{-4\left(\widehat{d}+\frac{c_{Vir}}{24}\right)}\right),
\nn\\&&\lb{chcosetVir}
\ena
where from \eqref{FockSOS}, we have 
\be
&&\Tr_{\F_{\pm m-2r^*j,\ell}(n)} \ 
x^{-4\left(\widehat{d}+\frac{c_{Vir}}{24}\right)}=
\sum_{M=0\atop {\rm mod} 2k }^{2k-1} c^{\la_\ell}_{\la_M}(\bar{\tau})\ 
x^{\frac{(\pm mr-(M\pm m-2r^*j+2kn))r^*-2rr^*j)^2}{krr^*}}.
\en
\eqref{chcosetVir} coincides with the branching function $b_{m,a}^{(\ell)}(\bar{\tau})$
 ,i.e., the 
 character of the irreducible coset Virasoro minimal module $Vir_{m,a}\ 
(a=m+\ell \ \mod 2)$. This is also known to equal to 
the one point function of the $k$ fusion RSOS model with  height restriction  $1\leq a\leq r-1$ and $1\leq m\leq r-k-1$\cite{DJMO86b,DJKMO87}. 
Hence in this case we can make the identification
\be  
&&H^0({\cal C}_{m,\ell})\ \longleftrightarrow\ \cH_{m,a}^{RSOS(\ell)}. 
\en

\section{Fusion of the $L$-matrix}\lb{Lfusion}

In this appendix, we give
explicit formulae for the fused $L$-matrix.
The $L$-matrix $L^{(1)}$ is defined by
\begin{eqnarray*}
L^{(1)}\left(\left.\begin{array}{cc}
a&b\\
c&d\end{array}\right|u\right)&=&\sum_{\vep}
\psi_\vep^{*(1)}(u)_{c}^{d}
\psi_\vep^{(1)}(u)_{b}^{a}.
\end{eqnarray*}

From \eqref{intertwinvec} and \eqref{dualintvec}, we obtain the formula
\begin{eqnarray*}
L^{(1)}\left(\left.\begin{array}{cc}
m&m\pm1\\
n&n\pm1\end{array}\right|u\right)&=&
\frac{[u\pm\frac{n-m}{2}][\frac{n+m}{2}]}{[u][n]},\\
L^{(1)}\left(\left.\begin{array}{cc}
m&m\mp1\\
n&n\pm1\end{array}\right|u\right)&=&
\frac{[u\pm\frac{n+m}{2}][\frac{n-m}{2}]}{[u][n]}.
\end{eqnarray*}
The  $k$-fused $L$-matrix $L^{(k)}$ is given by
\begin{eqnarray*}
L^{(k)}\left(\left.\begin{array}{cc}
m_0&m_k\\
n_0&n_k
\end{array}\right|u\right)
&=&\sum_{\epsilon}
\psi_\epsilon^{*(k)}(u)_{n_0}^{n_k}
\psi_\epsilon^{(k)}(u)_{m_k}^{m_0}.
\end{eqnarray*}
According to the fusion formulae for $\psi^{(k)}(u)$ and $\psi^{*(k)}(u)$ \eqref{kpsi}
 and \eqref{kpsis}, the $L^{(k)}$
satisfies the following fusion formula.
\begin{eqnarray}
L^{(k)}\left(\left.\begin{array}{cc}
m_0&m_k\\
n_0&n_k
\end{array}\right|u\right)
&=&\sum_{n_1,n_2,\cdots,n_{k-1}}
L^{(1)}\left(\left.\begin{array}{cc}
m_0&m_1\\
n_0&n_1
\end{array}\right|u+k-1\right)
L^{(1)}\left(\left.\begin{array}{cc}
m_1&m_2\\
n_1&n_2
\end{array}\right|u\right)\nonumber\\
&\cdots&
L^{(1)}\left(\left.\begin{array}{cc}
m_{k-1}&m_k\\
n_{k-1}&n_k
\end{array}\right|u\right),
\lb{Lfusform}
\end{eqnarray}
where the right hand side is independent of
the dynamical variables $m_1,m_2,\cdots,m_{k-1}$.
By induction making use of \mref{Lfusform}, 
we obtain the following compact expressions. 
\begin{eqnarray*}
L^{(k)}\left(\left.\begin{array}{cc}
m&m-k+2i\\
n&n-k+2j
\end{array}\right|u\right)
&=&\sum_{l=Max(0,i-j)}^{Min(i,k-j)}
L^{(k-i)}\left(\left.\begin{array}{cc}
m&m-k+i\\
n&n-k+2j-i+2l
\end{array}\right|u+i \right)\nonumber\\
&&\qquad\qquad \times L^{(i)}\left(\left.\begin{array}{cc}
m-k+i&m-k+2i\\
n-k+2j-i+2l&n-k+2j
\end{array}\right|u\right)\nonumber\\
\\
&=&\sum_{l=Max(0,i+j-k)}^{Min(i,j)}
L^{(i)}\left(\left.\begin{array}{cc}
m&m+i\\
n&n-i+2l
\end{array}\right|u+k-i \right)\nonumber\\
&&\qquad\qquad \times L^{(k-i)}\left(\left.\begin{array}{cc}
m+i&m-k+2i\\
n-i+2l&n-k+2j
\end{array}\right|u\right).
\end{eqnarray*}
Here  $0\leq i,j\leq k$ and  we have 
\begin{eqnarray}
&&L^{(k)}\left(\left.\begin{array}{cc}
m&m+k\\
n&n+k-2j
\end{array}\right|u\right)\nonumber
\\
&=&
\frac{
\left[\begin{array}{c}
\frac{1}{2}(n+m)+k-1-j\\
k-j
\end{array}
\right]
\left[\begin{array}{c}
\frac{1}{2}(n-m)\\j
\end{array}
\right]
\left[\begin{array}{c}
-u+\frac{1}{2}(n+m)\\j
\end{array}
\right]
\left[\begin{array}{c}
-u+\frac{1}{2}(m-n)\\k-j
\end{array}
\right]
}{
\left[\begin{array}{c}
n+k-1-2j\\k-j
\end{array}
\right]
\left[\begin{array}{c}
n+k-j\\j
\end{array}
\right]
\left[\begin{array}{c}
-u\\k
\end{array}
\right]},\nonumber\\
\label{L1}
\end{eqnarray}
%and
\begin{eqnarray}
&&L^{(k)}\left(\left.\begin{array}{cc}
m&m-k\\
n&n+k-2j
\end{array}\right|u\right)\nonumber
\\
&=&
\frac{
\left[\begin{array}{c}
\frac{1}{2}(n-m)+k-1-j\\
k-j
\end{array}
\right]
\left[\begin{array}{c}
\frac{1}{2}(n+m)\\j
\end{array}
\right]
\left[\begin{array}{c}
-u+\frac{1}{2}(n-m)\\j
\end{array}
\right]
\left[\begin{array}{c}
-u-\frac{1}{2}(m+n)\\k-j
\end{array}
\right]
}{
\left[\begin{array}{c}
n+k-1-2j\\k-j
\end{array}
\right]
\left[\begin{array}{c}
n+k-j\\j
\end{array}
\right]
\left[\begin{array}{c}
-u\\k
\end{array}
\right]}.\nonumber\\
\label{L2}
\end{eqnarray}
In particular, we use the following formulae in Section \ref{tailoprel}.
\begin{eqnarray}
\hspace*{-10mm}L^{(k)}\left(\left.\begin{array}{cc}
m&m+k-2j\\
n&n+k
\end{array}\right|u\right)
&=&\frac{[\frac{n+m}{2}+k-1+j]_{k-j}
[\frac{n-m}{2}-1+j]_j}{
[n+k-1]_k}\nn\\
&\times&
\frac{[-u+\frac{m-n}{2}-j]_{k-j}
[-u-\frac{m+n}{2}+j-k]_j}{
[-u]_k},\\
\hspace*{-10mm}L^{(k)}\left(\left.\begin{array}{cc}
m&m+k-2j\\
n&n-k
\end{array}\right|u\right)
&=&\frac{[\frac{n+m}{2}]_{j}
[\frac{n-m}{2}]_{k-j}[-u+\frac{m+n}{2}-j]_{k-j}
[-u+\frac{n-m}{2}+j-k]_j}{
[n]_k[-u]_k}.\lb{Lmatk2}
\end{eqnarray}

Recently one of the authors obtained an explicit expression of $L^{(k)}$ in terms of the 
very-well-poised elliptic hyper geometric series\cite{Konno05}.

\section{Commutation Relations of the Tail  and Vertex Operators}\lb{commrelapp}

In this appendix, we 
check the commutation relations between
the tail operator 
${\Lambda}_{\vep}(u)$ and
the vertex operator ${\Phi}_{\vep}(u)$.
Let us consider an integral of the form
\begin{eqnarray*}
\oint \frac{dw_1}{2\pi i w_1}
\oint \frac{dw_2}{2\pi i w_2}
F(v_1)F(v_2)f(v_1,v_2),
\end{eqnarray*}
where the integration contours for 
$w_1$ and $w_2$ are the same.
The commutation relation
$F(v_1)F(v_2)=\frac{[v_1-v_2-1]}{[v_1-v_2+1]}F(v_2)F(v_2)$,
implies that this integral is equal to
\begin{eqnarray*}
\oint \frac{dw_1}{2\pi i w_1}
\oint \frac{dw_2}{2\pi i w_2}
F(v_1)F(v_2)f(v_2,v_1)\frac{[v_2-v_1-1]}{[v_2-v_1+1]}.
\end{eqnarray*}
Observing this, we define the notion of `weak equality'.
The functions $f(v_1,v_2)$ and $g(v_1,v_2)$
are equal in weak sense if
\begin{eqnarray*}
f(v_1,v_2)+\frac{[v_2-v_1-1]}{[v_2-v_1+1]}f(v_2,v_1)=
g(v_1,v_2)+\frac{[v_2-v_1-1]}{[v_2-v_1+1]}g(v_2,v_1).
\end{eqnarray*}
We write
$f(v_1,v_2) \sim g(v_1,v_2)$ 
to denote weak equality.

Let us consider the commutation relation
\begin{eqnarray}
&&\hspace*{-20mm}{\Lambda}_{-2k}(u_2)
\Phi_k(u_1)
%\nonumber\\&=&
=\sum_{s=0}^k
L^{(k)}\left(\left.\begin{array}{cc}
P+h&P+h+k-2s\\
P+h+2k&P+h+k
\end{array}\right|
u_1-u_2\right)
\Phi_{-k+2s}(u_1)
{\Lambda}_{-2s}(u_2).
%\nonumber\\
\label{com5}
\end{eqnarray}
This reduces to the following weak equality
\begin{eqnarray*}
I_{{\Lambda}}(v_1, v_2, \cdots,v_k)\sim 0,
\end{eqnarray*}
where
\begin{eqnarray*}
&&I_{{\Lambda}}(v_1,v_2,\cdots,v_k)\nonumber\\
&=&(-1)^k\sqrt{\frac{[n]}{[n+2k]}}
\frac{1}{(n+k,n+2k)_k}
\prod_{j=1}^k
\frac{[-u_2-v_j+n-\frac{k}{2}+2j-1]}{
[-u_2-v_j-\frac{k}{2}]}\nonumber\\
&-&\sum_{s=0}^k
L^{(k)}\left(\left.\begin{array}{cc}
n&n+k-2s\\
n+2k&n+k
\end{array}\right|
u_1+u_2+1\right)(-1)^s \sqrt{\frac{[n+k-2s]}{[n+k]}}
\frac{1}{(n,n+k-2s)_k}\nonumber
\\
&\times&
\prod_{j=1}^{k-s}
\frac{[u_1-v_j+n-\frac{k}{2}+2j-1]}{[u_1-v_j+\frac{k}{2}]}
\prod_{j=k-s+1}^k
\frac{[-u_2-v_j+n-\frac{3k}{2}+2j-1]}{
[-u_2-v_j-\frac{k}{2}]}.\nonumber\\
\end{eqnarray*}
Let us consider
\begin{eqnarray}
\widehat{\Phi}_{-k}(u_2)
\widehat{\Phi}_k(u_1)=\sum_{s=0}^k
W^{(k,k)}\left(\left.\begin{array}{cc}
P+h&P+h+k-2s\\
P+h+k&P+h
\end{array}\right|
u_1-u_2\right)
\widehat{\Phi}_{-k+2s}(u_1)
\widehat{\Phi}_{k-2s}(u_2).\nonumber\\
\label{com6}
\end{eqnarray}
This in turn reduces to the weak equality
\begin{eqnarray*}
I_{\hat{\Phi}}
(v_1,v_2,\cdots,v_k)\sim 0,
\end{eqnarray*}
where \begin{eqnarray*}
\hspace*{-10mm}I_{\hat{\Phi}}
(v_1,v_2,\cdots,v_k)
%\nonumber\\
&=&
\prod_{j=1}^k
\frac{[u_1-v_j-\frac{k}{2}]}{[u_1-v_j+\frac{k}{2}]}
\prod_{j=1}^k
\frac{[u_2-v_j+n-\frac{k}{2}+2j-1]}
{[u_2-v_j+\frac{k}{2}]}\nonumber\\
&-&\sum_{s=0}^k
\bW_{k,k}\left(\left.\begin{array}{cc}
n&n+k\\
n+k-2s&n
\end{array}\right|
u_1-u_2\right)\prod_{j=1}^{k-s}
\frac{[u_2-v_j-\frac{k}{2}]}
{[u_2-v_j+\frac{k}{2}]}
\\
&\times&
\prod_{j=1}^{k-s}
\frac{[u_1-v_j+n-\frac{k}{2}+2j-1]}
{[u_1-v_j+\frac{k}{2}]}
\prod_{j=k-s+1}^k
\frac{[u_2-v_j+n-\frac{3k}{2}+2j-1]}
{[u_2-v_j+\frac{k}{2}]}.\nonumber
\end{eqnarray*}
We have checked 
$I_{\widehat{\Phi}}
(v_1,v_2,\cdots,v_k)\sim 0$
for the case $k=1,2,3$.

\begin{prop}~~~
The $L^{(k)}$-matrix and the Boltzmann Weight $W_{k,k}$ 
are related by
\begin{eqnarray*}
&&\bW_{k,k}
\left(\left.\begin{array}{cc}
n&n+k\\
n+k-2s&n
\end{array}\right|
u\right)\nonumber\\
&=&
\sqrt{\frac{[n-k][n+2k]}{[n+k][n]}}
\frac{(n+2k,n+k)_k
}{(n,n+k-2s)_k}
L^{(k)}\left(\left.\begin{array}{cc}
n&n+k-2s\\
n+2k&n+k
\end{array}\right|
u+1\right).\nonumber\\
\end{eqnarray*}
\end{prop}
By using this proposition, we have
\begin{eqnarray*}
&&I_{\hat{\Phi}}(v_1,v_2,\cdots,v_k|u_1,u_2)\nonumber\\
&=&(-1)^k \sqrt{\frac{[n+2k]}{[n]}}
(n+k,n+2k)_k
I_{{\Lambda}}
(v_1,v_2,\cdots,v_k|u_1,-u_2)
\prod_{j=1}^k
\frac{[u_2-v_j-\frac{k}{2}]}
{[u_2-v_j+\frac{k}{2}]}.\nonumber\\
\end{eqnarray*}
Therefore (\ref{com5}) and (\ref{com6})
reduce to the same identity,
which we have checked for $k=1,2,3$.
This supports our conjecture \mref{tail1} for the explicit form of the
tail operator.

Let us check a commutation relation of the tail and vertex
operators which does {\it not} reduce to one for just vertex operators
in
the above way.
Let us consider the commutation relation
for $s \geq c$,
\begin{eqnarray*}
&&
{\Lambda}_{-2s}(u_2)
\Phi_k(u_1)\\
&=&
\sum_{t=0}^k
L^{(k)}\left(\left.\begin{array}{cc}
P+h&P+h+k-2t\\
P+h+2s&P+h-k+2s
\end{array}\right|u_1-u_2\right)
\Phi_{-k+2t}(u_1)
{\Lambda}_{2k-2s-2t}(u_2).
\nonumber
\end{eqnarray*}
Taking the residue at $u_1=-u_2-k$,
a necessary condition becomes the
following theta identity
\begin{eqnarray*}
&&\sum_{t=0}^k
[n+k-s]_t[s]_{k-t}
[2k+n-1-s-t]_{k-t}[s+t-1]_t\nonumber\\
&\times&
(-1)^{s+t-k}\sqrt{\frac{[n-2s-2t+2k]}{[n]}}
\frac{1}{
(n+k-2s,n+2k-2s-2t)_k}=0.
\end{eqnarray*}
\end{appendix}

%\newpage
\small
\baselineskip=12pt
%\bibliography{vfusion}   % Bibliography

\begin{thebibliography}{10}

\bibitem{Bax72a}
R.J. Baxter.
\newblock {Partition Function of the Eight-Vertex Lattice Model}.
\newblock {\em {Annals of Physics}}, 70:193--228, 1972.

\bibitem{Bax73aI}
R.J. Baxter.
\newblock {Eight-Vertex Model in Lattice Statistics and One-Dimenisional
  Anisotropic Heisenberg Chain. 1. Some Fundamental Eigenvectors}.
\newblock {\em Annals of Physics}, 76:1--24, 1973.

\bibitem{Bax73aII}
R.J. Baxter.
\newblock {Eight-Vertex Model in Lattice Statistics and One-Dimenisional
  Anisotropic Heisenberg Chain. 11. Equivalence to a Generalized Ice-type
  Model}.
\newblock {\em Annals of Physics}, 76:25--47, 1973.

\bibitem{Bax73aIII}
R.J. Baxter.
\newblock {Eight-Vertex Model in Lattice Statistics and One-Dimenisional
  Anisotropic Heisenberg Chain. 111. Eigenvectors of the Transfer Matrix and
  Hamitonian}.
\newblock {\em Annals of Physics}, 76:48--71, 1973.

\bibitem{ABF}
G.~E. Andrews, R.~J. Baxter, and P.~J. Forrester.
\newblock {Eight-vertex SOS Model and Generalized Rogers-Ramanujan-Type
  Identities}.
\newblock {\em J. Stat. Phys.}, 35:193--266, 1984.

\bibitem{Huse84}
D.A. Huse.
\newblock {Exact exponents for infinitely many new multicritical points}.
\newblock {\em Phys. Rev.}, B30:3908--3915, 1984.

\bibitem{KRYS81}
P.~P. Kulish, N.~Yu. Reshetikhin, and E.K. Sklyanin.
\newblock {Yang-Baxter Equation and Representation Theory: I}.
\newblock {\em Lett. Math. Phys}, 5:393--403, 1981.

\bibitem{KS82}
P.~P. Kulish and E.K. Sklyanin.
\newblock {Solutions of the Yang-Baxter Equation}.
\newblock {\em J. Soviet Math.}, 19:1596--1620, 1982.

\bibitem{DJMO86}
E.~Date, M.~Jimbo, T.~Miwa, and M.~Okado.
\newblock {Fusion of the Eight Vertex SOS Model}.
\newblock {\em {Lett. Math. Phys.}}, 12:209--215, 1986.

\bibitem{DJKMO87}
E.~Date, M.~Jimbo, A.~Kuniba, T.~Miwa, and M.~Okado.
\newblock {Exactly Solvable SOS Models}.
\newblock {\em Nucl. Phys.}, B290 [FS20]:231--273, 1987.

\bibitem{DJKMO88}
E.~Date, M.~Jimbo, A.~Kuniba, T.~Miwa, and M.~Okado.
\newblock {Exactly Solvable SOS Models II}.
\newblock {\em Adv. Studies in Pure Math.}, 16:17--122, 1988.

\bibitem{Daval}
B.~Davies, O.~Foda, M.~Jimbo, T.~Miwa, and A.~Nakayashiki.
\newblock {Diagonalization of the XXZ Hamiltonian by Vertex Operators}.
\newblock {\em Comm. Math. Phys.}, 151:89--153, 1993.

\bibitem{JM}
M.~Jimbo and T.~Miwa.
\newblock {\em Algebraic Analysis of Solvable Lattice Models}.
\newblock CBMS Regional Conference Series in Mathematics, vol. 85. Amer. Math.
  Soc., 1994.

\bibitem{FIJKMYa}
O.~Foda, K.~Iohara, M.~Jimbo, R.~Kedem, T.~Miwa, and H.~Yan.
\newblock An elliptic quantum algebra for $\widehat{sl}_2$.
\newblock {\em Lett. Math. Phys.}, 32:259--268, 1994.

\bibitem{FIJKMYb}
O.~Foda, K.~Iohara, M.~Jimbo, R.~Kedem, T.~Miwa, and H.~Yan.
\newblock Notes on highest weight modules of the elliptic algebra
  ${A}_{q,p}(\widehat{sl}_2)$.
\newblock {\em Prog. Theor. Phys. Suppl.}, 118:1--34, 1995.

\bibitem{JKOS1}
M.~Jimbo, H.~Konno, S.~Odake, and J.~Shiraishi.
\newblock {Quasi-Hopf twistors for elliptic quantum groups}.
\newblock {\em {Transformation Groups}}, 4:303--327, 1999.

\bibitem{JMN}
M.~Jimbo, T.~Miwa, and A.~Nakayashiki.
\newblock {Difference Equations for the Correlation Functions of the
  Eight-Vertex Model}.
\newblock {\em J. Phys.}, A26:2199--2209, 1993.

\bibitem{JKKMW}
M.~Jimbo, R.~Kedem, H.~Konno, T.~Miwa, and R.A. Weston.
\newblock Difference equations in spin chains with a boundary.
\newblock {\em Nucl. Phys.}, B448:429--456, 1995.

\bibitem{JMO93}
M.~Jimbo, T.~Miwa, and Y.~Ohta.
\newblock {Structure of the Space of States in RSOS Models}.
\newblock {\em {Int. J. Mod. Phys.}}, A8:1457--1477, 1993.

\bibitem{LP96}
S.~Lukyanov and Y.~Pugai.
\newblock {Multi-point Local Height Probabilities in the Integrable RSOS
  Models}.
\newblock {\em Nucl. Phys.}, B473:631--658, 1996.

\bibitem{LaP98}
M.~Lashkevich and Y.~Pugai.
\newblock {Free Field Construction for Correlation Functions of the
  Eight-Vertex Model}.
\newblock {\em Nucl. Phys.}, B516:623--651, 1998.

\bibitem{La02}
M.~Lashkevich.
\newblock {Free field construction for the eight-vertex model: representation
  for form factors}.
\newblock {\em Nucl. Phys.}, B621:587--621, 2002.

\bibitem{Ko98}
H.~Konno.
\newblock {An Elliptic Algebra $U_{q,p}(\widehat{sl}_2)$ and the Fusion RSOS
  Models}.
\newblock {\em Comm. Math. Phys.}, 195:373--403, 1998.

\bibitem{JKOS}
M.~Jimbo, H.~Konno, S.~Odake, and J.~Shiraishi.
\newblock {Elliptic algebra $U\sb {q,p}(\widehat{ \mathfrak{sl}}\sb 2)$:
  Drinfeld currents and vertex operators}.
\newblock {\em {Comm. Math. Phys. }}, 199:605--647, 1999.

\bibitem{KKW04b}
T.~Kojima, H.~Konno, and R.~Weston.
\newblock {The Vertex-Face Ccorrespondence and Correlation Functions of the
  Fusion Eight-Vertex Models II: The 21 Vertex Model}.
\newblock To be published.

\bibitem{Konno04}
H.~Konno.
\newblock Fusion of {B}axter's elliptic ${R}$-matrix and the vertex-face
  correspondence.
\newblock {\em RIMS Koukyuroku}, page to appear, 2005. math.QA/0503726.

\bibitem{Fodal93}
O.~Foda, M.~Jimbo, T.~Miwa, K.~Miki, and A.~Nakayashiki.
\newblock {Vertex Operators of Solvable Lattice Models}.
\newblock {\em J. Math. Phys.}, 35:13--46, 1994.

\bibitem{DateLMP89}
E.~Date, M.~Jimbo, A.~Kuniba, T.~Miwa, and M.~Okado.
\newblock One-dimensional configuration sums in vertex models and affine {L}ie
  algebra characters.
\newblock {\em Lett. Math. Phys.}, 17:69--77, 1989.

\bibitem{KaPe}
V.~G. Kac and D.~H. Peterson.
\newblock Infinite dimensional {L}ie algebras, theta-functions and modular
  forms.
\newblock {\em Adv. in Math.}, 53:125--264, 1984.

\bibitem{DJKMO89}
E.~Date, M.~Jimbo, A.~Kuniba, T.~Miwa, and M.~Okado.
\newblock {Paths, Maya Diagrams and Representations of $\slth(r,\C)$}.
\newblock {\em Adv. Studies in Pure Math.}, 19:149--191, 1989.

\bibitem{DJMO86b}
E.~Date, M.~Jimbo, T.~Miwa, and M.~Okado.
\newblock Automorphic properties of local height probabilities for integrable
  solid-on-solid models.
\newblock {\em Phys.Rev.B}, 35:2105--2107, 1986.

\bibitem{idzal93}
M.~Idzumi, T.~Tokihiro, K.~Iohara, M.~Jimbo, T.~Miwa, and T.~Nakashima.
\newblock {Quantum Affine Symmetry in Vertex Models}.
\newblock {\em Int. J. Mod. Phys.}, A8:1479--1511, 1993.

\bibitem{Idz94}
M.~Idzumi.
\newblock {Level 2 Irreducible Representations of $U_q(\widehat{sl_2})$, Vertex
  Operators, and their Correlation Functions}.
\newblock {\em Int. J. Mod. Phys.}, A9:4449--4484, 1994.

\bibitem{BoWe94b}
A.~H. Bougourzi and R.~A. Weston.
\newblock N-point correlation functions of the spin-1 {X}{X}{Z} model.
\newblock {\em Nucl. Phys}, B417:439--462, 1994.

\bibitem{Mat94}
A.~Matsuo.
\newblock {A q-Deformation of Wakimoto Modules, Primary Fields and Screening
  Operators}.
\newblock {\em Comm. Math. Phys.}, 160:33--48, 1994.

\bibitem{GeQiu}
D.~Gepner and Z.~Qiu.
\newblock {Modular invariant partition functions for parafermionic field
  theories}.
\newblock {\em Nucl. Phys.}, B285:423--, 1987.

\bibitem{LePr85}
J.~Lepowsky and M.~Primc.
\newblock {Structure of the standard modules for the affine Lie algebra $A\sp
  {[1]}\sb 1$}.
\newblock {\em Contemporary Mathematics}, 46, 1985.

\bibitem{Jing}
N.~H. Jing.
\newblock Higher level representations of the quantum affine algebra
  ${U}_q(\widehat{sl}_2)$.
\newblock {\em J. Algebra}, 182:448--468, 1996.

\bibitem{FrRe}
E.V. Frenkel and N.~Yu Reshetikhin.
\newblock Deformations of {$W$}-algebras associated to simple {L}ie algebras.
\newblock {\em Comm. Math. Phys.}, 197:1--32, 1998.

\bibitem{Konno05}
H.~Konno.
\newblock The vertex-face correspondence and the elliptic $6j$-symbols,
  math.QA/0503725.
\newblock To be published in Lett.Math. Phys.

\end{thebibliography}
%%%%%%%%%%%%%%%%%%%%%%%%%%%%%%%%%%  

\end{document}